\documentclass[12pt, onecolumn, draftclsnofoot]{IEEEtran}

\newcommand{\paren}[1]{\left(#1\right)}
\newcommand{\sqparen}[1]{\left[#1\right]}
\newcommand{\brparen}[1]{\left\{#1\right\}}
\newcommand{\field}[1]{\ensuremath{\mathbb{#1}}}
\newcommand{\R}{\ensuremath{\field{R}}} 
\newcommand{\I}[1]{\ensuremath{\mathsf{1}_{\left\{#1\right\}}}} 
\newcommand{\ra}{\ensuremath{\rightarrow}} 
\newcommand{\PR}[1]{\ensuremath{\mathsf{Pr}\left\{#1\right\}}} 
\newcommand{\PRP}[1]{\ensuremath{\mathsf{Pr}\left(#1\right)}} 
\newcommand{\ES}[1]{\ensuremath{\mathsf{E}\left[#1 \right]}} 
\newcommand{\V}[1]{\ensuremath{\mathsf{Var}\left[#1 \right]}} 
\renewcommand{\Re}{\ensuremath{\R}} %
\newcommand{\e}[1]{\ensuremath{{\rm e}^{#1}}} 
\newcommand{\diff}{\ensuremath{{\rm d}}} 
\newcommand{\disteq}{\ensuremath{\stackrel{\rm d}{=}}}

\newcommand{\sinr}{\ensuremath{{\rm SINR}}}
\newcommand{\snr}{\ensuremath{{\rm SNR}}}
\newcommand{\pg}{\ensuremath{{\rm PG}}}

\newcommand{\Lap}[1]{\ensuremath{\mathcal{L}_{#1}}}

\newcommand{\Ceil}[1]{\ensuremath{\left\lceil #1\right\rceil}}
\newcommand{\distconv}{\ensuremath{\stackrel{\rm d}{\ra}}}

\newcommand{\BO}[1]{\ensuremath{O\paren{#1}}}
\newcommand{\LO}[1]{\ensuremath{o\paren{#1}}}
\newcommand{\TO}[1]{\ensuremath{\Theta\paren{#1}}}
\newcommand{\OMO}[1]{\ensuremath{\Omega\paren{#1}}}
\renewcommand{\vec}[1]{\ensuremath{\boldsymbol{#1}}}

\newtheorem{theorem}{Theorem}
\newtheorem{lemma}{Lemma}

\newtheorem{corollary}{Corollary}

\usepackage{amsmath, latexsym, amsfonts, amssymb}
\usepackage[mathscr]{eucal}
\usepackage{graphicx}
\usepackage{array, tabularx, colortbl, color, threeparttable}
\usepackage[table]{xcolor}
\usepackage{psfrag, pstricks}
\usepackage{multirow}
\usepackage{booktabs}

\begin{document}
\title{Gaussian Approximation for the Wireless Multi-access Interference Distribution and Its Applications}

\author{Hazer Inaltekin, \IEEEmembership{Member, IEEE} 
\thanks{This research was supported in part by the European Union Research Executive Agency Marie Curie FP7-Reintegration-Grants under Grant PCIG10-GA-2011-303713 and in part by the Australian Research Council under Grant DP-11-0102729.}
\thanks{H. Inaltekin is with the Department of Electrical and Electronics Engineering, Antalya International University, Dosemealti, Antalya, Turkey. E-mail: hazeri@antalya.edu.tr. Phone: +902423216074. Fax: +902423216072.}
}
\date{}
\maketitle
\thispagestyle{empty}
\begin{abstract} 
\textbf{
This paper investigates the problem of Gaussian approximation for the wireless multi-access interference distribution in large spatial wireless networks. First, a principled methodology is presented to establish rates of convergence of the multi-access interference distribution to a Gaussian distribution for general bounded and power-law decaying path-loss functions.  The model is general enough to also include various random wireless channel dynamics such as fading and shadowing arising from multipath propagation and obstacles existing in the communication environment.  It is shown that the wireless multi-access interference distribution converges to the Gaussian distribution with the same mean and variance at a rate $\frac{1}{\sqrt{\lambda}}$, where $\lambda>0$ is a parameter controlling the intensity of the planar (possibly non-stationary) Poisson point process generating node locations.  An explicit expression for the scaling coefficient is obtained as a function of fading statistics and the path-loss function.  Second, an extensive numerical and simulation study is performed to illustrate the accuracy of the derived Gaussian approximation bounds.  A good statistical fit between the interference distribution and its Gaussian approximation is observed for moderate to high values of $\lambda$.  Finally, applications of these approximation results to upper and lower bound the outage capacity and ergodic sum capacity for spatial wireless networks are illustrated.  The derived performance bounds on these capacity metrics track the network performance within one nats per second per hertz.}

\textcolor{black}{{\bf EDICS}: WIN-PHLY, WIN-INFO}
\end{abstract}
\setcounter{page}{0}
\section{Introduction}
\subsection{Background and Motivation}
Wireless communication technologies have been evolved at an unprecedented pace over the last decade.  This has, in turn, given rise to the birth of many next generation wireless systems including 4G networks, femtocells and cognitive radio networks.  The main design philosophy underpinning most of such emerging classes of wireless systems, if not all, is the more proficient utilization of wireless spectrum than ever before to accommodate high volume of data traffic from increasingly overcrowded user populations within a given frequency band as well as to support increasingly more data-rate-intense multimedia applications over the wireless links.  
Interference mitigation and characterization are the primary design challenges to overcome for achieving this goal and for meeting the high target data rates ({\em e.g.,} $1$ Gbit/s for low mobility and $100$ Mbit/s for high mobility environments for 4G networks \cite{ITU08}) set for future wireless networks consisting of many interfering links.  

However, even as a first step to characterize network performance and relevant performance metrics, computation of the exact wireless multi-access interference (WMAI) distributions turns out to be a mathematically intractable problem in most practical scenarios.  This is mainly because the level of WMAI at a receiver node in a wireless network is a random function of  the overall network geometry as well as numerous other wireless channel dynamics such as path-loss, fading, shadowing and so on, {\em e.g.,} see \cite{JSAC09-Tutorial}-\cite{WPS09}.  This motivates the search for tight bounds on the WMAI distributions that can accurately track the interference behavior \cite{WA06}-\cite{AY10}.    

This paper focuses on the statistical characterization of WMAI for {\em spatial} wireless networks, and establishes tight Gaussian approximation bounds for the WMAI distributions.  To this end, the underlying spatial stochastic process determining transmitter locations is assumed to be {\em Poisson} but not necessarily stationary.  The signal power attenuation due to path-loss is modeled by means of a general bounded and power-law decaying path-loss function.  Other random wireless channel dynamics such as fading affecting the received signal power are also accounted for in the employed signal propagation model. 

Our main contribution is the derivation of the rate of convergence of the WMAI distributions to the Gaussian distribution with the same mean and variance, which is formally stated in Theorem \ref{Thm: Rates of Convergence}.  Briefly, this rate is equal to $\frac{c(x)}{\sqrt{\lambda}}$, where $\lambda$ is a modeling parameter enabling us to control the ``intensity" of the planar Poisson point process (PPP) generating transmitter locations ({\em i.e.,} see Section \ref{Section: Network Model} for details), and $c: \R \mapsto \R_+$ is a positive function which depends on the shape of the path-loss function and the point $x \in \Re$ at which we want to estimate the value of the (centered and normalized) WMAI distribution function.  $c(x)$ approaches zero for large absolute values of $x$ at a rate $|x|^{-3}$.  This behavior makes the derived bounds on the tails of the WMAI distributions tight for any given value of $\lambda$.  Moreover, the supremum of $c(x)$ over $x$ is a small constant, which allows us to obtain {\em uniform} rates of convergence as a function of $\lambda$. 

An extensive simulation and numerical study is performed to illustrate the accuracy of derived Gaussian approximation results estimating WMAI distributions for specific instances of both stationary ({\em i.e.,} see Section \ref{Section: Approximation Bounds for Stationary PPPs}) and non-stationary ({\em i.e.,} see Appendix \ref{Appendix: Normal Approximation Non-stationary PPP}) PPPs.  It is observed that the Gaussian distribution  with the same mean and variance can accurately track the behavior of the WMAI distribution even for small values of $\lambda$.  As predicted by our bounds, the approximation accuracy increases further as $\lambda$ increases, and the WMAI distribution becomes {\em almost} undistinguishable from its Gaussian approximation when $\lambda$ is around $10$.  On the other hand, random fading effects have an adverse impact on this Gaussian behavior, and the bounds given in Theorem \ref{Thm: Rates of Convergence} become looser when the fading distributions have larger dynamic ranges.      

The utility of the derived Gaussian approximation bounds on the WMAI distributions is also illustrated by characterizing wireless network performance in tangible communication scenarios.  In particular, tight upper and lower bounds on the outage capacity for a given victim link in a spatial wireless network and those on the sum capacity for spatial wireless multiple-access networks are obtained.  The approximation accuracy of the derived performance bounds on these network capacity measures lies within {\em one} nats per second per hertz for moderate to high values of $\lambda$.  For high values of $\lambda$, our bounds become very tight, and they {\em almost} coincide with the outage capacity ({\em i.e.,} see Subsection \ref{Subsection: Outage Capacity}) and the sum capacity ({\em i.e.,} see Subsection \ref{Subsection: Network Capacity}).                              

\subsection{Related Work and Paper Organization}
Wireless networks are often interference limited due to the broadcast nature of the wireless medium.  Hence, the statistical characterization of interference for large spatial wireless networks has been a key area of research for more than a decade.  Many of the latest developments in the field are summarized in the recent review articles \cite{JSAC09-Tutorial}-\cite{WPS09}.  

Historically,  the early efforts for characterizing the structure of WMAI in wireless networks by using stochastic geometry can be traced back to as early as 1978 \cite{Musa78}.  Sousa et al. applied similar techniques, {\em e.g.,} \cite{Sousa90} and \cite{Sousa92}, in the 1990s to assess the performance of spread spectrum wireless networks as well as to find optimum transmission ranges in these networks.  Subsequently, several approximation techniques appeared in the field to approximate the WMAI distributions, and then to use these approximation results for the network performance analysis, {\em e.g.,} \cite{WA06}, \cite{IH98} and \cite{CH01}.  However, except for one special case \cite{Sousa90}, {\em i.e.,} the case in which the signal power decays according to the unbounded power-law decaying path-loss function $t^{-4}$, there are still no closed form expressions available for the WMAI distributions.    

Therefore, it becomes necessary to resort to numerical methods to calculate the WMAI distributions by modeling WMAI as a power-law shot-noise process and then inverting its computed characteristic function.  Even though algorithmic perspectives based on fast Fourier transforms \cite{Gubner96} to numerically compute power-law shot-noise distributions and densities are promising, they are of limited interest and importance in the context of wireless networking.  The numerical computation cannot provide closed form expressions revealing structural dependencies between the WMAI distributions and network design parameters to assess the wireless network performance under candidate/existing wireless communication technologies.  

Alternative approaches in the field include various approximation methods based on geometrical considerations and distribution approximation techniques \cite{WA06}-\cite{AY10}, LePage series \cite{IH98} or Edgeworth expansion \cite{CH01} to estimate the WMAI distributions.  The main objective is to obtain simple but insightful upper and lower bounds on the related performance quantities of interest in spatial wireless networks by utilizing such approximations.  In particular, our results in this paper are close to those in \cite{WA06}-\cite{AY10} in that the authors of these works have also obtained upper and lower bounds on the normalized WMAI distributions with provably small gaps between the computed bounds and the actual WMAI distribution.\footnote{To be more precise, the authors in \cite{WA06} focused on the distribution of the normalized inverted signal-to-interference-ratio, which is, in essence, the same as computing the WMAI distribution.}  

In contrast to \cite{WA06} as well as most other earlier works in the field such as \cite{Sousa90} and \cite{Sousa92}, one distinctive aspect of our analysis in this paper is that we work with general {\em bounded} and power-law decaying path-loss functions to calculate the WMAI distributions.  
The underlying motivation for working with bounded path-loss models is recent findings pointing out that the unrealistic singularity of the unbounded path-loss model at $0$ leads to unexpected deviations on the final computed WMAI distributions, {\em e.g.,} see \cite{Inaltekin09} and \cite{GH09}.  When compared with the results reported in \cite{InaltekinHanly10} and \cite{AY10}, we derive tighter Gaussian approximation bounds for more general spatial node distributions including non-stationary PPPs and for more general signal propagation models including fading and shadowing effects.  

Since interference in a wireless network is a specific instance of a shot-noise process, the results of this paper are also related to a more general body of work on shot-noise processes such as \cite{Lowen90} and \cite{HS85}.  The paper \cite{Lowen90} established many properties for power-law shot-noise processes on the line such as its moment generating functions, moments and cumulants.  For a very specific bounded and power-law decaying impulse response function driving the power-law shot-noise process, the convergence of the amplitude distribution of the power-law shot-noise process to the Gaussian form was briefly mentioned in \cite{Lowen90}, but without any formal proof for this convergence and without establishing rates of convergence.  

Our results are closely related to those reported in \cite{HS85} since the authors in \cite{HS85} also established the uniform rates of convergence for the amplitude distribution of the shot-noise process.  When compared with \cite{HS85}, our results are different than those of \cite{HS85} in three aspects.  Firstly, we give both uniform and non-uniform Berry-Esseen types of bounds on the WMAI distributions.  The non-uniform bounds allow us to tightly approximate the WMAI distributions for extreme interference values even for small finite values of $\lambda$.  Secondly, error terms appearing in our uniform Gaussian approximation formulas are eleven times smaller than those in \cite{HS85}.  This implies much tighter performance bounds for characterizing various performance measures for spatial wireless networks.  Finally, this paper introduces a principled and simpler methodology to establish Gaussian approximation results.  Hence, the distribution approximation techniques presented here are expected to find other potential applications in the analysis of emerging wireless network architectures. 

The remainder of the paper is organized as follows.  Section \ref{Section: Network Model} describes the network configuration along with our modeling assumptions.  In Section \ref{Section: WMAI Distribution}, we establish the main Gaussian approximation result for the WMAI distributions.  Section \ref{Section: Approximation Bounds for Stationary PPPs} illustrates the applications of the derived Gaussian approximation bounds for stationary PPPs, and presents an extensive numerical and simulation study to verify these bounds.  Similar analysis is also performed for a non-stationary PPP in Appendix \ref{Appendix: Normal Approximation Non-stationary PPP}.  Section \ref{Section: Applications} presents further applications of our Gaussian approximation bounds to characterize outage capacity and sum capacity for spatial wireless networks. Finally, Section \ref{Section: Conclusions} concludes the paper.  Most of our proofs are relegated to appendices for the sake of fluency of the paper.         

\section{Network Model} \label{Section: Network Model}

We consider a planar network in which transmitters are distributed according to a planar PPP with {\em mean measure} $\Lambda$ (alternatively called: intensity 
measure), denoted by $\Phi_{\Lambda}$, over $\R^2$.  Here, for any (Borel) subset $\mathcal{A}$ of $\R^2$, $\Lambda\paren{\mathcal{A}}$ gives us the {\em average} number of transmitters lying in $\mathcal{A}$.  We will assume that $\Lambda\paren{\mathcal{A}}$ is {\em locally finite}, {\em i.e.,} $\Lambda\paren{\mathcal{A}} < \infty$ for bounded subsets $\mathcal{A}$ of $\R^2$, and $\Lambda\paren{\R^2} = \infty$, {\em i.e.,} there is an infinite population of transmitters scattered all around in $\R^2$.  The location of the $k$th transmitter is represented by $\vec{X}_k$.  We will often represent the transmitter location process $\Phi_\Lambda$ as a discrete sum of Dirac measures as $\Phi_\Lambda = \sum_{k \geq 1} \delta_{\vec{X}_k}$, where $\delta_{\vec{X}_k}\paren{\mathcal{A}} = 1$ if $\vec{X}_k \in \mathcal{A} \subseteq \Re^2$, and zero otherwise.  $\vec{X}_k$'s can be interpreted as {\em points} (or, atoms) of $\Phi_\Lambda$, and therefore, we also use the notation $\vec{X}_k \in \Phi_\Lambda$ to symbolize this interpretation.  For a point $\vec{x} \in \R^2$, $x^{(i)}$, $i = 1, 2$, represents the $i$th component of $\vec{x}$.  $\mathcal{B}\paren{\vec{x}, r}$ represents the (planar) ball centered at $\vec{x}$ with radius $r$.       

We consider the case in which all transmitters transmit with the same power $P$, {\em i.e.,} the case of a non-power controlled wireless network.  For the signal power attenuation in the wireless medium, we consider a {\em bounded} monotone non-increasing path-loss function $G: [0, \infty) \mapsto [0, \infty)$, which asymptotically decays to zero at least as fast as $t^{-\alpha}$ for some path-loss exponent $\alpha > 2$.  In addition to the signal power attenuation due to path-loss, it is assumed that transmitted signals are also corrupted by fading.  The random (power) fading coefficient at transmitter $k$ is given by $H_k$, and is assumed to be independent of $\Phi_\Lambda$.\footnote{For simplicity, we only assign a {\em single} fading coefficient to each transmitter.  In reality, it is expected that the channels between a transmitter and all potential receivers (intended or unintended) experience different and possibly independent fading processes.  Our simplified notation does not cause any ambiguity here since we focus on the total interference power at a given arbitrary position in $\R^2$ in the remainder of the paper.}  The received signal power at a distance $t_k$ is given by $P H_k G\paren{t_k}$ for transmitter $k$.  $H_k$'s are independent and identically distributed (i.i.d.) with a common density $q(h)$, $h \geq 0$, and have finite first, second and third order moments, {\em i.e.,} $\int_0^\infty h^k q(h) dh < \infty$, $k = 1, 2 \mbox{ and } 3$.  The first, second and third order moments of the fading coefficients are denoted by $m_H$, $m_{H^2}$ and $m_{H^3}$, respectively.  We note that the employed signal propagation model is general enough that $H_k$'s could also be thought to incorporate {\em shadow fading} effects due to blocking of signals by large obstacles existing in the communication environment although we do not model such random factors explicitly and separately in this paper.  We also note that our model is general enough to include the case where transmitted signal powers are also random due to power control, {\em i.e.,} we just need to scale the fading process with transmitted signal powers in this case.  As it is common in most of the earlier works, we will assume a contention based medium access control (MAC) layer such as ALOHA mediating node transmissions, and giving rise to the observed distribution of active transmitters over $\R^2$.     

We place a test receiver node at an arbitrary point $\vec{X}_o = \paren{X_o^{(1)}, X_o^{(2)}} \in \R^2$, and consider signals coming from all other transmitters, whose locations are given by $\Phi_\Lambda$, as interference to this test receiver node.  Without loss of generality, we focus on the distribution of WMAI seen by such a test receiver node placed at arbitrary $\vec{X}_o$ for the rest of the paper.  If $\Phi_\Lambda$ is {\em stationary}, the interference statistics seen from any other point in $\Re^2$ will be the same.  

The level of WMAI at $\vec{X}_o$ depends on the distances between the points of $\Phi_\Lambda$ and $\vec{X}_o$.  Hence, the transformed process $\sum_{k \geq 1} \delta_{T\paren{\vec{X}_k}}$, whose points lie on the positive real line, is of particular importance to derive interference statistics at $\vec{X}_o$, where $T: \R^2 \mapsto \R$ is given by $T\paren{\vec{x}} = \|\vec{x} - \vec{X}_o\|_2 = \sqrt{\paren{x^{(1)} - X_o^{(1)}}^2 + \paren{x^{(2)} - X_o^{(2)}}^2}$.  The mean measure of the transformed process is equal to $\Lambda \circ T^{-1}$, where $T^{-1}\paren{\mathcal{A}} = \brparen{\vec{x} \in \R^2: T\paren{\vec{x}} \in \mathcal{A}}$ for all $\mathcal{A} \subseteq \R$, and $\Lambda \circ T^{-1}\paren{\mathcal{A}}$ has the same interpretation above, {\em i.e.,} it gives us the average number of points of the transformed process in $\mathcal{A}$ (see the Mapping Theorem on page 18 in \cite{Kingman93}).  We assume that $\Lambda \circ T^{-1}$ has a density in the form 
$
\Lambda \circ T^{-1}\paren{\mathcal{A}} = \lambda \int_{\mathcal{A}} p(t) \diff t \nonumber
$     
such that $p(t) = \BO{t^{\alpha - 1 - \epsilon}}$ as $t \ra \infty$ for some $\epsilon > 0$.  This assumption on the rate of growth of $p(t)$ is necessary to ensure the finiteness of WMAI at $\vec{X}_o$.  Here, $\lambda$ is a modeling parameter, which can be interpreted as the {\em transmitter intensity parameter}, that will enable us to control the average number of transmitters lying in $\mathcal{A}$ and interfering with the signal reception at the test receiver at $\vec{X}_o$.   

The level of WMAI at $\vec{X}_o$ is equal to    
$
I_\lambda = \sum_{k \geq 1} P H_k G\paren{T\paren{\vec{X}_k}}. \nonumber
$            
$I_\lambda$ is a random variable since transmitter locations and associated fading coefficients, {\em i.e.,} $\brparen{\paren{\vec{X}_k, H_k}, k \geq 1}$, are random variables.  Therefore, different node configurations and fading states result in different levels of interference at the test receiver.  In the next section, we will show that the distribution of $I_\lambda$ can be approximated by a Gaussian distribution.

{\em A Note about Notation:} We use boldface and calligraphic letters to denote vector quantities and sets, respectively. $| \cdot |$ notation is used to measure the magnitudes of scalar quantities, whereas $\| \cdot \|_2$ notation is used to measure the Euclidean norms of vector quantities.  As is standard, when we write $f(t) = \BO{g(t)}$, $f(t) = \OMO{g(t)}$ and $f(t) = \LO{g(t)}$ as $t \ra t_0$ for two positive functions $f(t)$ and $g(t)$, we mean $\limsup_{t \ra t_0} \frac{f(t)}{g(t)} < \infty$,  $\liminf_{t \ra t_0} \frac{f(t)}{g(t)} >0$ and $\lim_{t \ra t_0} \frac{f(t)}{g(t)} = 0$, respectively. $f(t)$ is said to be $\TO{g(t)}$ as $t \ra t_0$ if $f(t) = \BO{g(t)}$ and $f(t) = \OMO{g(t)}$ as $t \ra t_0$.  With a slight abuse of notation, we sometimes use $I_\lambda(P)$ notation to represent the level of WMAI at $\vec{X}_o$ when we need to put emphasis on transmission powers to explain some performance results in Section \ref{Section: Applications}.                   

\section{WMAI Distribution and Rates of Convergence to the Gaussian Distribution} \label{Section: WMAI Distribution}
This section presents calculations for approximating the WMAI distributions as a Gaussian distribution, and establishes the rates of convergence for this approximation as $\lambda$ grows large.  By using Laplace functionals of Poisson processes (see \cite{Kingman93} for details), we have the following Laplace transform for $I_\lambda$:
\begin{eqnarray}
\Lap{I_\lambda}(s) = \ES{\e{-s I_\lambda}} = \exp\paren{-\lambda \int_0^\infty \int_0^\infty \paren{1 - \e{-s P h G(t)}} q(h) p(t) \diff t \diff h}, \nonumber   
\end{eqnarray} 
where $s \geq 0$.  The next lemma shows that $I_\lambda$ has a non-degenerate distribution. 
\begin{lemma} \label{Lemma: Non-degenerate Distribution}
For all $s \geq 0$, $\int_0^\infty \int_0^\infty \paren{1 - \e{-s P h G(t)}} q(h) p(t) \diff h \diff t < \infty$. 
\end{lemma} 
\begin{IEEEproof}
Please see Appendix \ref{Appendix: Non-degenerate Distribution}. 
\end{IEEEproof}

We will need some auxiliary results to prove the main approximation result of the paper.  The proofs of these auxiliary results are relegated to the appendices at the end for the sake of fluency of the paper.  The next lemma shows that the distribution of $I_\lambda$ can be approximated as a limit distribution of a sequence of random variables $I_n$, {\em i.e.,} $I_n \distconv I_\lambda$ as $n \ra \infty$. 

\begin{lemma} \label{Lemma: Distribution Convergence}
For each $n$, let $U_{1, n}, \ldots, U_{\Ceil{\Lambda_n}, n}$ be a sequence of i.i.d. random variables with a common density $f(t) = \frac{\lambda p(t)}{\Lambda_n} \I{0 \leq t \leq n}$,  where $\Lambda_n = \lambda \int_0^n p(t) \diff t$ and $\Ceil{\cdot}$ is the smallest integer greater than or equal to its argument.  Let 
\begin{eqnarray}
I_n = \sum_{k=1}^{\Ceil{\Lambda_n}} P H_k G\paren{U_{k, n}}. \label{Eqn: I_n}
\end{eqnarray}
Then, $I_n$ converges in distribution to $I_\lambda$, which is shown as $I_n \distconv I_\lambda$, as $n \ra \infty$.      
\end{lemma} 
\begin{IEEEproof}
Please see Appendix \ref{Appendix: Distribution Convergence}.  
\end{IEEEproof} 

The next lemma shows that the mean value and variance of $I_\lambda$ can also be approximated by the mean value and variance of $I_n$.
\begin{lemma} \label{Lemma: Mean Variance Convergence}
Let $I_n$ be defined as in \eqref{Eqn: I_n}.  Then,  
\begin{eqnarray}
\lim_{n \ra \infty} \ES{I_n} = \ES{I_\lambda} \nonumber
\end{eqnarray}
and
\begin{eqnarray}
\lim_{n \ra \infty} \V{I_n} = \V{I_\lambda} \nonumber
\end{eqnarray}
\end{lemma} 
\begin{IEEEproof}
Please see Appendix \ref{Appendix: Mean Variance Convergence}.
\end{IEEEproof}

\begin{lemma} \label{Thm: Berry Esseen}
Let $\xi_1, \ldots, \xi_j$ be a sequence of independent and real-valued random variables such that each of which has zero mean and $\sum_{i=1}^j\ES{\xi_i^2} = 1$.  Let $\chi = \sum_{i=1}^j \ES{\left| \xi_i \right|^3}$.  Then,  
\begin{eqnarray} 
\left| \PR{\sum_{i=1}^j \xi_i \leq x} - \Psi(x) \right| \leq \chi \min\paren{0.4785, \frac{31.935}{1+|x|^3}} \nonumber
\end{eqnarray} 
for all $x \in \Re$.
\end{lemma}
\begin{IEEEproof}
The main theorem stated in page 1 in \cite{Tyurin10} implies that $\left| \PR{\sum_{i=1}^j \xi_i \leq x} - \Psi(x) \right| \leq \chi 0.4785$. Theorem 1 in \cite{Paditz89} after some simplifications implies that $\left| \PR{\sum_{i=1}^j \xi_i \leq x} - \Psi(x) \right| \leq \chi \frac{31.935}{1+|x|^3}$. Combining these two bounds completes the proof.    
\end{IEEEproof}

By using the above assisting lemmas, the main approximation result of the paper is established in the next theorem.  

\begin{theorem} \label{Thm: Rates of Convergence}
For all $x \in \Re$,
\begin{eqnarray}
\left| \PR{\frac{I_\lambda - \ES{I_\lambda}}{\sqrt{\V{I_\lambda}}} \leq x} - \Psi(x)\right| \leq \frac{c(x)}{\sqrt{\lambda}},
\end{eqnarray}
where $\Psi(x) = \frac{1}{\sqrt{2 \pi}}\int_{-\infty}^x \e{-\frac{t^2}{2}} dt$, which is the standard normal cumulative distribution function (CDF), and $c(x) = \frac{m_{H^3}}{\paren{m_{H^2}}^\frac32}\frac{\int_0^\infty G^3(t) p(t) dt}{\paren{\int_0^\infty G^2(t) p(t) dt}^\frac32} \min\paren{0.4785, \frac{31.935}{1+|x|^3}}$. 
\end{theorem}
\begin{IEEEproof} We let $\xi_{k, n} = \frac{PH_kG\paren{U_{k, n}} - m_{k, n}}{\sigma_n}$ for $n \geq 1$ and $1 \leq k \leq \Ceil{\Lambda_n}$, where $U_{k, n}$ and $\Lambda_n$ are as defined above, $m_{k, n} = \ES{P H_k G\paren{U_{k, n}}}$ and $\sigma_n = \sqrt{\V{I_n}}$.  We note that $\ES{\xi_{k, n}} = 0$ and $\sum_{k=1}^{\Ceil{\Lambda_n}} \ES{\xi_{k, n}^2} = 1$.  Hence, the collection of random variables $\brparen{\xi_{k, n}}_{k=1}^{\Ceil{\Lambda_n}}$ is in the correct form to apply Lemma \ref{Thm: Berry Esseen}.  We need to calculate $\chi_n = \sum_{k=1}^{\Ceil{\Lambda_n}} \ES{\left| \xi_{k, n} \right|^3}$ to finish the proof.  We can upper bound this summation as follows.  
\begin{eqnarray}
\chi_n &\leq& \frac{1}{\sigma_n^3}\sum_{k=1}^{\Ceil{\Lambda_n}} \ES{\left| P H_k G\paren{U_{k, n}} + m_{k, n} \right|^3} \nonumber \\
&=& \frac{\Ceil{\Lambda_n}}{\sigma_n^3} \ES{\left| P H_1 G\paren{U_{1, n}} + m_{1, n} \right|^3} \nonumber \\
&=& \frac{\Ceil{\Lambda_n}}{\sigma_n^3} \ES{P^3 H_1^3 G^3\paren{U_{1, n}} + 3 P^2 H_1^2 G^2\paren{U_{1, n}} m_{1, n} + 3PH_1 G\paren{U_{1, n}} m_{1, n}^2 + m_{1, n}^3} \nonumber
\end{eqnarray}
\begin{eqnarray}
&=&  \frac{\Ceil{\Lambda_n}P^3 m_{H^3} \lambda}{\Lambda_n \sigma_n^3} \int_0^n G^3(t)p(t)\diff t + \frac{3\Ceil{\Lambda_n}P^2m_{H^2} m_{1,n} \lambda}{\Lambda_n \sigma_n^3}\int_0^n G^2(t) p(t) \diff t \nonumber \\
& & + \frac{3 \Ceil{\Lambda_n} P m_H m_{1, n}^2 \lambda}{\Lambda_n \sigma_n^3} \int_0^n G(t)p(t)\diff t + \frac{\Ceil{\Lambda_n} m_{1, n}^3}{\sigma_n^3}. \nonumber
\end{eqnarray}
Observing $m_{1, n} = \LO{1}$ and $m_{1, n}^3 \Ceil{\Lambda_n} = \LO{1}$ as $n \ra \infty$, {\em i.e.,} see Appendix \ref{Appendix: Mean Variance Convergence}, and the convergence of $\sigma_n^2$ to $\V{I_\lambda}$ as $n \ra \infty$, {\em i.e.,} see Lemma \ref{Lemma: Mean Variance Convergence}, we have
\begin{eqnarray}
\limsup_{n \ra \infty} \chi_n \leq \frac{P^3 m_{H^3} \lambda}{\paren{\V{I_\lambda}}^\frac32} \int_0^\infty G^3(t) p(t) \diff t. \label{Eqn: Limsup Bound 1}
\end{eqnarray} 
Substituting the expression given for $\V{I_\lambda}$ in \eqref{Eqn: Interference Variance} in Appendix \ref{Appendix: Mean Variance Convergence}, we can rewrite \eqref{Eqn: Limsup Bound 1} as 
\begin{eqnarray}
\limsup_{n \ra \infty} \chi_n \leq \frac{m_{H^3}}{\sqrt{\lambda} \paren{m_{H^2}}^\frac32} \frac{\int_0^\infty G^3(t) p(t) \diff t}{\paren{\int_0^\infty G^2(t) p(t) \diff t}^\frac32}.  \label{Eqn: Limsup Bound 2}
\end{eqnarray}   
By using Theorem \ref{Thm: Berry Esseen}, we have 
\begin{eqnarray}
\left| \PR{\sum_{k=1}^{\Ceil{\Lambda_n}} \xi_{k, n} \leq x} - \Psi(x) \right| \leq \chi_n \min\paren{0.4785, \frac{31.935}{1+|x|^3}} \label{Eqn: Berry Esseen Bound}
\end{eqnarray} 
for all $n \geq 1$ and $x \in \R$.  Lemmas \ref{Lemma: Distribution Convergence} and \ref{Lemma: Mean Variance Convergence} imply that
\begin{eqnarray}
\sum_{k=1}^{\Ceil{\Lambda_n}} \xi_{k,n}  \distconv \frac{I_\lambda - \ES{I_\lambda}}{\sqrt{\V{I_\lambda}}} \mbox{ as } n \ra \infty.  \nonumber
\end{eqnarray}  

Hence, using \eqref{Eqn: Limsup Bound 2} and taking the $\limsup$ of both sides in \eqref{Eqn: Berry Esseen Bound}, we have  
\begin{eqnarray}
\lefteqn{\limsup_{n \ra \infty} \left| \PR{\sum_{k=1}^{\Ceil{\Lambda_n}} \xi_{k, n} \leq x} - \Psi(x) \right|} \hspace{10cm} \nonumber \\
\lefteqn{= \left| \PR{\frac{I_\lambda - \ES{I_\lambda}}{\sqrt{\V{I_\lambda}}} \leq x} - \Psi(x) \right|} \hspace{7cm} \nonumber \\
\lefteqn{\leq \frac{1}{\sqrt{\lambda}} \frac{m_{H^3}}{\paren{m_{H^2}}^\frac32} \frac{\int_0^\infty G^3(t)p(t)\diff t}{\paren{\int_0^\infty G^2(t)p(t)\diff t}^\frac32} \min\paren{0.4785, \frac{31.935}{1+|x|^3}},} \hspace{7cm} \nonumber
\end{eqnarray}
which completes the proof. 
\end{IEEEproof}      

Firstly, we note that $m_{H^3} \geq \paren{m_{H^2}}^\frac32$ by Jensen's inequality, with equality for deterministic fading coefficients.  Therefore, the bounds given in Theorem \ref{Thm: Rates of Convergence} become tighter if the fading distribution is more concentrated around a point, {\em i.e.,} fading distributions with more restricted dynamic ranges.  Secondly, the Gaussian approximation bound derived in Theorem \ref{Thm: Rates of Convergence} is a combination of two different types of Berry-Esseen bounds, one of which is a {\em uniform} bound and the other one is a {\em non-uniform} bound.  The non-uniform bound is designed to be tight for large values of $|x|$.  On the other hand, the uniform bound is tighter for moderate values of $|x|$.  This point will be further illustrated in detail in the next section by using tangible examples for the transmitter location process $\Phi_\Lambda$.  It should also be noted that these convergence rates depend on the reference point $\vec{X}_0$ at which we measure the interference power level.  One easy corollary of Theorem \ref{Thm: Rates of Convergence} is the following uniform approximation bound.   
\begin{corollary}
\begin{eqnarray}
\sup_{x \in \Re} \left| \PR{\frac{I_\lambda - \ES{I_\lambda}}{\sqrt{\V{I_\lambda}}} \leq x} - \Psi(x)\right| \leq \frac{c}{\sqrt{\lambda}}, \label{Eqn: Uniform Rates}
\end{eqnarray}
where $c =0.4785 \frac{m_{H^3}}{\paren{m_{H^2}}^\frac32}\frac{\int_0^\infty G^3(t) p(t) \diff t}{\paren{\int_0^\infty G^2(t) p(t) \diff t}^\frac32}$.
\end{corollary}

\section{Gaussian Approximation Bounds for Stationary PPPs} \label{Section: Approximation Bounds for Stationary PPPs}
In this part, we take $\Phi_\Lambda$ to be a {\em stationary} PPP on $\R^2$ with intensity $\lambda$ [nodes per unit area], {\em i.e.,} $\Lambda\paren{\mathcal{A}} = \lambda \nu\paren{\mathcal{A}}$ and $\nu\paren{\mathcal{A}}$ is the area of the set $\mathcal{A}$, to numerically illustrate the validity of our bounds derived in Section \ref{Section: WMAI Distribution}.  Further applications to approximate the WMAI distribution for a non-stationary PPP are demonstrated in Appendix \ref{Appendix: Normal Approximation Non-stationary PPP}.  

\subsection{Theoretical Results}
There are several equivalent ways to represent $\Phi_\Lambda$ in this case.  The most convenient representation for our purposes in this paper is the one obtained by transforming and marking (see \cite{Kingman93} for the details of marking and transforming of PPPs) a {\em stationary} PPP with {\em intensity} $1$ on $[0, \infty)$, which is given as
\begin{eqnarray}
\Phi_\Lambda \disteq \sum_{k \geq 1} \delta_{\paren{\sqrt{\frac{\Gamma_k}{\lambda \pi}} \cos\paren{U_k}, \sqrt{\frac{\Gamma_k}{\lambda \pi}} \sin\paren{U_k}}}, 
\end{eqnarray}
where $X \disteq Y$ means two random variables $X$ and $Y$ are equal in distribution, $U_k$'s are i.i.d. random variables uniformly distributed over $[0, 2\pi]$, and $\Gamma_k = \sum_{i=1}^k E_i$, where $E_i$'s are i.i.d. random variables with unit exponential distribution.  This representation allows us to take $\vec{X}_k \disteq \paren{\sqrt{\frac{\Gamma_k}{\lambda \pi}} \cos\paren{U_k}, \sqrt{\frac{\Gamma_k}{\lambda \pi}} \sin\paren{U_k}}$, and $\|\vec{X}_k\|_2 \disteq \sqrt{\frac{\Gamma_k}{\lambda \pi}}$.  We assume that the test receiver is placed at the origin, {\em i.e.,} $\vec{X}_o = \vec{0}$.  By using Poisson process  transformations one more time, one can further show that the distances between the origin and the points of $\Phi_\Lambda$ 
form a PPP on $[0, \infty)$ with mean measure $\Lambda \circ T^{-1}\paren{[0,t]} = \lambda \pi t^2$, and the density $p_\lambda(t) = 2 \lambda \pi t \I{t \geq 0}$.  Hence, the WMAI distribution for stationary PPPs can be approximated as in the following theorem.
\begin{theorem} \label{Thm: Rates of Convergence - Stationary}
Assume $\Phi_\Lambda$ is a stationary PPP with intensity $\lambda$ transmitters per unit area.  Then, for all $x \in \Re$,
\begin{eqnarray}
\left| \PR{\frac{I_\lambda - \ES{I_\lambda}}{\sqrt{\V{I_\lambda}}} \leq x} - \Psi(x)\right| \leq \frac{c(x)}{\sqrt{\lambda}}, \label{Eqn: Uniform Rates Stationary PPP}
\end{eqnarray}
where $\Psi(x) = \frac{1}{\sqrt{2 \pi}}\int_{-\infty}^x \e{-\frac{t^2}{2}} dt$ and $c(x) = \frac{1}{\sqrt{2 \pi}}\frac{m_{H^3}}{\paren{m_{H^2}}^\frac32}\frac{\int_0^\infty G^3(t) t \diff t}{\paren{\int_0^\infty G^2(t) t \diff t}^\frac32} \min\paren{0.4785, \frac{31.935}{1+|x|^3}}$. 
\end{theorem}
\begin{IEEEproof}
Directly follows from Theorem \ref{Thm: Rates of Convergence} after substituting $2\pi t\I{t\geq 0}$ for $p(t)$. 
\end{IEEEproof} 

\subsection{Simulation and Numerical Study}
Now, we present our numerically computed Gaussian approximation bounds and simulation results confirming the theoretical predictions above.  We will use two different path-loss models $G_1(t) = \frac{1}{(1+t)^\alpha}$ and $G_2(t) = \frac{1}{1+t^\alpha}$ with various values of $\alpha$.  Similar conclusions continue to hold for other path-loss models.  

\begin{figure*}[!t]
\begin{minipage}[t]{\textwidth}
\begin{center}
%
%
\begin{psfrags}%
\psfragscanon%
%
\psfrag{s09}[t][t]{\color[rgb]{0,0,0}\setlength{\tabcolsep}{0pt}\scriptsize \begin{tabular}{c}Centered and Normalized Interference Power\end{tabular}}%
\psfrag{s10}[b][b]{\color[rgb]{0,0,0}\setlength{\tabcolsep}{0pt}\scriptsize \begin{tabular}{c}Normal CDF and Bounds\end{tabular}}%
\psfrag{s12}[b][b]{\color[rgb]{0,0,0}\setlength{\tabcolsep}{0pt}\scriptsize \begin{tabular}{c}$G_1(t) = \frac{1}{(1+t)^{4}}$ (without fading) \end{tabular}}%
\psfrag{s14}[][]{\color[rgb]{0,0,0}\setlength{\tabcolsep}{0pt}\begin{tabular}{c} $ $ \end{tabular}}%
\psfrag{s15}[][]{\color[rgb]{0,0,0}\setlength{\tabcolsep}{0pt}\begin{tabular}{c} $ $\end{tabular}}%
\psfrag{s17}[l][l]{\color[rgb]{0,0,0}\scriptsize Normal}%
\psfrag{s18}[l][l]{\color[rgb]{0,0,0}\scriptsize $\lambda = 5$ (LB)}%
\psfrag{s19}[l][l]{\color[rgb]{0,0,0}\scriptsize $\lambda = 5$ (UB)}%
\psfrag{s20}[l][l]{\color[rgb]{0,0,0}\scriptsize $\lambda = 25$ (LB)}%
\psfrag{s21}[l][l]{\color[rgb]{0,0,0}\scriptsize $\lambda = 25$ (UB)}%
\psfrag{s22}[l][l]{\color[rgb]{0,0,0}\scriptsize $\lambda = 100$ (LB)}%
\psfrag{s23}[l][l]{\color[rgb]{0,0,0}\scriptsize $\lambda = 100$ (UB)}%
%
\psfrag{x01}[t][t]{\scriptsize 0}%
\psfrag{x02}[t][t]{\scriptsize 0.1}%
\psfrag{x03}[t][t]{\scriptsize 0.2}%
\psfrag{x04}[t][t]{\scriptsize 0.3}%
\psfrag{x05}[t][t]{\scriptsize 0.4}%
\psfrag{x06}[t][t]{\scriptsize 0.5}%
\psfrag{x07}[t][t]{\scriptsize 0.6}%
\psfrag{x08}[t][t]{\scriptsize 0.7}%
\psfrag{x09}[t][t]{\scriptsize 0.8}%
\psfrag{x10}[t][t]{\scriptsize 0.9}%
\psfrag{x11}[t][t]{\scriptsize 1}%
\psfrag{x12}[t][t]{\scriptsize -6}%
\psfrag{x13}[t][t]{\scriptsize -5}%
\psfrag{x14}[t][t]{\scriptsize -4}%
\psfrag{x15}[t][t]{\scriptsize -3}%
\psfrag{x16}[t][t]{\scriptsize -2}%
\psfrag{x17}[t][t]{\scriptsize -1}%
\psfrag{x18}[t][t]{\scriptsize 0}%
\psfrag{x19}[t][t]{\scriptsize 1}%
\psfrag{x20}[t][t]{\scriptsize 2}%
\psfrag{x21}[t][t]{\scriptsize 3}%
\psfrag{x22}[t][t]{\scriptsize 4}%
\psfrag{x23}[t][t]{\scriptsize 5}%
\psfrag{x24}[t][t]{\scriptsize 6}%
%
\psfrag{v01}[r][r]{\scriptsize 0}%
\psfrag{v02}[r][r]{\scriptsize 0.1}%
\psfrag{v03}[r][r]{\scriptsize 0.2}%
\psfrag{v04}[r][r]{\scriptsize 0.3}%
\psfrag{v05}[r][r]{\scriptsize 0.4}%
\psfrag{v06}[r][r]{\scriptsize 0.5}%
\psfrag{v07}[r][r]{\scriptsize 0.6}%
\psfrag{v08}[r][r]{\scriptsize 0.7}%
\psfrag{v09}[r][r]{\scriptsize 0.8}%
\psfrag{v10}[r][r]{\scriptsize 0.9}%
\psfrag{v11}[r][r]{\scriptsize 1}%
\psfrag{v12}[r][r]{\scriptsize 0}%
\psfrag{v13}[r][r]{\scriptsize 0.2}%
\psfrag{v14}[r][r]{\scriptsize 0.4}%
\psfrag{v15}[r][r]{\scriptsize 0.6}%
\psfrag{v16}[r][r]{\scriptsize 0.8}%
\psfrag{v17}[r][r]{\scriptsize 1}%
%
\includegraphics[scale=0.55]{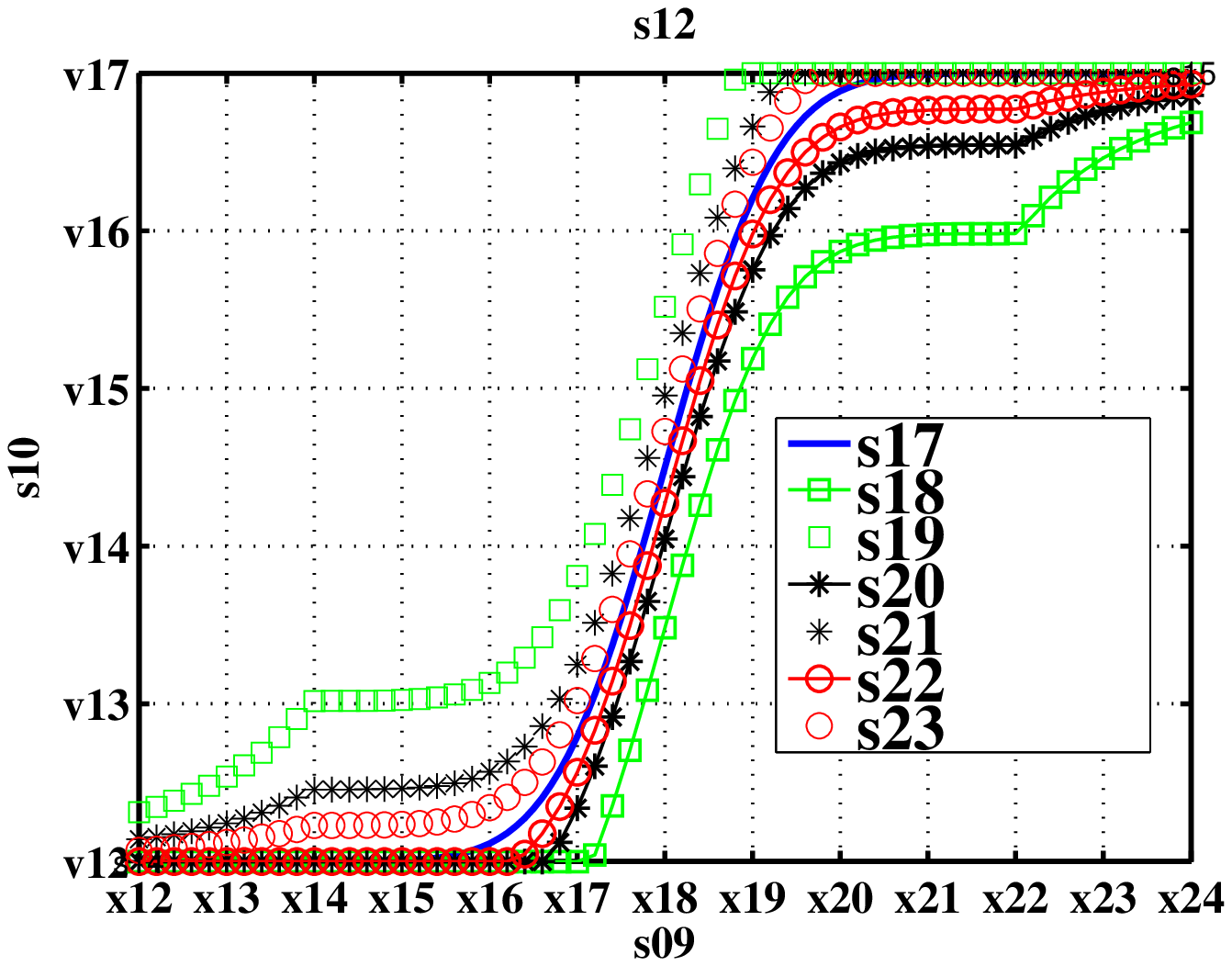}%
\end{psfrags}%
%

%
\hspace{\fill}
%
%
%
\begin{psfrags}%
\psfragscanon%
%
\psfrag{s09}[t][t]{\color[rgb]{0,0,0}\setlength{\tabcolsep}{0pt}\scriptsize \begin{tabular}{c}Centered and Normalized Interference Power\end{tabular}}%
\psfrag{s10}[b][b]{\color[rgb]{0,0,0}\setlength{\tabcolsep}{0pt}\scriptsize \begin{tabular}{c}Normal CDF and Bounds\end{tabular}}%
\psfrag{s12}[b][b]{\color[rgb]{0,0,0}\setlength{\tabcolsep}{0pt}\scriptsize \begin{tabular}{c}$G_2(t) = \frac{1}{1+t^{4}}$ (without fading) \end{tabular}}%
\psfrag{s14}[][]{\color[rgb]{0,0,0}\setlength{\tabcolsep}{0pt}\begin{tabular}{c} $ $ \end{tabular}}%
\psfrag{s15}[][]{\color[rgb]{0,0,0}\setlength{\tabcolsep}{0pt}\begin{tabular}{c} $ $ \end{tabular}}%
\psfrag{s17}[l][l]{\color[rgb]{0,0,0}\scriptsize Normal}%
\psfrag{s18}[l][l]{\color[rgb]{0,0,0}\scriptsize $\lambda = 5$ (LB)}%
\psfrag{s19}[l][l]{\color[rgb]{0,0,0}\scriptsize $\lambda = 5$ (UB)}%
\psfrag{s20}[l][l]{\color[rgb]{0,0,0}\scriptsize $\lambda = 25$ (LB)}%
\psfrag{s21}[l][l]{\color[rgb]{0,0,0}\scriptsize $\lambda = 25$ (UB)}%
\psfrag{s22}[l][l]{\color[rgb]{0,0,0}\scriptsize $\lambda = 100$ (LB)}%
\psfrag{s23}[l][l]{\color[rgb]{0,0,0}\scriptsize $\lambda = 100$ (UB)}%
%
\psfrag{x01}[t][t]{\scriptsize 0}%
\psfrag{x02}[t][t]{\scriptsize 0.1}%
\psfrag{x03}[t][t]{\scriptsize 0.2}%
\psfrag{x04}[t][t]{\scriptsize 0.3}%
\psfrag{x05}[t][t]{\scriptsize 0.4}%
\psfrag{x06}[t][t]{\scriptsize 0.5}%
\psfrag{x07}[t][t]{\scriptsize 0.6}%
\psfrag{x08}[t][t]{\scriptsize 0.7}%
\psfrag{x09}[t][t]{\scriptsize 0.8}%
\psfrag{x10}[t][t]{\scriptsize 0.9}%
\psfrag{x11}[t][t]{\scriptsize 1}%
\psfrag{x12}[t][t]{\scriptsize -6}%
\psfrag{x13}[t][t]{\scriptsize -5}%
\psfrag{x14}[t][t]{\scriptsize -4}%
\psfrag{x15}[t][t]{\scriptsize -3}%
\psfrag{x16}[t][t]{\scriptsize -2}%
\psfrag{x17}[t][t]{\scriptsize -1}%
\psfrag{x18}[t][t]{\scriptsize 0}%
\psfrag{x19}[t][t]{\scriptsize 1}%
\psfrag{x20}[t][t]{\scriptsize 2}%
\psfrag{x21}[t][t]{\scriptsize 3}%
\psfrag{x22}[t][t]{\scriptsize 4}%
\psfrag{x23}[t][t]{\scriptsize 5}%
\psfrag{x24}[t][t]{\scriptsize 6}%
%
\psfrag{v01}[r][r]{\scriptsize 0}%
\psfrag{v02}[r][r]{\scriptsize 0.1}%
\psfrag{v03}[r][r]{\scriptsize 0.2}%
\psfrag{v04}[r][r]{\scriptsize 0.3}%
\psfrag{v05}[r][r]{\scriptsize 0.4}%
\psfrag{v06}[r][r]{\scriptsize 0.5}%
\psfrag{v07}[r][r]{\scriptsize 0.6}%
\psfrag{v08}[r][r]{\scriptsize 0.7}%
\psfrag{v09}[r][r]{\scriptsize 0.8}%
\psfrag{v10}[r][r]{\scriptsize 0.9}%
\psfrag{v11}[r][r]{\scriptsize 1}%
\psfrag{v12}[r][r]{\scriptsize 0}%
\psfrag{v13}[r][r]{\scriptsize 0.2}%
\psfrag{v14}[r][r]{\scriptsize 0.4}%
\psfrag{v15}[r][r]{\scriptsize 0.6}%
\psfrag{v16}[r][r]{\scriptsize 0.8}%
\psfrag{v17}[r][r]{\scriptsize 1}%
%
\includegraphics[scale=0.55]{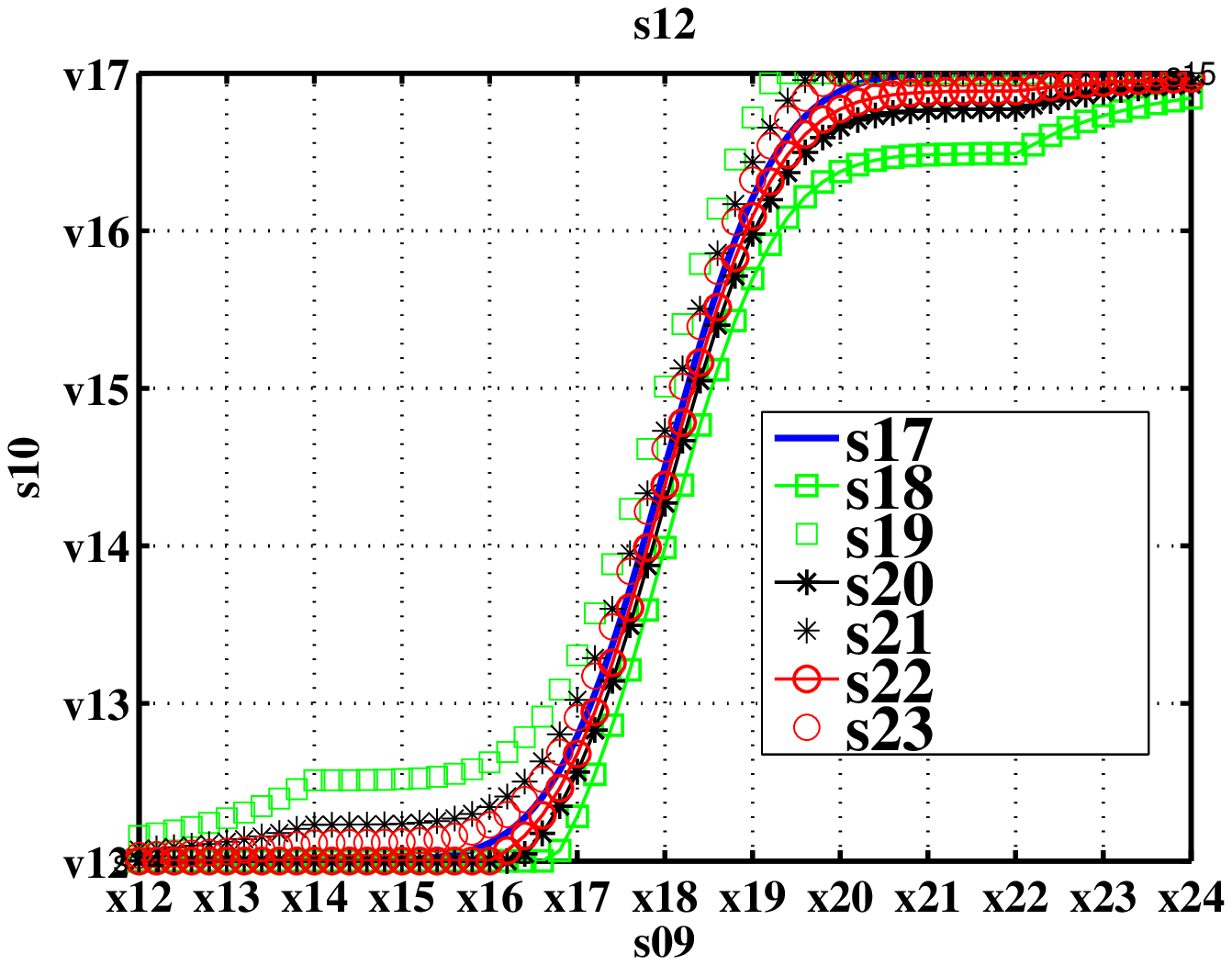}%
\end{psfrags}%
%

\end{center}
\end{minipage} 
\\ 
\begin{minipage}[t]{\textwidth}
\vspace{-2mm} 
\begin{center}
%
%
\begin{psfrags}%
\psfragscanon%
%
\psfrag{s09}[t][t]{\color[rgb]{0,0,0}\setlength{\tabcolsep}{0pt}\scriptsize \begin{tabular}{c}Centered and Normalized Interference Power\end{tabular}}%
\psfrag{s10}[b][b]{\color[rgb]{0,0,0}\setlength{\tabcolsep}{0pt}\scriptsize \begin{tabular}{c}Normal CDF and Bounds\end{tabular}}%
\psfrag{s12}[b][b]{\color[rgb]{0,0,0}\setlength{\tabcolsep}{0pt}\scriptsize \begin{tabular}{c}$G_1(t) = \frac{1}{(1+t)^4}$ (with fading)\end{tabular}}%
\psfrag{s14}[][]{\color[rgb]{0,0,0}\setlength{\tabcolsep}{0pt}\begin{tabular}{c} $ $ \end{tabular}}%
\psfrag{s15}[][]{\color[rgb]{0,0,0}\setlength{\tabcolsep}{0pt}\begin{tabular}{c} $ $ \end{tabular}}%
\psfrag{s17}[l][l]{\color[rgb]{0,0,0}\scriptsize Normal}%
\psfrag{s18}[l][l]{\color[rgb]{0,0,0}\scriptsize $\lambda = 5$ (LB)}%
\psfrag{s19}[l][l]{\color[rgb]{0,0,0}\scriptsize $\lambda = 5$ (UB)}%
\psfrag{s20}[l][l]{\color[rgb]{0,0,0}\scriptsize $\lambda = 25$ (LB)}%
\psfrag{s21}[l][l]{\color[rgb]{0,0,0}\scriptsize $\lambda = 25$ (UB)}%
\psfrag{s22}[l][l]{\color[rgb]{0,0,0}\scriptsize $\lambda = 100$ (LB)}%
\psfrag{s23}[l][l]{\color[rgb]{0,0,0}\scriptsize $\lambda = 100$ (UB)}%
%
\psfrag{x01}[t][t]{\scriptsize 0}%
\psfrag{x02}[t][t]{\scriptsize 0.1}%
\psfrag{x03}[t][t]{\scriptsize 0.2}%
\psfrag{x04}[t][t]{\scriptsize 0.3}%
\psfrag{x05}[t][t]{\scriptsize 0.4}%
\psfrag{x06}[t][t]{\scriptsize 0.5}%
\psfrag{x07}[t][t]{\scriptsize 0.6}%
\psfrag{x08}[t][t]{\scriptsize 0.7}%
\psfrag{x09}[t][t]{\scriptsize 0.8}%
\psfrag{x10}[t][t]{\scriptsize 0.9}%
\psfrag{x11}[t][t]{\scriptsize 1}%
\psfrag{x12}[t][t]{\scriptsize -6}%
\psfrag{x13}[t][t]{\scriptsize -5}%
\psfrag{x14}[t][t]{\scriptsize -4}%
\psfrag{x15}[t][t]{\scriptsize -3}%
\psfrag{x16}[t][t]{\scriptsize -2}%
\psfrag{x17}[t][t]{\scriptsize -1}%
\psfrag{x18}[t][t]{\scriptsize 0}%
\psfrag{x19}[t][t]{\scriptsize 1}%
\psfrag{x20}[t][t]{\scriptsize 2}%
\psfrag{x21}[t][t]{\scriptsize 3}%
\psfrag{x22}[t][t]{\scriptsize 4}%
\psfrag{x23}[t][t]{\scriptsize 5}%
\psfrag{x24}[t][t]{\scriptsize 6}%
%
\psfrag{v01}[r][r]{\scriptsize 0}%
\psfrag{v02}[r][r]{\scriptsize 0.1}%
\psfrag{v03}[r][r]{\scriptsize 0.2}%
\psfrag{v04}[r][r]{\scriptsize 0.3}%
\psfrag{v05}[r][r]{\scriptsize 0.4}%
\psfrag{v06}[r][r]{\scriptsize 0.5}%
\psfrag{v07}[r][r]{\scriptsize 0.6}%
\psfrag{v08}[r][r]{\scriptsize 0.7}%
\psfrag{v09}[r][r]{\scriptsize 0.8}%
\psfrag{v10}[r][r]{\scriptsize 0.9}%
\psfrag{v11}[r][r]{\scriptsize 1}%
\psfrag{v12}[r][r]{\scriptsize 0}%
\psfrag{v13}[r][r]{\scriptsize 0.2}%
\psfrag{v14}[r][r]{\scriptsize 0.4}%
\psfrag{v15}[r][r]{\scriptsize 0.6}%
\psfrag{v16}[r][r]{\scriptsize 0.8}%
\psfrag{v17}[r][r]{\scriptsize 1}%
%
\includegraphics[scale=0.55]{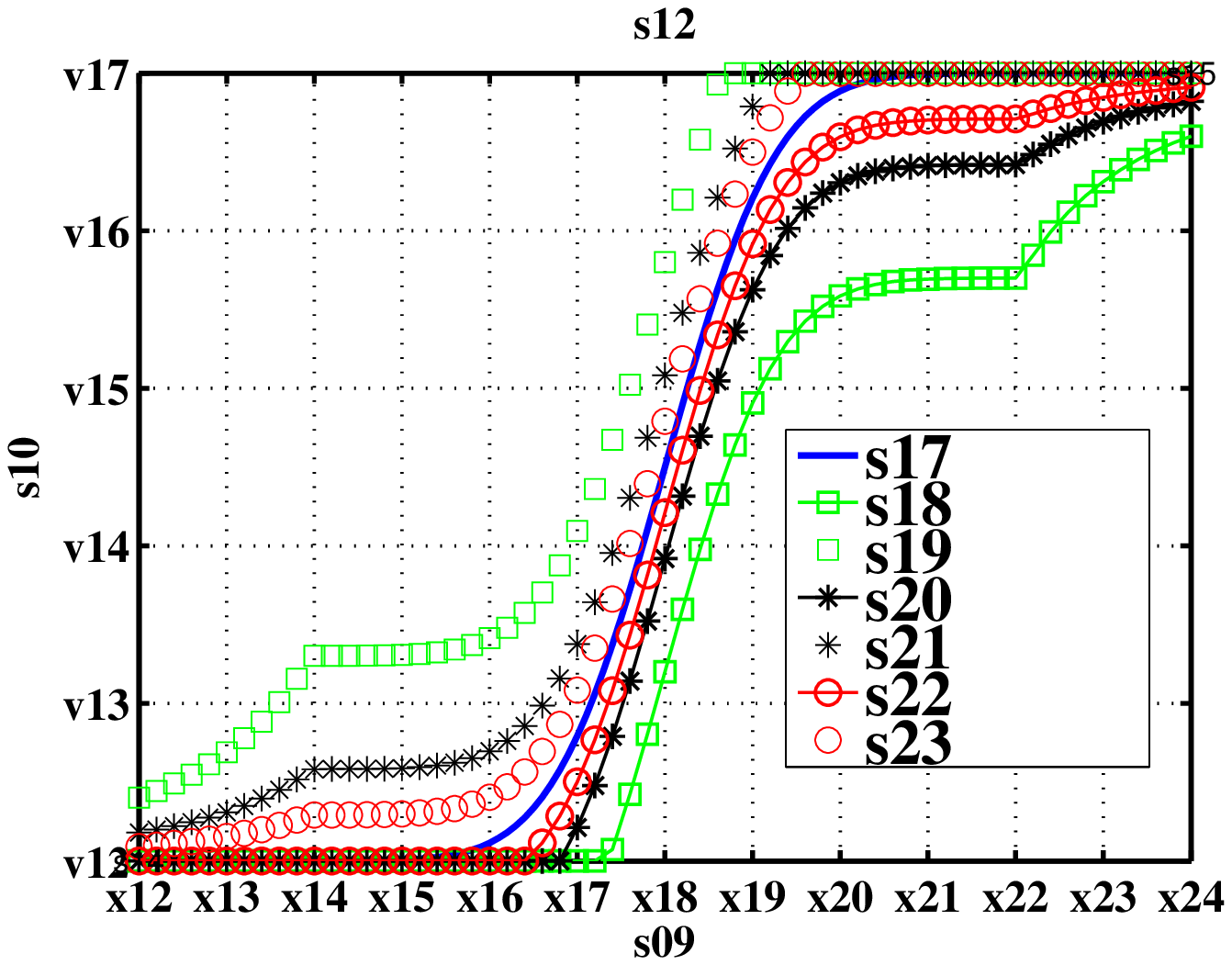}%
\end{psfrags}%
%

%
\hfill
%
%
%
\begin{psfrags}%
\psfragscanon%
%
\psfrag{s09}[t][t]{\color[rgb]{0,0,0}\setlength{\tabcolsep}{0pt}\scriptsize \begin{tabular}{c}Centered and Normalized Interference Power\end{tabular}}%
\psfrag{s10}[b][b]{\color[rgb]{0,0,0}\setlength{\tabcolsep}{0pt}\scriptsize \begin{tabular}{c}Normal CDF and Bounds\end{tabular}}%
\psfrag{s12}[b][b]{\color[rgb]{0,0,0}\setlength{\tabcolsep}{0pt}\scriptsize \begin{tabular}{c}$G_2(t) = \frac{1}{1+t^{4}}$ (with fading)\end{tabular}}%
\psfrag{s14}[][]{\color[rgb]{0,0,0}\setlength{\tabcolsep}{0pt}\begin{tabular}{c} $ $ \end{tabular}}%
\psfrag{s15}[][]{\color[rgb]{0,0,0}\setlength{\tabcolsep}{0pt}\begin{tabular}{c} $ $ \end{tabular}}%
\psfrag{s17}[l][l]{\color[rgb]{0,0,0}\scriptsize Normal}%
\psfrag{s18}[l][l]{\color[rgb]{0,0,0}\scriptsize $\lambda = 5$ (LB)}%
\psfrag{s19}[l][l]{\color[rgb]{0,0,0}\scriptsize $\lambda = 5$ (UB)}%
\psfrag{s20}[l][l]{\color[rgb]{0,0,0}\scriptsize $\lambda = 25$ (LB)}%
\psfrag{s21}[l][l]{\color[rgb]{0,0,0}\scriptsize $\lambda = 25$ (UB)}%
\psfrag{s22}[l][l]{\color[rgb]{0,0,0}\scriptsize $\lambda = 100$ (LB)}%
\psfrag{s23}[l][l]{\color[rgb]{0,0,0}\scriptsize $\lambda = 100$ (UB)}%
%
\psfrag{x01}[t][t]{\scriptsize 0}%
\psfrag{x02}[t][t]{\scriptsize 0.1}%
\psfrag{x03}[t][t]{\scriptsize 0.2}%
\psfrag{x04}[t][t]{\scriptsize 0.3}%
\psfrag{x05}[t][t]{\scriptsize 0.4}%
\psfrag{x06}[t][t]{\scriptsize 0.5}%
\psfrag{x07}[t][t]{\scriptsize 0.6}%
\psfrag{x08}[t][t]{\scriptsize 0.7}%
\psfrag{x09}[t][t]{\scriptsize 0.8}%
\psfrag{x10}[t][t]{\scriptsize 0.9}%
\psfrag{x11}[t][t]{\scriptsize 1}%
\psfrag{x12}[t][t]{\scriptsize -6}%
\psfrag{x13}[t][t]{\scriptsize -5}%
\psfrag{x14}[t][t]{\scriptsize -4}%
\psfrag{x15}[t][t]{\scriptsize -3}%
\psfrag{x16}[t][t]{\scriptsize -2}%
\psfrag{x17}[t][t]{\scriptsize -1}%
\psfrag{x18}[t][t]{\scriptsize 0}%
\psfrag{x19}[t][t]{\scriptsize 1}%
\psfrag{x20}[t][t]{\scriptsize 2}%
\psfrag{x21}[t][t]{\scriptsize 3}%
\psfrag{x22}[t][t]{\scriptsize 4}%
\psfrag{x23}[t][t]{\scriptsize 5}%
\psfrag{x24}[t][t]{6}%
%
\psfrag{v01}[r][r]{\scriptsize 0}%
\psfrag{v02}[r][r]{\scriptsize 0.1}%
\psfrag{v03}[r][r]{\scriptsize 0.2}%
\psfrag{v04}[r][r]{\scriptsize 0.3}%
\psfrag{v05}[r][r]{\scriptsize 0.4}%
\psfrag{v06}[r][r]{\scriptsize 0.5}%
\psfrag{v07}[r][r]{\scriptsize 0.6}%
\psfrag{v08}[r][r]{\scriptsize 0.7}%
\psfrag{v09}[r][r]{\scriptsize 0.8}%
\psfrag{v10}[r][r]{\scriptsize 0.9}%
\psfrag{v11}[r][r]{\scriptsize 1}%
\psfrag{v12}[r][r]{\scriptsize 0}%
\psfrag{v13}[r][r]{\scriptsize 0.2}%
\psfrag{v14}[r][r]{\scriptsize 0.4}%
\psfrag{v15}[r][r]{\scriptsize 0.6}%
\psfrag{v16}[r][r]{\scriptsize 0.8}%
\psfrag{v17}[r][r]{\scriptsize 1}%
%
\includegraphics[scale=0.55]{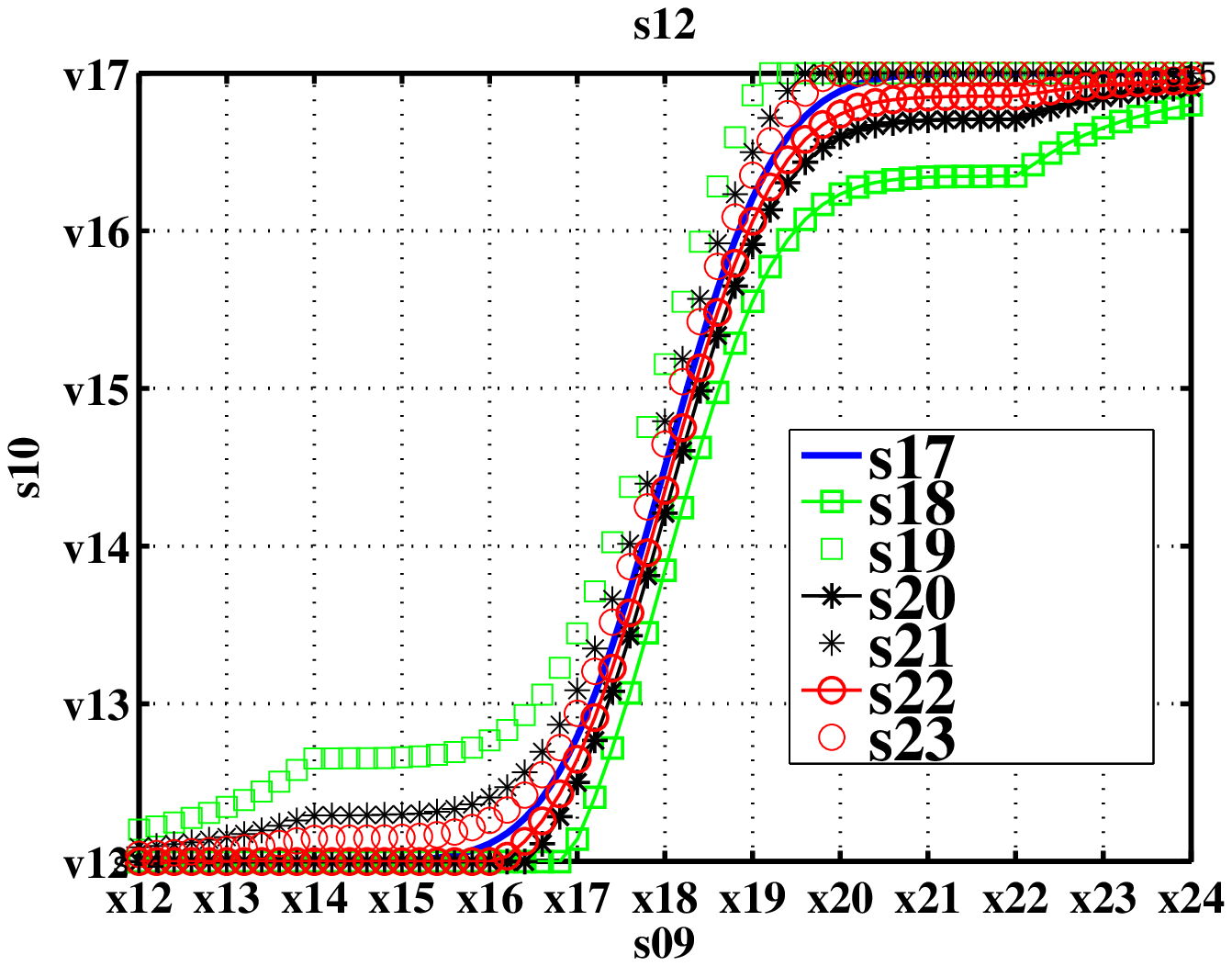}%
\end{psfrags}%
%

\end{center}
\end{minipage}
\caption{Upper and lower bounds on the centered and normalized WMAI CDFs for the path-loss functions $G_1(t) = \frac{1}{\paren{1+t}^\alpha}$ (lefthand side figures) and $G_2(t) = \frac{1}{1+t^\alpha}$ (righthand side figures).  The effect of fading is also illustrated in the bottom figures by assuming Nakagami-$m$ fading with $m$ parameter set to $5$. ($\alpha = 4$)} \label{Fig: CDF Bounds}
\end{figure*}

In Fig. \ref{Fig: CDF Bounds}, we plot our numerically computed bounds for the WMAI distributions both with (bottom figures) and without (top figures) fading. Nakagami-$m$ fading model \cite{Stuber96} with unit mean power gain and $m$ parameter set to $5$ is used to model the likely fading effects existing in the communication environment.  
We set the path-loss exponent $\alpha$ to $4$ in Fig.  \ref{Fig: CDF Bounds}. 

We observe two different regimes in our computed bounds for the WMAI distributions in Fig. \ref{Fig: CDF Bounds}.  For the moderate values of the centered and normalized WMAI, {\em i.e.,} $\frac{I_\lambda - \ES{I_\lambda}}{\sqrt{\V{I_\lambda}}}$, our uniform Berry-Esseen bound $\frac{0.4785}{\sqrt{2 \lambda \pi}}\frac{\int_0^\infty G^3(t) t dt}{\paren{\int_0^\infty G^2(t) t dt}^\frac32}$ gives better upper and lower bounds around the normal CDF for the interference distribution.  On the other hand, for large ({\em i.e.,} greater than $4$) absolute values of the centered and normalized WMAI, our non-uniform Berry-Esseen bound $\frac{31.935}{\sqrt{2 \lambda \pi}\paren{1+|x|^3}}\frac{\int_0^\infty G^3(t) t dt}{\paren{\int_0^\infty G^2(t) t dt}^\frac32}$ becomes a better estimator for the interference distribution.  Our bounds can be used to bound the probability of outage, outage capacity and etc. in a wireless communications setting as illustrated in Section \ref{Section: Applications}.  We remark that \eqref{Eqn: Uniform Rates Stationary PPP} is given with an almost eleven times larger constant ({\em i.e.,} $2.21$ rather than $\frac{0.4785}{\sqrt{2\pi}}$) in \cite{HS85} ({\em i.e.,} see Equation 7.1 in \cite{HS85}).  Hence, such bounds will be at least eleven times sharper than those based on the convergence results in \cite{HS85}. 

For any fixed value of $\lambda$, the gap between the upper and lower bounds vanish at a rate $\BO{|x|^{-3}}$ as $ x \ra \infty$, {\em i.e.,} as the interference power increases.  When $\lambda$ increases, the upper and lower bounds approach the normal CDF at a rate $\frac{1}{\sqrt{\lambda}}$, and we start to approximate the WMAI distribution as a Gaussian distribution increasingly more accurately.  When the upper and lower bounds on the WMAI distribution are compared for different path-loss models, we see that they become tighter for $G_2(t)$.  This is because the path-loss dependent constant, {\em i.e.,} $\frac{\int_0^\infty G^3(t) t dt}{\paren{\int_0^\infty G^2(t) t dt}^\frac32}$, appearing in Theorem \ref{Thm: Rates of Convergence - Stationary} is smaller for $G_2(t)$ than that for $G_1(t)$ ({\em i.e.,} see Table \ref{Table: 1} below).        

\begin{figure*}[!t]
\begin{minipage}[t]{\textwidth}
\begin{center}
%
%
\begin{psfrags}%
\psfragscanon%
%
\psfrag{s09}[t][t]{\color[rgb]{0,0,0}\setlength{\tabcolsep}{0pt}\scriptsize \begin{tabular}{c} $ $ \end{tabular}}%
\psfrag{s10}[b][b]{\color[rgb]{0,0,0}\setlength{\tabcolsep}{0pt}\scriptsize \begin{tabular}{c} \hspace{-25mm} Cumulative Distribution Function \end{tabular}}%
\psfrag{s12}[b][b]{\color[rgb]{0,0,0}\setlength{\tabcolsep}{0pt}\scriptsize \begin{tabular}{c} $G_1(t) = \frac{1}{(1+t)^\alpha}$ (without fading) \end{tabular}}%
\psfrag{s14}[][]{\color[rgb]{0,0,0}\setlength{\tabcolsep}{0pt}\scriptsize \begin{tabular}{c} $ $ \end{tabular}}%
\psfrag{s15}[][]{\color[rgb]{0,0,0}\setlength{\tabcolsep}{0pt}\begin{tabular}{c} $ $ \end{tabular}}%
\psfrag{s16}[t][t]{\color[rgb]{0,0,0}\setlength{\tabcolsep}{0pt}\scriptsize \begin{tabular}{c} Centered and Normalized Interference Power \end{tabular}}%
\psfrag{s17}[b][b]{\color[rgb]{0,0,0}\setlength{\tabcolsep}{0pt}\scriptsize \begin{tabular}{c} $ $ \end{tabular}}%
\psfrag{s21}[][]{\color[rgb]{0,0,0}\setlength{\tabcolsep}{0pt}\begin{tabular}{c} $ $ \end{tabular}}%
\psfrag{s22}[][]{\color[rgb]{0,0,0}\setlength{\tabcolsep}{0pt}\begin{tabular}{c} $ $ \end{tabular}}%
\psfrag{s23}[l][l]{\color[rgb]{0,0,0}\scriptsize $\lambda = 10$}%
\psfrag{s24}[l][l]{\color[rgb]{0,0,0}\scriptsize Normal}%
\psfrag{s25}[l][l]{\color[rgb]{0,0,0}\scriptsize $\lambda = 0.1$}%
\psfrag{s26}[l][l]{\color[rgb]{0,0,0}\scriptsize $\lambda = 1$}%
\psfrag{s27}[l][l]{\color[rgb]{0,0,0}\scriptsize $\lambda = 10$}%
\psfrag{s28}[l][l]{\color[rgb]{0,0,0}\scriptsize $\lambda = 10$}%
\psfrag{s29}[l][l]{\color[rgb]{0,0,0}\scriptsize Normal}%
\psfrag{s30}[l][l]{\color[rgb]{0,0,0}\scriptsize $\lambda = 0.1$}%
\psfrag{s31}[l][l]{\color[rgb]{0,0,0}\scriptsize $\lambda = 1$}%
\psfrag{s32}[l][l]{\color[rgb]{0,0,0}\scriptsize $\lambda = 10$}%
%
\psfrag{x01}[t][t]{\scriptsize 0}%
\psfrag{x02}[t][t]{\scriptsize 0.1}%
\psfrag{x03}[t][t]{\scriptsize 0.2}%
\psfrag{x04}[t][t]{\scriptsize 0.3}%
\psfrag{x05}[t][t]{\scriptsize 0.4}%
\psfrag{x06}[t][t]{\scriptsize 0.5}%
\psfrag{x07}[t][t]{\scriptsize 0.6}%
\psfrag{x08}[t][t]{\scriptsize 0.7}%
\psfrag{x09}[t][t]{\scriptsize 0.8}%
\psfrag{x10}[t][t]{\scriptsize 0.9}%
\psfrag{x11}[t][t]{\scriptsize 1}%
\psfrag{x12}[t][t]{\scriptsize -3}%
\psfrag{x13}[t][t]{\scriptsize -2}%
\psfrag{x14}[t][t]{\scriptsize -1}%
\psfrag{x15}[t][t]{\scriptsize 0}%
\psfrag{x16}[t][t]{\scriptsize 1}%
\psfrag{x17}[t][t]{\scriptsize 2}%
\psfrag{x18}[t][t]{\scriptsize 3}%
\psfrag{x19}[t][t]{\scriptsize -3}%
\psfrag{x20}[t][t]{\scriptsize -2}%
\psfrag{x21}[t][t]{\scriptsize -1}%
\psfrag{x22}[t][t]{\scriptsize 0}%
\psfrag{x23}[t][t]{\scriptsize 1}%
\psfrag{x24}[t][t]{\scriptsize 2}%
\psfrag{x25}[t][t]{\scriptsize 3}%
%
\psfrag{v01}[r][r]{\scriptsize 0}%
\psfrag{v02}[r][r]{\scriptsize 0.1}%
\psfrag{v03}[r][r]{\scriptsize 0.2}%
\psfrag{v04}[r][r]{\scriptsize 0.3}%
\psfrag{v05}[r][r]{\scriptsize 0.4}%
\psfrag{v06}[r][r]{\scriptsize 0.5}%
\psfrag{v07}[r][r]{\scriptsize 0.6}%
\psfrag{v08}[r][r]{\scriptsize 0.7}%
\psfrag{v09}[r][r]{\scriptsize 0.8}%
\psfrag{v10}[r][r]{\scriptsize 0.9}%
\psfrag{v11}[r][r]{\scriptsize 1}%
\psfrag{v12}[r][r]{\scriptsize 0}%
\psfrag{v13}[r][r]{\scriptsize 0.2}%
\psfrag{v14}[r][r]{\scriptsize 0.4}%
\psfrag{v15}[r][r]{\scriptsize 0.6}%
\psfrag{v16}[r][r]{\scriptsize 0.8}%
\psfrag{v17}[r][r]{\scriptsize 1}%
\psfrag{v18}[r][r]{\scriptsize 0}%
\psfrag{v19}[r][r]{\scriptsize 0.2}%
\psfrag{v20}[r][r]{\scriptsize 0.4}%
\psfrag{v21}[r][r]{\scriptsize 0.6}%
\psfrag{v22}[r][r]{\scriptsize 0.8}%
\psfrag{v23}[r][r]{\scriptsize 1}%
%
\setlength{\unitlength}{1cm}
\begin{picture}(8, 6)(0, 0)
\put(0, 0){\includegraphics[scale=0.55]{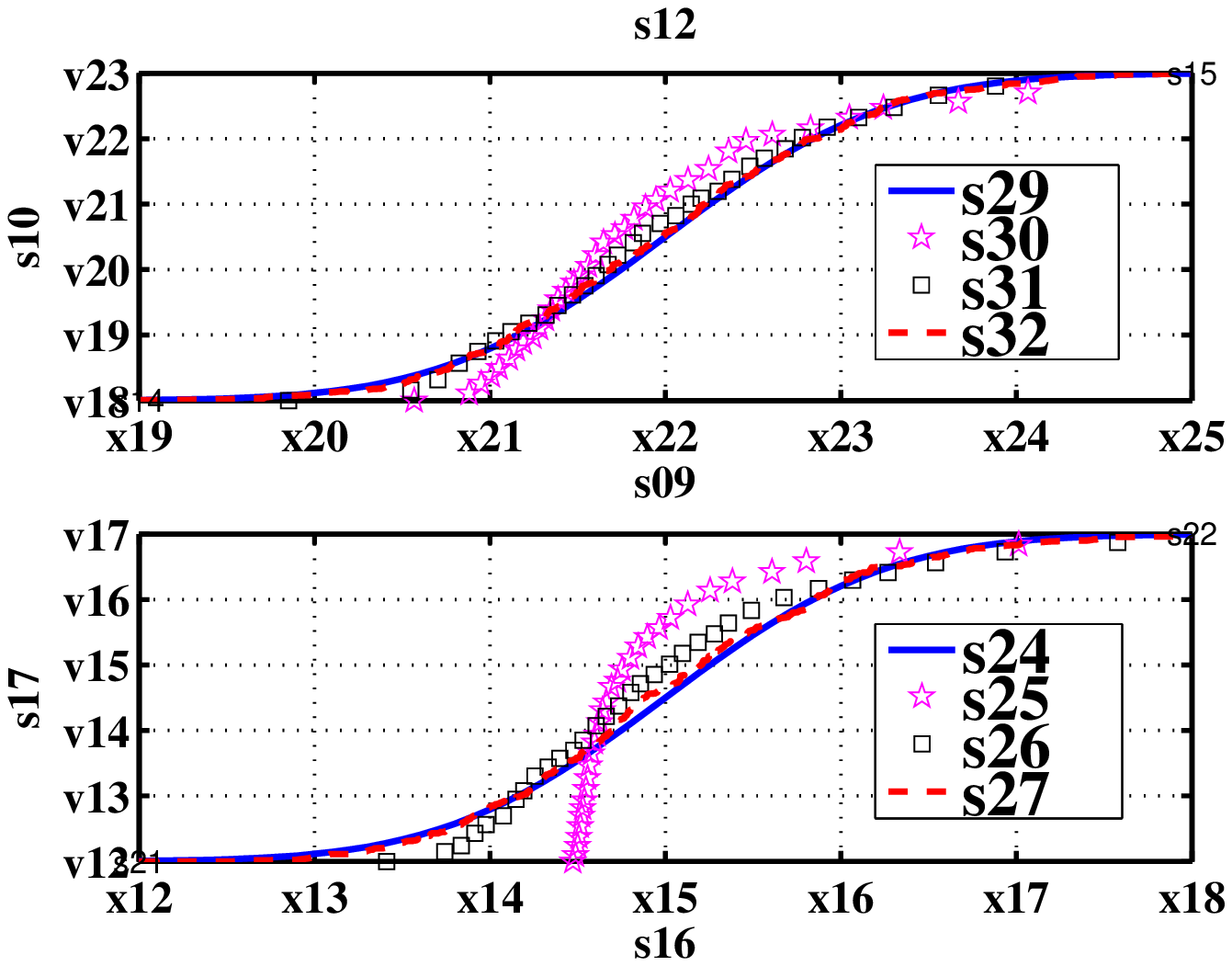}}%
\put(1.92, 4.35){\color[rgb]{0, 0, 0} \scriptsize \colorbox{white}{\fbox{$\alpha = 3$}}}
\put(1.92, 1.54){\color[rgb]{0, 0, 0} \scriptsize \colorbox{white}{\fbox{$\alpha = 5$}}}
\end{picture}
\end{psfrags}%
%

%
\hspace{\fill}
%
%
%
\begin{psfrags}%
\psfragscanon%
%
\psfrag{s13}[t][t]{\color[rgb]{0,0,0}\setlength{\tabcolsep}{0pt}\scriptsize \begin{tabular}{c} $ $ \end{tabular}}%
\psfrag{s14}[b][b]{\color[rgb]{0,0,0}\setlength{\tabcolsep}{0pt}\scriptsize \begin{tabular}{c} \hspace{-25mm} Cumulative Distribution Function \end{tabular}}%
\psfrag{s16}[b][b]{\color[rgb]{0,0,0}\setlength{\tabcolsep}{0pt}\scriptsize \begin{tabular}{c} $G_2(t) = \frac{1}{1+t^{\alpha}}$ (without fading)\end{tabular}}%
\psfrag{s18}[][]{\color[rgb]{0,0,0}\setlength{\tabcolsep}{0pt}\scriptsize \begin{tabular}{c} $ $ \end{tabular}}%
\psfrag{s19}[][]{\color[rgb]{0,0,0}\setlength{\tabcolsep}{0pt}\scriptsize \begin{tabular}{c} $ $ \end{tabular}}%
\psfrag{s20}[l][l]{\color[rgb]{0,0,0}\scriptsize $\lambda = 10$}%
\psfrag{s21}[l][l]{\color[rgb]{0,0,0}\scriptsize Normal}%
\psfrag{s22}[l][l]{\color[rgb]{0,0,0}\scriptsize $\lambda = 0.1$}%
\psfrag{s23}[l][l]{\color[rgb]{0,0,0}\scriptsize $\lambda = 1$}%
\psfrag{s24}[l][l]{\color[rgb]{0,0,0}\scriptsize $\lambda = 10$}%
\psfrag{s25}[t][t]{\color[rgb]{0,0,0}\scriptsize \setlength{\tabcolsep}{0pt}\begin{tabular}{c} Centered and Normalized Interference Power\end{tabular}}%
\psfrag{s26}[b][b]{\color[rgb]{0,0,0}\scriptsize \setlength{\tabcolsep}{0pt}\begin{tabular}{c} $ $ \end{tabular}}%
\psfrag{s30}[][]{\color[rgb]{0,0,0}\scriptsize \setlength{\tabcolsep}{0pt}\begin{tabular}{c} $ $ \end{tabular}}%
\psfrag{s31}[][]{\color[rgb]{0,0,0}\scriptsize \setlength{\tabcolsep}{0pt}\begin{tabular}{c} $ $ \end{tabular}}%
\psfrag{s32}[l][l]{\color[rgb]{0,0,0}\scriptsize $\lambda = 10$}%
\psfrag{s33}[l][l]{\color[rgb]{0,0,0}\scriptsize Normal}%
\psfrag{s34}[l][l]{\color[rgb]{0,0,0}\scriptsize $\lambda = 0.1$}%
\psfrag{s35}[l][l]{\color[rgb]{0,0,0}\scriptsize $\lambda = 1$}%
\psfrag{s36}[l][l]{\color[rgb]{0,0,0}\scriptsize $\lambda = 10$}%
%
\psfrag{x01}[t][t]{\scriptsize 0}%
\psfrag{x02}[t][t]{\scriptsize 0.1}%
\psfrag{x03}[t][t]{\scriptsize 0.2}%
\psfrag{x04}[t][t]{\scriptsize 0.3}%
\psfrag{x05}[t][t]{\scriptsize 0.4}%
\psfrag{x06}[t][t]{\scriptsize 0.5}%
\psfrag{x07}[t][t]{\scriptsize 0.6}%
\psfrag{x08}[t][t]{\scriptsize 0.7}%
\psfrag{x09}[t][t]{\scriptsize 0.8}%
\psfrag{x10}[t][t]{\scriptsize 0.9}%
\psfrag{x11}[t][t]{\scriptsize 1}%
\psfrag{x12}[t][t]{\scriptsize -3}%
\psfrag{x13}[t][t]{\scriptsize -2}%
\psfrag{x14}[t][t]{\scriptsize -1}%
\psfrag{x15}[t][t]{\scriptsize 0}%
\psfrag{x16}[t][t]{\scriptsize 1}%
\psfrag{x17}[t][t]{\scriptsize 2}%
\psfrag{x18}[t][t]{\scriptsize 3}%
\psfrag{x19}[t][t]{\scriptsize -3}%
\psfrag{x20}[t][t]{\scriptsize -2}%
\psfrag{x21}[t][t]{\scriptsize -1}%
\psfrag{x22}[t][t]{\scriptsize 0}%
\psfrag{x23}[t][t]{\scriptsize 1}%
\psfrag{x24}[t][t]{\scriptsize 2}%
\psfrag{x25}[t][t]{\scriptsize 3}%
%
\psfrag{v01}[r][r]{\scriptsize 0}%
\psfrag{v02}[r][r]{\scriptsize 0.1}%
\psfrag{v03}[r][r]{\scriptsize 0.2}%
\psfrag{v04}[r][r]{\scriptsize 0.3}%
\psfrag{v05}[r][r]{\scriptsize 0.4}%
\psfrag{v06}[r][r]{\scriptsize 0.5}%
\psfrag{v07}[r][r]{\scriptsize 0.6}%
\psfrag{v08}[r][r]{\scriptsize 0.7}%
\psfrag{v09}[r][r]{\scriptsize 0.8}%
\psfrag{v10}[r][r]{\scriptsize 0.9}%
\psfrag{v11}[r][r]{\scriptsize 1}%
\psfrag{v12}[r][r]{\scriptsize 0}%
\psfrag{v13}[r][r]{\scriptsize 0.2}%
\psfrag{v14}[r][r]{\scriptsize 0.4}%
\psfrag{v15}[r][r]{\scriptsize 0.6}%
\psfrag{v16}[r][r]{\scriptsize 0.8}%
\psfrag{v17}[r][r]{\scriptsize 1}%
\psfrag{v18}[r][r]{\scriptsize 0}%
\psfrag{v19}[r][r]{\scriptsize 0.2}%
\psfrag{v20}[r][r]{\scriptsize 0.4}%
\psfrag{v21}[r][r]{\scriptsize 0.6}%
\psfrag{v22}[r][r]{\scriptsize 0.8}%
\psfrag{v23}[r][r]{\scriptsize 1}%
%
\setlength{\unitlength}{1cm}
\begin{picture}(8, 6)(0, 0)
\includegraphics[scale=0.55]{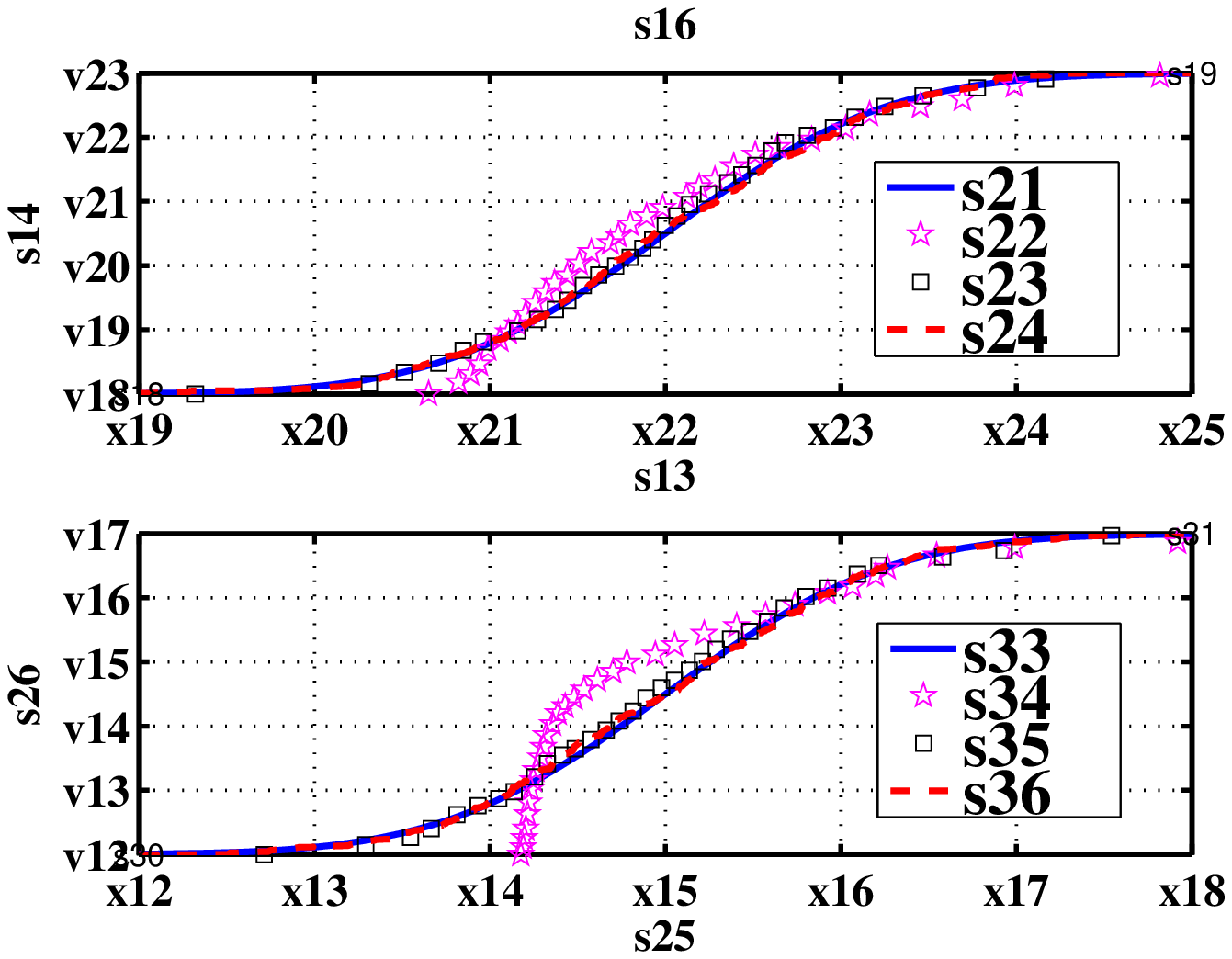}%
\put(-5.67, 4.35){\color[rgb]{0, 0, 0} \scriptsize \colorbox{white}{\fbox{$\alpha = 3$}}}
\put(-5.67, 1.54){\color[rgb]{0, 0, 0} \scriptsize \colorbox{white}{\fbox{$\alpha = 5$}}}
\end{picture}
\end{psfrags}%
%

\end{center}
\end{minipage} 
\\ 
\begin{minipage}[t]{\textwidth}
\vspace{-2mm} 
\begin{center}
%
%
\begin{psfrags}%
\psfragscanon%
%
\psfrag{s13}[t][t]{\color[rgb]{0,0,0}\setlength{\tabcolsep}{0pt}\scriptsize \begin{tabular}{c} $ $ \end{tabular}}%
\psfrag{s14}[b][b]{\color[rgb]{0,0,0}\setlength{\tabcolsep}{0pt}\scriptsize \begin{tabular}{c}\hspace{-25mm} Cumulative Distribution Function\end{tabular}}%
\psfrag{s16}[b][b]{\color[rgb]{0,0,0}\setlength{\tabcolsep}{0pt}\scriptsize \begin{tabular}{c}$G_1(t) = \frac{1}{(1+t)^\alpha}$ (with fading) \end{tabular}}%
\psfrag{s18}[][]{\color[rgb]{0,0,0}\setlength{\tabcolsep}{0pt}\scriptsize \begin{tabular}{c}$ $ \end{tabular}}%
\psfrag{s19}[][]{\color[rgb]{0,0,0}\setlength{\tabcolsep}{0pt}\scriptsize \begin{tabular}{c} $ $ \end{tabular}}%
\psfrag{s20}[l][l]{\color[rgb]{0,0,0}\scriptsize $\lambda = 10$}%
\psfrag{s21}[l][l]{\color[rgb]{0,0,0}\scriptsize Normal}%
\psfrag{s22}[l][l]{\color[rgb]{0,0,0}\scriptsize $\lambda = 0.1$}%
\psfrag{s23}[l][l]{\color[rgb]{0,0,0}\scriptsize $\lambda = 1$}%
\psfrag{s24}[l][l]{\color[rgb]{0,0,0}\scriptsize $\lambda = 10$}%
\psfrag{s25}[t][t]{\color[rgb]{0,0,0}\setlength{\tabcolsep}{0pt}\scriptsize \begin{tabular}{c}Centered and Normalized Interference Power\end{tabular}}%
\psfrag{s26}[b][b]{\color[rgb]{0,0,0}\setlength{\tabcolsep}{0pt}\scriptsize \begin{tabular}{c}$ $ \end{tabular}}%
\psfrag{s30}[][]{\color[rgb]{0,0,0}\setlength{\tabcolsep}{0pt}\scriptsize \begin{tabular}{c} \end{tabular}}%
\psfrag{s31}[][]{\color[rgb]{0,0,0}\setlength{\tabcolsep}{0pt}\scriptsize \begin{tabular}{c} $ $ \end{tabular}}%
\psfrag{s32}[l][l]{\color[rgb]{0,0,0}\scriptsize $\lambda = 10$}%
\psfrag{s33}[l][l]{\color[rgb]{0,0,0}\scriptsize Normal}%
\psfrag{s34}[l][l]{\color[rgb]{0,0,0}\scriptsize $\lambda = 0.1$}%
\psfrag{s35}[l][l]{\color[rgb]{0,0,0}\scriptsize $\lambda = 1$}%
\psfrag{s36}[l][l]{\color[rgb]{0,0,0}\scriptsize $\lambda = 10$}%
%
\psfrag{x01}[t][t]{\scriptsize 0}%
\psfrag{x02}[t][t]{\scriptsize 0.1}%
\psfrag{x03}[t][t]{\scriptsize 0.2}%
\psfrag{x04}[t][t]{\scriptsize 0.3}%
\psfrag{x05}[t][t]{\scriptsize 0.4}%
\psfrag{x06}[t][t]{\scriptsize 0.5}%
\psfrag{x07}[t][t]{\scriptsize 0.6}%
\psfrag{x08}[t][t]{\scriptsize 0.7}%
\psfrag{x09}[t][t]{\scriptsize 0.8}%
\psfrag{x10}[t][t]{\scriptsize 0.9}%
\psfrag{x11}[t][t]{\scriptsize 1}%
\psfrag{x12}[t][t]{\scriptsize -3}%
\psfrag{x13}[t][t]{\scriptsize -2}%
\psfrag{x14}[t][t]{\scriptsize -1}%
\psfrag{x15}[t][t]{\scriptsize 0}%
\psfrag{x16}[t][t]{\scriptsize 1}%
\psfrag{x17}[t][t]{\scriptsize 2}%
\psfrag{x18}[t][t]{\scriptsize 3}%
\psfrag{x19}[t][t]{\scriptsize -3}%
\psfrag{x20}[t][t]{\scriptsize -2}%
\psfrag{x21}[t][t]{\scriptsize -1}%
\psfrag{x22}[t][t]{\scriptsize 0}%
\psfrag{x23}[t][t]{\scriptsize 1}%
\psfrag{x24}[t][t]{\scriptsize 2}%
\psfrag{x25}[t][t]{\scriptsize 3}%
%
\psfrag{v01}[r][r]{\scriptsize 0}%
\psfrag{v02}[r][r]{\scriptsize 0.1}%
\psfrag{v03}[r][r]{\scriptsize 0.2}%
\psfrag{v04}[r][r]{\scriptsize 0.3}%
\psfrag{v05}[r][r]{\scriptsize 0.4}%
\psfrag{v06}[r][r]{\scriptsize 0.5}%
\psfrag{v07}[r][r]{\scriptsize 0.6}%
\psfrag{v08}[r][r]{\scriptsize 0.7}%
\psfrag{v09}[r][r]{\scriptsize 0.8}%
\psfrag{v10}[r][r]{\scriptsize 0.9}%
\psfrag{v11}[r][r]{\scriptsize 1}%
\psfrag{v12}[r][r]{\scriptsize 0}%
\psfrag{v13}[r][r]{\scriptsize 0.2}%
\psfrag{v14}[r][r]{\scriptsize 0.4}%
\psfrag{v15}[r][r]{\scriptsize 0.6}%
\psfrag{v16}[r][r]{\scriptsize 0.8}%
\psfrag{v17}[r][r]{\scriptsize 1}%
\psfrag{v18}[r][r]{\scriptsize 0}%
\psfrag{v19}[r][r]{\scriptsize 0.2}%
\psfrag{v20}[r][r]{\scriptsize 0.4}%
\psfrag{v21}[r][r]{\scriptsize 0.6}%
\psfrag{v22}[r][r]{\scriptsize 0.8}%
\psfrag{v23}[r][r]{\scriptsize 1}%
%
\setlength{\unitlength}{1cm}
\begin{picture}(8, 6)(0, 0)
\put(0, 0){\includegraphics[scale=0.55]{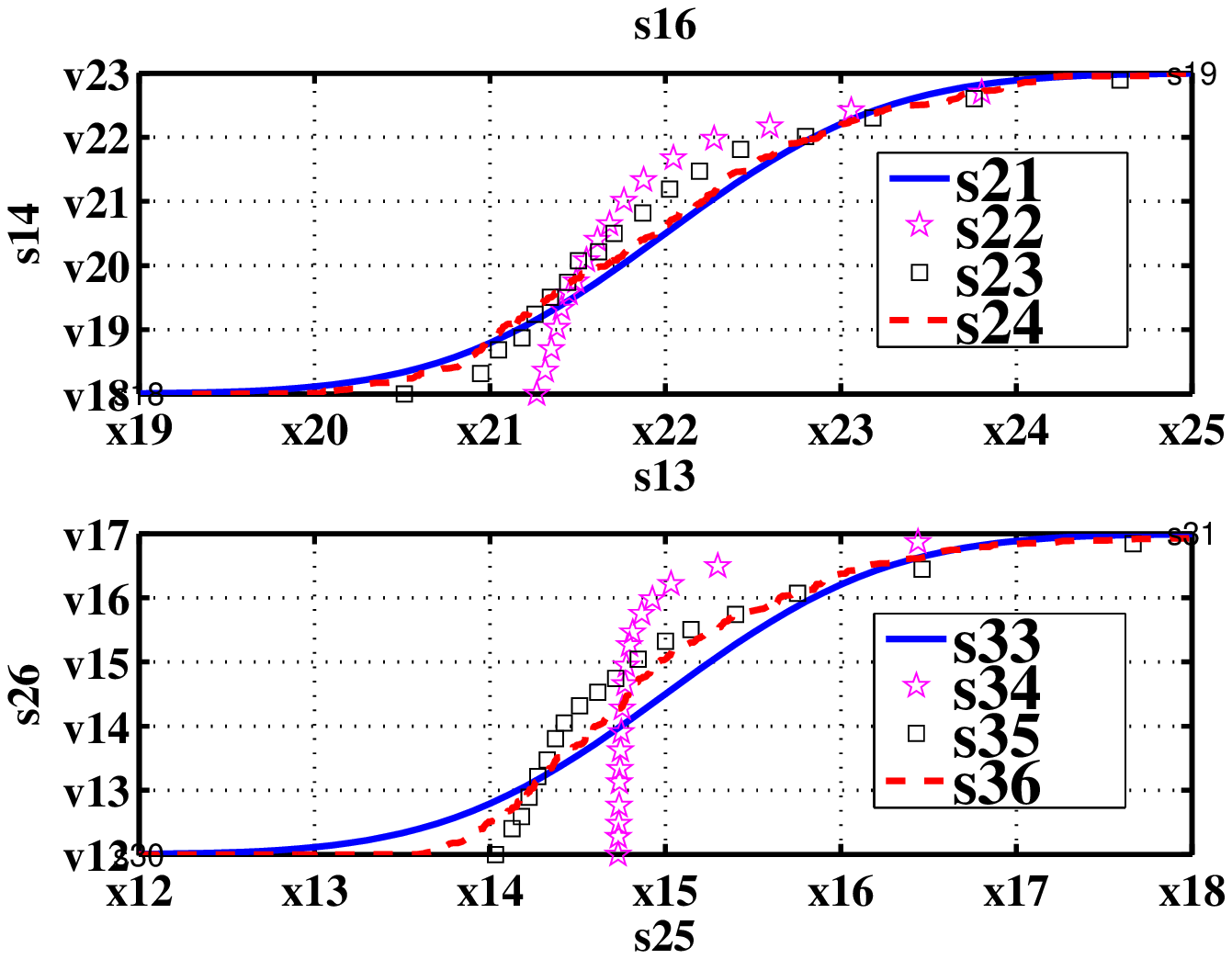}}%
\put(1.92, 4.35){\color[rgb]{0, 0, 0} \scriptsize \colorbox{white}{\fbox{$\alpha = 3$}}}
\put(1.92, 1.54){\color[rgb]{0, 0, 0} \scriptsize \colorbox{white}{\fbox{$\alpha = 5$}}}
\end{picture}
\end{psfrags}%
%

%
\hfill
%
%
%
\begin{psfrags}%
\psfragscanon%
%
\psfrag{s13}[t][t]{\color[rgb]{0,0,0}\setlength{\tabcolsep}{0pt}\scriptsize \begin{tabular}{c}$ $\end{tabular}}%
\psfrag{s14}[b][b]{\color[rgb]{0,0,0}\setlength{\tabcolsep}{0pt}\scriptsize \begin{tabular}{c}\hspace{-25mm} Cumulative Distribution Function \end{tabular}}%
\psfrag{s16}[b][b]{\color[rgb]{0,0,0}\setlength{\tabcolsep}{0pt} \scriptsize \begin{tabular}{c}$G_2(t) = \frac{1}{1+t^\alpha}$ (with fading) \end{tabular}}%
\psfrag{s18}[][]{\color[rgb]{0,0,0}\setlength{\tabcolsep}{0pt}\scriptsize \begin{tabular}{c} $ $ \end{tabular}}%
\psfrag{s19}[][]{\color[rgb]{0,0,0}\setlength{\tabcolsep}{0pt}\scriptsize \begin{tabular}{c} $ $ \end{tabular}}%
\psfrag{s20}[l][l]{\color[rgb]{0,0,0}\scriptsize $\lambda = 10$}%
\psfrag{s21}[l][l]{\color[rgb]{0,0,0}\scriptsize Normal}%
\psfrag{s22}[l][l]{\color[rgb]{0,0,0}\scriptsize $\lambda = 0.1$}%
\psfrag{s23}[l][l]{\color[rgb]{0,0,0}\scriptsize $\lambda = 1$}%
\psfrag{s24}[l][l]{\color[rgb]{0,0,0}\scriptsize $\lambda = 10$}%
\psfrag{s25}[t][t]{\color[rgb]{0,0,0}\setlength{\tabcolsep}{0pt}\scriptsize \begin{tabular}{c} Centered and Normalized Interference Power \end{tabular}}%
\psfrag{s26}[b][b]{\color[rgb]{0,0,0}\setlength{\tabcolsep}{0pt}\scriptsize \begin{tabular}{c} $ $ \end{tabular}}%
\psfrag{s30}[][]{\color[rgb]{0,0,0}\setlength{\tabcolsep}{0pt}\scriptsize \begin{tabular}{c} $ $ \end{tabular}}%
\psfrag{s31}[][]{\color[rgb]{0,0,0}\setlength{\tabcolsep}{0pt}\scriptsize \begin{tabular}{c} $ $ \end{tabular}}%
\psfrag{s32}[l][l]{\color[rgb]{0,0,0}\scriptsize $\lambda = 10$}%
\psfrag{s33}[l][l]{\color[rgb]{0,0,0}\scriptsize Normal}%
\psfrag{s34}[l][l]{\color[rgb]{0,0,0}\scriptsize $\lambda = 0.1$}%
\psfrag{s35}[l][l]{\color[rgb]{0,0,0}\scriptsize $\lambda = 1$}%
\psfrag{s36}[l][l]{\color[rgb]{0,0,0}\scriptsize $\lambda = 10$}%
%
\psfrag{x01}[t][t]{\scriptsize 0}%
\psfrag{x02}[t][t]{\scriptsize 0.1}%
\psfrag{x03}[t][t]{\scriptsize 0.2}%
\psfrag{x04}[t][t]{\scriptsize 0.3}%
\psfrag{x05}[t][t]{\scriptsize 0.4}%
\psfrag{x06}[t][t]{\scriptsize 0.5}%
\psfrag{x07}[t][t]{\scriptsize 0.6}%
\psfrag{x08}[t][t]{\scriptsize 0.7}%
\psfrag{x09}[t][t]{\scriptsize 0.8}%
\psfrag{x10}[t][t]{\scriptsize 0.9}%
\psfrag{x11}[t][t]{\scriptsize 1}%
\psfrag{x12}[t][t]{\scriptsize -3}%
\psfrag{x13}[t][t]{\scriptsize -2}%
\psfrag{x14}[t][t]{\scriptsize -1}%
\psfrag{x15}[t][t]{\scriptsize 0}%
\psfrag{x16}[t][t]{\scriptsize 1}%
\psfrag{x17}[t][t]{\scriptsize 2}%
\psfrag{x18}[t][t]{\scriptsize 3}%
\psfrag{x19}[t][t]{\scriptsize -3}%
\psfrag{x20}[t][t]{\scriptsize -2}%
\psfrag{x21}[t][t]{\scriptsize -1}%
\psfrag{x22}[t][t]{\scriptsize 0}%
\psfrag{x23}[t][t]{\scriptsize 1}%
\psfrag{x24}[t][t]{\scriptsize 2}%
\psfrag{x25}[t][t]{\scriptsize 3}%
%
\psfrag{v01}[r][r]{\scriptsize 0}%
\psfrag{v02}[r][r]{\scriptsize 0.1}%
\psfrag{v03}[r][r]{\scriptsize 0.2}%
\psfrag{v04}[r][r]{\scriptsize 0.3}%
\psfrag{v05}[r][r]{\scriptsize 0.4}%
\psfrag{v06}[r][r]{\scriptsize 0.5}%
\psfrag{v07}[r][r]{\scriptsize 0.6}%
\psfrag{v08}[r][r]{\scriptsize 0.7}%
\psfrag{v09}[r][r]{\scriptsize 0.8}%
\psfrag{v10}[r][r]{\scriptsize 0.9}%
\psfrag{v11}[r][r]{\scriptsize 1}%
\psfrag{v12}[r][r]{\scriptsize 0}%
\psfrag{v13}[r][r]{\scriptsize 0.2}%
\psfrag{v14}[r][r]{\scriptsize 0.4}%
\psfrag{v15}[r][r]{\scriptsize 0.6}%
\psfrag{v16}[r][r]{\scriptsize 0.8}%
\psfrag{v17}[r][r]{\scriptsize 1}%
\psfrag{v18}[r][r]{\scriptsize 0}%
\psfrag{v19}[r][r]{\scriptsize 0.2}%
\psfrag{v20}[r][r]{\scriptsize 0.4}%
\psfrag{v21}[r][r]{\scriptsize 0.6}%
\psfrag{v22}[r][r]{\scriptsize 0.8}%
\psfrag{v23}[r][r]{\scriptsize 1}%
%
\setlength{\unitlength}{1cm}
\begin{picture}(8, 6)(0, 0)
\includegraphics[scale=0.55]{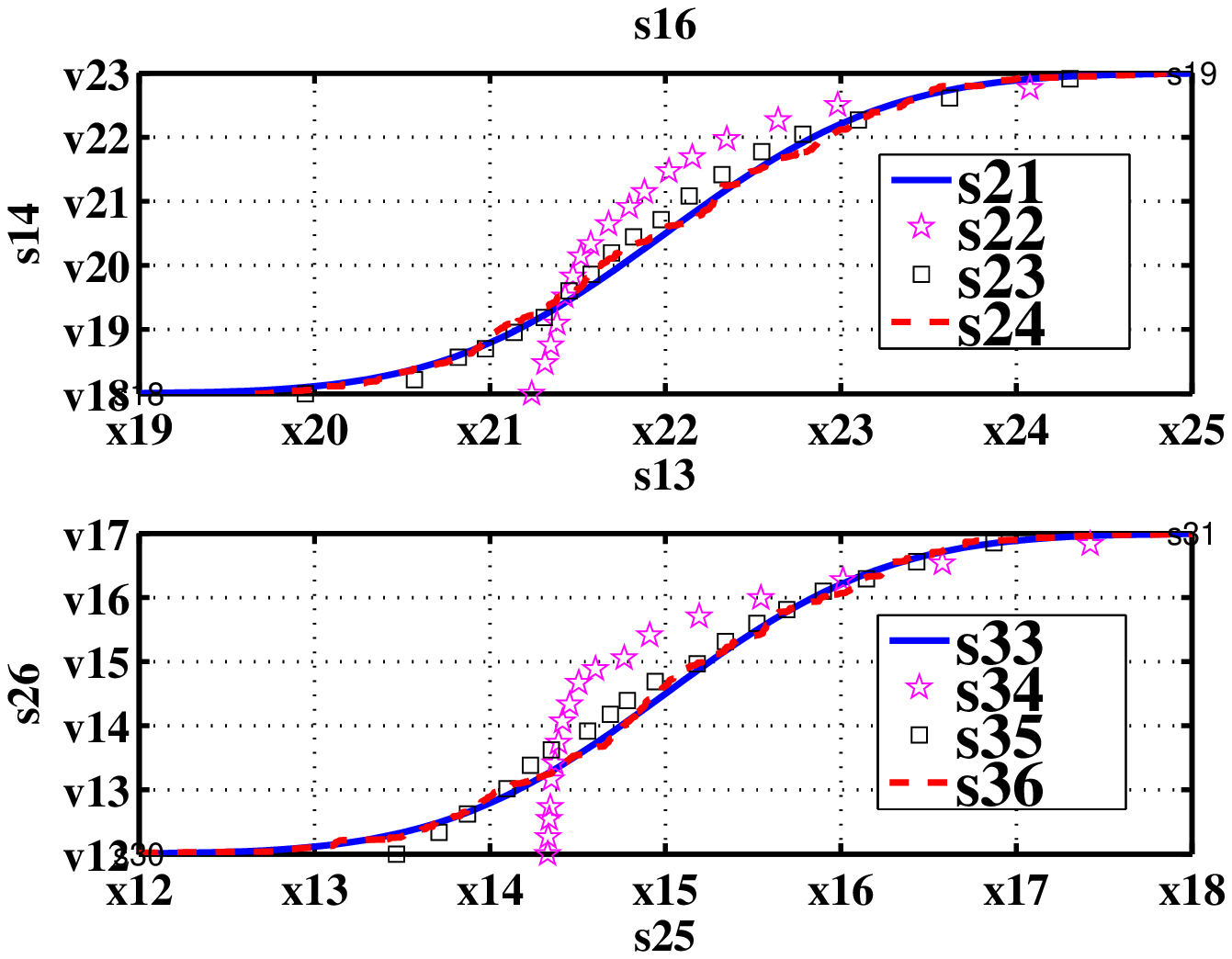}%
\put(-5.67, 4.35){\color[rgb]{0, 0, 0} \scriptsize \colorbox{white}{\fbox{$\alpha = 3$}}}
\put(-5.67, 1.54){\color[rgb]{0, 0, 0} \scriptsize \colorbox{white}{\fbox{$\alpha = 5$}}}
\end{picture}
\end{psfrags}%
%

\end{center}
\end{minipage}
\caption{Comparison of the simulated centered and normalized WMAI CDFs with the normal CDF for the path-loss functions $G_1(t) = \frac{1}{\paren{1+t}^\alpha}$ (lefthand side figures) and $G_2(t) = \frac{1}{1+t^\alpha}$ (righthand side figures).  The effect of fading is also illustrated in the bottom figures by assuming Nakagami-$m$ fading with $m$ parameter set to $1$.} \label{Fig: Simulated CDFs}
\end{figure*}

For the simulation study, we built a C-Simulator in order to perform Monte-Carlo simulations deriving WMAI CDFs numerically.  In Fig. \ref{Fig: Simulated CDFs}, the simulated WMAI distributions are plotted for $G_1(t)$ and $G_2(t)$ for two different choices of $\alpha$, {\em i.e.,} $3$ and $5$, under various values of $\lambda$.  Fading effects are also illustrated in the bottom figures by using the Nakagami-$m$ fading model with unit mean power gain.  $m$ parameter is set to $1$, which is the Rayleigh fading case.  Similar results continue to hold for other bounded path-loss models and different values of $\alpha$ greater than $2$.  For each different choice of the path-loss model, path-loss exponent and node intensity, we compute the interference power level at the origin for $10^4$ random node configurations in order to estimate the WMAI distributions. 

We focus on small to moderate values of $\lambda$ to illustrate the Gaussian convergence result predicted by Theorem \ref{Thm: Rates of Convergence - Stationary}, and to understand the effect of small values of $\lambda$ on the WMAI distributions.  As observed in Fig. \ref{Fig: Simulated CDFs}, the deviations between the normal distribution and the simulated WMAI distributions are prominent for sparse networks, {\em i.e.,} $\lambda = 0.1$.  On the other hand,  the match between the normal CDF and the simulated WMAI distributions is promising for small to moderate values of $\lambda$, {\em i.e.,} $\lambda=1$ and $\lambda=10$.  In the non-fading case, for example, the match between the simulated distributions and the normal distribution is {\em almost} perfect for both path-loss models and path-loss exponents when $\lambda$ is  around $10$ nodes per unit area.  Even when $\lambda$ is around $1$, it is still very good.  These observations in conjunction with Theorem \ref{Thm: Rates of Convergence - Stationary} illustrate the utility of the Gaussian approximation of the WMAI distributions for small and large values of $\lambda$.  They also indicate the potential to further tighten the upper and lower bounds that we derive in this paper. 

For communication environments with fading, the deviations between the simulated distributions and the normal distribution tend to be larger when compared to those for communication environments without fading.  This is also in accordance with our bounds in Theorem \ref{Thm: Rates of Convergence - Stationary}.  That is, the factor $m_{H^3}/\paren{m_{H^2}}^{\frac32}$ appearing in Theorem \ref{Thm: Rates of Convergence - Stationary} is always greater than $1$ due to Jensen's inequality, as mentioned above, which implies looser approximation bounds in the presence of random fading effects.  On the other hand, as $m$ increases, the fading process corrupting transmitted signals becomes more deterministic, and therefore our bounds become tighter, and we start to observe better matches between the simulated WMAI distributions and the normal CDF.           

When the effect of small and large values of $\alpha$ on the WMAI distributions is analyzed,  it is seen that the match between the simulated WMAI distributions and the normal distribution is slightly better for small values of $\alpha$ ({\em e.g.,} for $\lambda = 0.1$ and $G_2(t)$ without fading in Fig. \ref{Fig: Simulated CDFs}, the maximum deviation between the simulated WMAI distribution and the normal distribution is $0.11$ and $0.21$ for $\alpha = 3$ and $5$, respectively.).  This is an expected result when we compare the path-loss model dependent constants appearing in Theorem \ref{Thm: Rates of Convergence - Stationary}, which are also numerically computed for various values of $\alpha$ for both path-loss models in Table \ref{Table: 1} for a comparative illustration.  When we compare the effect of different path-loss models on the Gaussian approximation, we observe that the match between the simulated WMAI distributions and the normal distribution is slightly better for $G_2(t)$ ({\em e.g.,} when $\lambda = 0.1$ and $\alpha=5$ without fading in Fig.  \ref{Fig: Simulated CDFs}, the maximum deviation between the simulated WMAI distribution and the normal distribution is $0.3$ and $0.21$ for $G_1(t)$ and $G_2(t)$, respectively.).  This is also an expected result when we compare the path-loss model dependent constants appearing in Theorem \ref{Thm: Rates of Convergence - Stationary} ({\em i.e.,} see Table \ref{Table: 1} again).

\begin{table}[t]
\caption{Path-loss Dependent Gaussian Approximation Constant $\paren{\frac{\int_0^\infty G^3(t) t dt}{\paren{\int_0^\infty G^2(t) t dt}^\frac32}}$} \label{Table: 1}
\newcommand{\m}{\hphantom{$-$}}
\setlength{\heavyrulewidth}{0.1em}
\newcommand{\otoprule}{\midrule[\heavyrulewidth]}
\renewcommand{\tabcolsep}{0.65pc} 
\renewcommand{\arraystretch}{1.5} 
\centering
\begin{threeparttable}[b]
\begin{tabularx}{8cm}{c c c c}
\toprule
Path-loss model & \multicolumn{3}{c}{Path-loss exponent ($\alpha$)} \\ \cmidrule(l){1-4}
& $\alpha = 3$ & $\alpha = 4$ & $\alpha=5$ \\ \cmidrule{2-4}
\begin{minipage}[t]{3cm} \centering $\displaystyle G_1(t) = \frac{1}{\paren{1+t}^\alpha}$ \end{minipage}& $1.564$& $2.3838$& $3.1688$\\
\begin{minipage}[t]{3cm} \centering $\displaystyle G_2(t) = \frac{1}{1+t^\alpha}$ \end{minipage}& $1.0501$& $1.1972$& $1.2713$ \\ \bottomrule
\end{tabularx}
\end{threeparttable}
\end{table}  

\section{Performance Bounds: Outage Capacity and Sum Capacity} \label{Section: Applications}
In this section of the paper, we present further applications of our Gaussian approximation bounds above to derive various performance limits and metrics for wireless networks.  In particular, we will illustrate two potential applications of our results to bound the {\em single} link outage capacity and the {\em ergodic} sum capacity for spatial wireless networks in the current and next subsections. These capacity measures can be computed in an exact form for Rayleigh fading channels \cite{BBM06} but there  does not exist such closed form expressions for general fading channel models.  Similar to our analysis above, we will employ two different path-loss models $G_1(t) = \frac{1}{(1+t)^\alpha}$ and $G_2(t) = \frac{1}{1+t^\alpha}$ with various values of $\alpha$.  Similar conclusions continue to hold for other path-loss models.   

\subsection{Single Link Outage Capacity} \label{Subsection: Outage Capacity}
We start our analysis with outage capacity calculations.  To this end, we introduce more structure and modeling parameters into the communication set-up under consideration to express the problem of bounding the single link outage capacity as general as possible.     

\subsubsection{Problem Set-up}  We consider a planar wireless network that contains a test transmitter-receiver (TX-RX) pair whose separation is $d$ [unit distance] (usually in kilometers).  Without loss of generality, the test receiver node is assumed to be located at the origin $\vec{0}$.  The wireless communication link between the TX-RX pair is subject to fading, path-loss, background noise, and interference signals emanating from other transmitters in the network.  Hence, the maximum (short-term) data rate of reliable communication supported by the victim link, which is to be defined shortly below, fluctuates as a function of all such effects.  We assume that the TX encodes data at a rate $R$ nats per second per hertz, and the outage event occurs whenever the maximum rate supported by the victim link is below $R$.  This formulation is appropriate for {\em delay sensitive} traffic for which the communication delay requirement is smaller than the time scale of channel variations.     

Let $N_0$ be the power of the (complex) Gaussian background noise present at the RX.  We define the signal-to-noise-ratio ($\snr$) of the communication system under consideration as the ratio $\snr = \frac{P}{N_0}$, where $P$, as defined in Section \ref{Section: Network Model}, is the transmission power common to all transmitters.  With a slight change of notation, we will denote the total interference power at the RX as $I_\lambda(P)$ to simplify the further notation below.  Further, we will denote the cumulative distribution function of $I_\lambda(P)$ by $F_\lambda$.  The same assumptions in Section \ref{Section: Network Model} continue to hold for the spatial distribution of interfering transmitters here.   

Then, the signal-to-interference-plus-noise ratio ($\sinr$) at the RX is given as
\begin{eqnarray}
\sinr = \frac{\widetilde{H} G(d)}{\snr^{-1} + \frac{1}{\pg} I_\lambda(1)},
\end{eqnarray}
where $\pg \geq 1$ is the processing gain of the system, and $\widetilde{H}$ is the fading coefficient, with a finite mean value, for the link between the TX-RX pair.  Roughly speaking, $\pg = 1$ case represents a narrowband communication scenario, whereas $\pg > 1$ case signifies a broadband ({\em e.g.,} a CDMA network) setting.  For a given realization of channel states, the maximum rate of reliable communication supported by the victim link in nats per second per hertz is equal to $\log\paren{1+\sinr}$, and therefore the communication between the TX-RX pair is said to be in {\em outage} if $\log\paren{1+\sinr} < R$.  For a given target outage probability $\gamma \in [0, 1]$, the outage capacity for the victim link is defined as
\begin{eqnarray}
C_{\lambda, {\rm outage}}\paren{\gamma} = \sup\brparen{R > 0: \PR{\log\paren{1+\sinr} < R} \leq \gamma}, 
\end{eqnarray}
which is the maximum data rate supported by the victim link with outage probability not exceeding $\gamma$.  

\subsubsection{Theoretical and Simulation Results}  The next theorem provides the upper and lower bounds on $C_{\lambda, {\rm outage}}\paren{\gamma}$.
\begin{theorem}\label{Thm: Outage Capacity Bounds}
Let $\zeta\paren{h, R} = \frac{\paren{\frac{h G(d)}{\e{R}-1}-\snr^{-1}}\pg - \ES{I_\lambda(1)}}{\sqrt{\V{I_\lambda(1)}}}$.  For the communication scenario above, $C_{\lambda, {\rm outage}}\paren{\gamma}$ is upper and lower bounded as
\begin{eqnarray}
C_{\lambda, {\rm outage}}\paren{\gamma} \geq \sup\brparen{R > 0: 1 - \ES{\max\brparen{0, Q_{\lambda, {\rm outage}}^-\paren{\zeta\paren{\widetilde{H}, R}}} \I{\widetilde{H} \geq \frac{\snr^{-1}\paren{\e{R}-1}}{G(d)}}} \leq \gamma} 
\end{eqnarray}
and
\begin{eqnarray}
C_{\lambda, {\rm outage}}\paren{\gamma} \leq \sup\brparen{R > 0: 1 - \ES{\min\brparen{1, Q_{\lambda, {\rm outage}}^+\paren{\zeta\paren{\widetilde{H}, R}}} \I{\widetilde{H} \geq \frac{\snr^{-1}\paren{\e{R}-1}}{G(d)}}} \leq \gamma},
\end{eqnarray}
where $Q_{\lambda, {\rm outage}}^\pm(x) = \Psi(x) \pm \frac{c(x)}{\sqrt{\lambda}}$, and $\Psi(x)$ and $c(x)$ are as given in Theorem \ref{Thm: Rates of Convergence}.   
\end{theorem}
\begin{IEEEproof}
Let $\tilde{q}$ be the probability density function for $\widetilde{H}$, possibly a different density function than $q$.  By straightforward manipulations, we can express the outage probability $\PRP{\mbox{Outage}} = \PR{\log\paren{1+\sinr} < R}$ as
\begin{eqnarray}
\PRP{\mbox{Outage}} &=& \PR{\sinr < \e{R}-1} \nonumber \\
&=& \int_0^\infty \PR{\paren{\frac{h G(d)}{\e{R}-1} - \snr^{-1}}\pg < I_\lambda(1)} \tilde{q}(h) \diff h \nonumber \\
&=& 1 - \int_{\frac{\snr^{-1}\paren{\e{R} - 1}}{G(d)}}^\infty \PR{I_\lambda(1) \leq \paren{\frac{h G(d)}{\e{R}-1} - \snr^{-1}}\pg}\tilde{q}(h) \diff h, \nonumber
\end{eqnarray}
where the last equality follows from the fact that $I_\lambda(1)$ is a positive random variable, and we have $\paren{\frac{h G(d)}{\e{R}-1} - \snr^{-1}}\pg < 0$ if and only if $h < \frac{\snr^{-1}\paren{\e{R} - 1}}{G(d)}$.  By using Theorem \ref{Thm: Rates of Convergence} and the natural bounds $0$ and $1$ on the probability, we can upper and lower bound $\PRP{\mbox{Outage}}$ as 
\begin{eqnarray}
\PRP{\mbox{Outage}} &\leq& 1 - \int_{\frac{\snr^{-1}\paren{\e{R} - 1}}{G(d)}}^\infty \max\brparen{0, Q_{\lambda, {\rm outage}}^-\paren{\zeta\paren{h, R}}} \tilde{q}(h) \diff h \nonumber \\
&=& 1 - \ES{\max\brparen{0, Q_{\lambda, {\rm outage}}^-\paren{\zeta\paren{\widetilde{H}, R}}} \I{\widetilde{H} \geq \frac{\snr^{-1}\paren{\e{R}-1}}{G(d)}}} \nonumber
\end{eqnarray}
and
\begin{eqnarray}
\PRP{\mbox{Outage}} &\geq& 1 - \int_{\frac{\snr^{-1}\paren{\e{R} - 1}}{G(d)}}^\infty \min\brparen{1, Q_{\lambda, {\rm outage}}^+\paren{\zeta\paren{h, R}}} \tilde{q}(h) \diff h \nonumber \\
&=& 1 - \ES{\min\brparen{1, Q_{\lambda, {\rm outage}}^+\paren{\zeta\paren{\widetilde{H}, R}}} \I{\widetilde{H} \geq \frac{\snr^{-1}\paren{\e{R}-1}}{G(d)}}}. \nonumber
\end{eqnarray}
The proof is completed by observing that the upper (lower) bound on the outage probability crosses the target outage probability $\gamma$ earlier (later) than $\PRP{\mbox{Outage}}$ as $R$ increases.        
\end{IEEEproof}
\begin{figure*}[!t]
\begin{minipage}[t]{\textwidth}
\begin{center}
%
%
\begin{psfrags}%
\psfragscanon%
%
\psfrag{s09}[t][t]{\color[rgb]{0,0,0}\setlength{\tabcolsep}{0pt}\scriptsize \begin{tabular}{c}Transmitter Intensity ($\lambda$)\end{tabular}}%
\psfrag{s10}[b][b]{\color[rgb]{0,0,0}\setlength{\tabcolsep}{0pt}\scriptsize \begin{tabular}{c}Outage Capacity and Bounds \\ (Nats per Second per Hertz) \end{tabular}}%
\psfrag{s12}[b][b]{\color[rgb]{0,0,0}\setlength{\tabcolsep}{0pt}\scriptsize \begin{tabular}{c}$G_1(t) = \frac{1}{(1+t)^{4}}$\end{tabular}}%
\psfrag{s14}[][]{\color[rgb]{0,0,0}\setlength{\tabcolsep}{0pt}\scriptsize \begin{tabular}{c} $ $ \end{tabular}}%
\psfrag{s15}[][]{\color[rgb]{0,0,0}\setlength{\tabcolsep}{0pt}\scriptsize \begin{tabular}{c} $ $ \end{tabular}}%
\psfrag{s16}[l][l]{\color[rgb]{0,0,0}\scriptsize Lower Bound}%
\psfrag{s17}[l][l]{\color[rgb]{0,0,0}\scriptsize Simulated Rate}%
\psfrag{s18}[l][l]{\color[rgb]{0,0,0}\scriptsize Upper Bound}%
\psfrag{s19}[l][l]{\color[rgb]{0,0,0}\scriptsize Lower Bound}%
%
\psfrag{x01}[t][t]{\scriptsize 0}%
\psfrag{x02}[t][t]{\scriptsize 0.1}%
\psfrag{x03}[t][t]{\scriptsize 0.2}%
\psfrag{x04}[t][t]{\scriptsize 0.3}%
\psfrag{x05}[t][t]{\scriptsize 0.4}%
\psfrag{x06}[t][t]{\scriptsize 0.5}%
\psfrag{x07}[t][t]{\scriptsize 0.6}%
\psfrag{x08}[t][t]{\scriptsize 0.7}%
\psfrag{x09}[t][t]{\scriptsize 0.8}%
\psfrag{x10}[t][t]{\scriptsize 0.9}%
\psfrag{x11}[t][t]{\scriptsize 1}%
\psfrag{x12}[t][t]{\scriptsize 1}%
\psfrag{x13}[t][t]{\scriptsize 2}%
\psfrag{x14}[t][t]{\scriptsize 3}%
\psfrag{x15}[t][t]{\scriptsize 4}%
\psfrag{x16}[t][t]{\scriptsize 5}%
\psfrag{x17}[t][t]{\scriptsize 6}%
\psfrag{x18}[t][t]{\scriptsize 7}%
\psfrag{x19}[t][t]{\scriptsize 8}%
\psfrag{x20}[t][t]{\scriptsize 9}%
\psfrag{x21}[t][t]{\scriptsize 10}%
%
\psfrag{v01}[r][r]{\scriptsize 0}%
\psfrag{v02}[r][r]{\scriptsize 0.1}%
\psfrag{v03}[r][r]{\scriptsize 0.2}%
\psfrag{v04}[r][r]{\scriptsize 0.3}%
\psfrag{v05}[r][r]{\scriptsize 0.4}%
\psfrag{v06}[r][r]{\scriptsize 0.5}%
\psfrag{v07}[r][r]{\scriptsize 0.6}%
\psfrag{v08}[r][r]{\scriptsize 0.7}%
\psfrag{v09}[r][r]{\scriptsize 0.8}%
\psfrag{v10}[r][r]{\scriptsize 0.9}%
\psfrag{v11}[r][r]{\scriptsize 1}%
\psfrag{v12}[r][r]{\scriptsize 0.2}%
\psfrag{v13}[r][r]{\scriptsize 0.4}%
\psfrag{v14}[r][r]{\scriptsize 0.6}%
\psfrag{v15}[r][r]{\scriptsize 0.8}%
\psfrag{v16}[r][r]{\scriptsize 1}%
\psfrag{v17}[r][r]{\scriptsize 1.2}%
%
\includegraphics[scale=0.55]{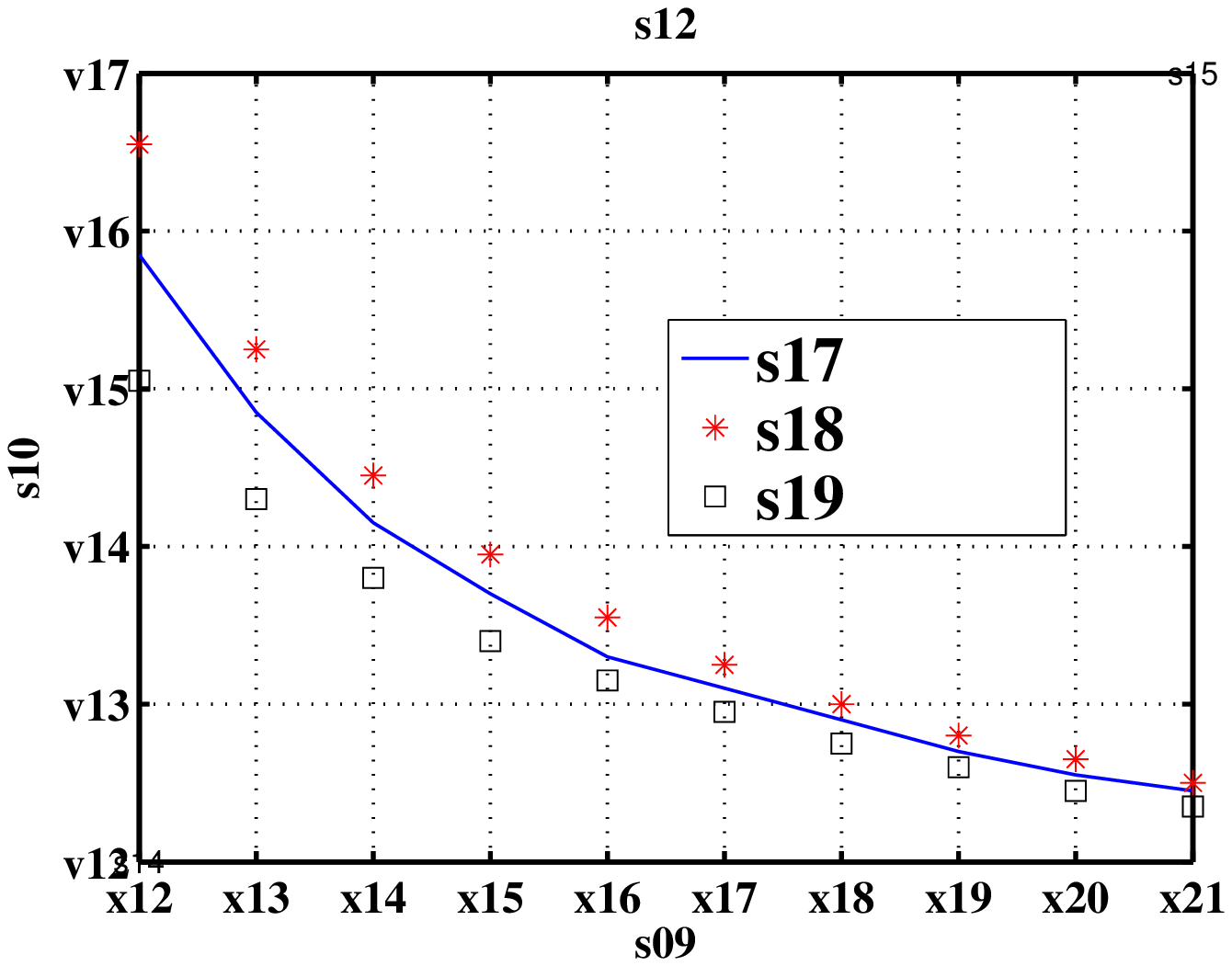}%
\end{psfrags}%
%

%
\hspace{\fill}
%
%
%
\begin{psfrags}%
\psfragscanon%
%
\psfrag{s09}[t][t]{\color[rgb]{0,0,0}\setlength{\tabcolsep}{0pt}\scriptsize \begin{tabular}{c} Transmitter Intensity ($\lambda$) \end{tabular}}%
\psfrag{s10}[b][b]{\color[rgb]{0,0,0}\setlength{\tabcolsep}{0pt}\scriptsize \begin{tabular}{c} Outage Capacity and Bounds \\ (Nats per Second per Hertz) \end{tabular}}%
\psfrag{s12}[b][b]{\color[rgb]{0,0,0}\setlength{\tabcolsep}{0pt}\scriptsize \begin{tabular}{c} $G_1(t) = \frac{1}{(1+t)^{4}}$ \end{tabular}}%
\psfrag{s14}[][]{\color[rgb]{0,0,0}\setlength{\tabcolsep}{0pt}\scriptsize \begin{tabular}{c} $ $ \end{tabular}}%
\psfrag{s15}[][]{\color[rgb]{0,0,0}\setlength{\tabcolsep}{0pt}\scriptsize \begin{tabular}{c} $ $ \end{tabular}}%
\psfrag{s16}[l][l]{\color[rgb]{0,0,0}\scriptsize Lower Bound}%
\psfrag{s17}[l][l]{\color[rgb]{0,0,0}\scriptsize Simulated Rate}%
\psfrag{s18}[l][l]{\color[rgb]{0,0,0}\scriptsize Upper Bound}%
\psfrag{s19}[l][l]{\color[rgb]{0,0,0}\scriptsize Lower Bound}%
%
\psfrag{x01}[t][t]{\scriptsize 0}%
\psfrag{x02}[t][t]{\scriptsize 0.1}%
\psfrag{x03}[t][t]{\scriptsize 0.2}%
\psfrag{x04}[t][t]{\scriptsize 0.3}%
\psfrag{x05}[t][t]{\scriptsize 0.4}%
\psfrag{x06}[t][t]{\scriptsize 0.5}%
\psfrag{x07}[t][t]{\scriptsize 0.6}%
\psfrag{x08}[t][t]{\scriptsize 0.7}%
\psfrag{x09}[t][t]{\scriptsize 0.8}%
\psfrag{x10}[t][t]{\scriptsize 0.9}%
\psfrag{x11}[t][t]{\scriptsize 1}%
\psfrag{x12}[t][t]{\scriptsize 100}%
\psfrag{x13}[t][t]{\scriptsize 20}%
\psfrag{x14}[t][t]{\scriptsize 30}%
\psfrag{x15}[t][t]{\scriptsize 40}%
\psfrag{x16}[t][t]{\scriptsize 50}%
\psfrag{x17}[t][t]{\scriptsize 60}%
\psfrag{x18}[t][t]{\scriptsize 70}%
\psfrag{x19}[t][t]{\scriptsize 80}%
\psfrag{x20}[t][t]{\scriptsize 90}%
%
\psfrag{v01}[r][r]{\scriptsize 0}%
\psfrag{v02}[r][r]{\scriptsize 0.1}%
\psfrag{v03}[r][r]{\scriptsize 0.2}%
\psfrag{v04}[r][r]{\scriptsize 0.3}%
\psfrag{v05}[r][r]{\scriptsize 0.4}%
\psfrag{v06}[r][r]{\scriptsize 0.5}%
\psfrag{v07}[r][r]{\scriptsize 0.6}%
\psfrag{v08}[r][r]{\scriptsize 0.7}%
\psfrag{v09}[r][r]{\scriptsize 0.8}%
\psfrag{v10}[r][r]{\scriptsize 0.9}%
\psfrag{v11}[r][r]{\scriptsize 1}%
\psfrag{v12}[r][r]{\scriptsize 0}%
\psfrag{v13}[r][r]{\scriptsize 0.05}%
\psfrag{v14}[r][r]{\scriptsize 0.1}%
\psfrag{v15}[r][r]{\scriptsize 0.15}%
\psfrag{v16}[r][r]{\scriptsize 0.2}%
%
\includegraphics[scale=0.55]{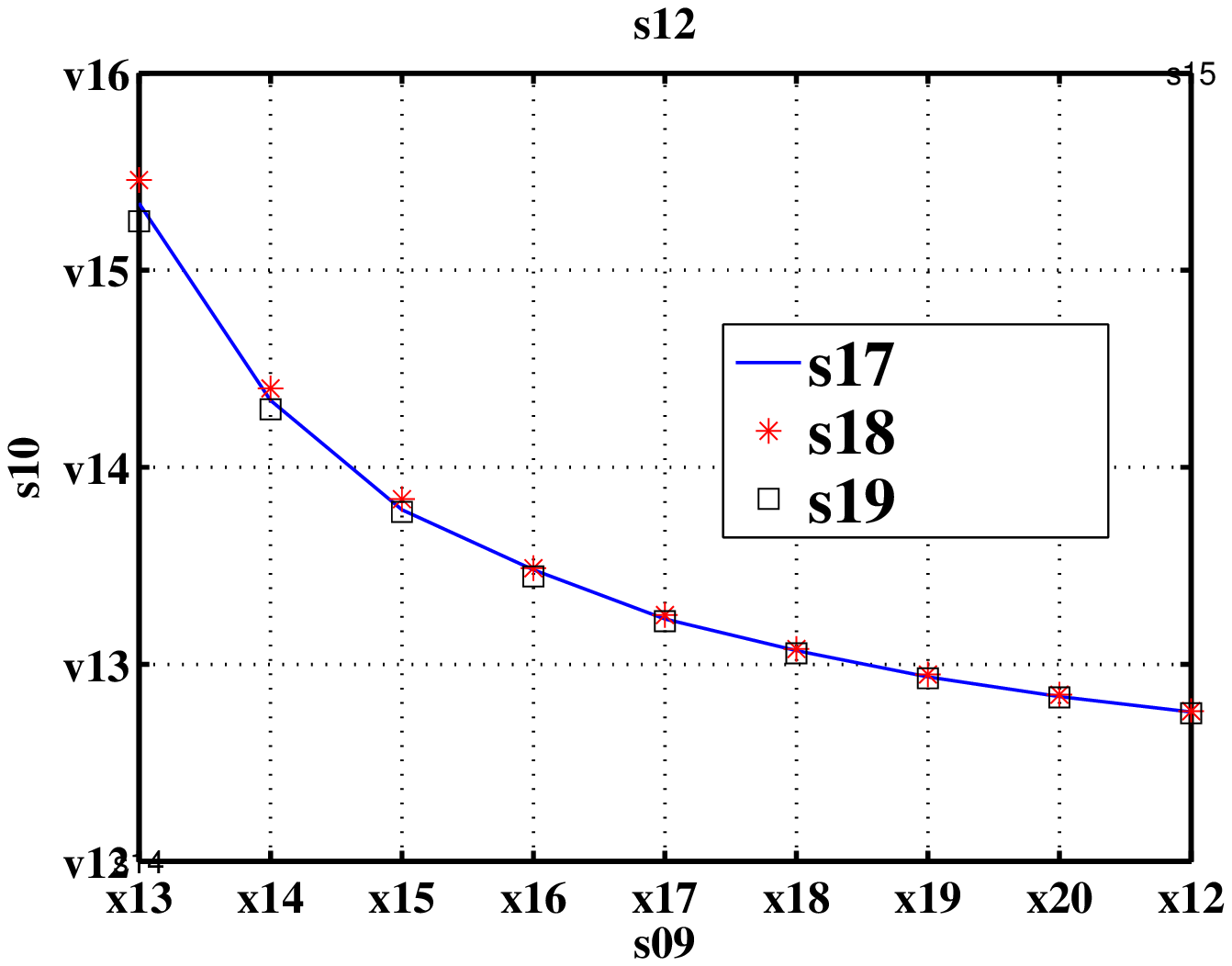}%
\end{psfrags}%
%

\end{center}
\end{minipage} 
\\ 
\begin{minipage}[t]{\textwidth}
\vspace{-2mm} 
\begin{center}
%
%
\begin{psfrags}%
\psfragscanon%
%
\psfrag{s09}[t][t]{\color[rgb]{0,0,0}\setlength{\tabcolsep}{0pt}\scriptsize \begin{tabular}{c}Transmitter Intensity ($\lambda$)\end{tabular}}%
\psfrag{s10}[b][b]{\color[rgb]{0,0,0}\setlength{\tabcolsep}{0pt}\scriptsize \begin{tabular}{c} Outage Capacity and Bounds \\ (Nats per Second per Hertz) \end{tabular}}%
\psfrag{s12}[b][b]{\color[rgb]{0,0,0}\setlength{\tabcolsep}{0pt}\scriptsize \begin{tabular}{c}$G_2(t) = \frac{1}{1+t^{4}}$\end{tabular}}%
\psfrag{s14}[][]{\color[rgb]{0,0,0}\setlength{\tabcolsep}{0pt}\scriptsize \begin{tabular}{c} $ $ \end{tabular}}%
\psfrag{s15}[][]{\color[rgb]{0,0,0}\setlength{\tabcolsep}{0pt}\scriptsize \begin{tabular}{c} $ $ \end{tabular}}%
\psfrag{s16}[l][l]{\color[rgb]{0,0,0}\scriptsize Lower Bound}%
\psfrag{s17}[l][l]{\color[rgb]{0,0,0}\scriptsize Simulated Rate}%
\psfrag{s18}[l][l]{\color[rgb]{0,0,0}\scriptsize Upper Bound}%
\psfrag{s19}[l][l]{\color[rgb]{0,0,0}\scriptsize Lower Bound}%
%
\psfrag{x01}[t][t]{\scriptsize 0}%
\psfrag{x02}[t][t]{\scriptsize 0.1}%
\psfrag{x03}[t][t]{\scriptsize 0.2}%
\psfrag{x04}[t][t]{\scriptsize 0.3}%
\psfrag{x05}[t][t]{\scriptsize 0.4}%
\psfrag{x06}[t][t]{\scriptsize 0.5}%
\psfrag{x07}[t][t]{\scriptsize 0.6}%
\psfrag{x08}[t][t]{\scriptsize 0.7}%
\psfrag{x09}[t][t]{\scriptsize 0.8}%
\psfrag{x10}[t][t]{\scriptsize 0.9}%
\psfrag{x11}[t][t]{\scriptsize 1}%
\psfrag{x12}[t][t]{\scriptsize 1}%
\psfrag{x13}[t][t]{\scriptsize 2}%
\psfrag{x14}[t][t]{\scriptsize 3}%
\psfrag{x15}[t][t]{\scriptsize 4}%
\psfrag{x16}[t][t]{\scriptsize 5}%
\psfrag{x17}[t][t]{\scriptsize 6}%
\psfrag{x18}[t][t]{\scriptsize 7}%
\psfrag{x19}[t][t]{\scriptsize 8}%
\psfrag{x20}[t][t]{\scriptsize 9}%
\psfrag{x21}[t][t]{\scriptsize 10}%
%
\psfrag{v01}[r][r]{\scriptsize 0}%
\psfrag{v02}[r][r]{\scriptsize 0.1}%
\psfrag{v03}[r][r]{\scriptsize 0.2}%
\psfrag{v04}[r][r]{\scriptsize 0.3}%
\psfrag{v05}[r][r]{\scriptsize 0.4}%
\psfrag{v06}[r][r]{\scriptsize 0.5}%
\psfrag{v07}[r][r]{\scriptsize 0.6}%
\psfrag{v08}[r][r]{\scriptsize 0.7}%
\psfrag{v09}[r][r]{\scriptsize 0.8}%
\psfrag{v10}[r][r]{\scriptsize 0.9}%
\psfrag{v11}[r][r]{\scriptsize 1}%
\psfrag{v12}[r][r]{\scriptsize 0}%
\psfrag{v13}[r][r]{\scriptsize 0.5}%
\psfrag{v14}[r][r]{\scriptsize 1}%
\psfrag{v15}[r][r]{\scriptsize 1.5}%
\psfrag{v16}[r][r]{\scriptsize 2}%
%
\includegraphics[scale=0.55]{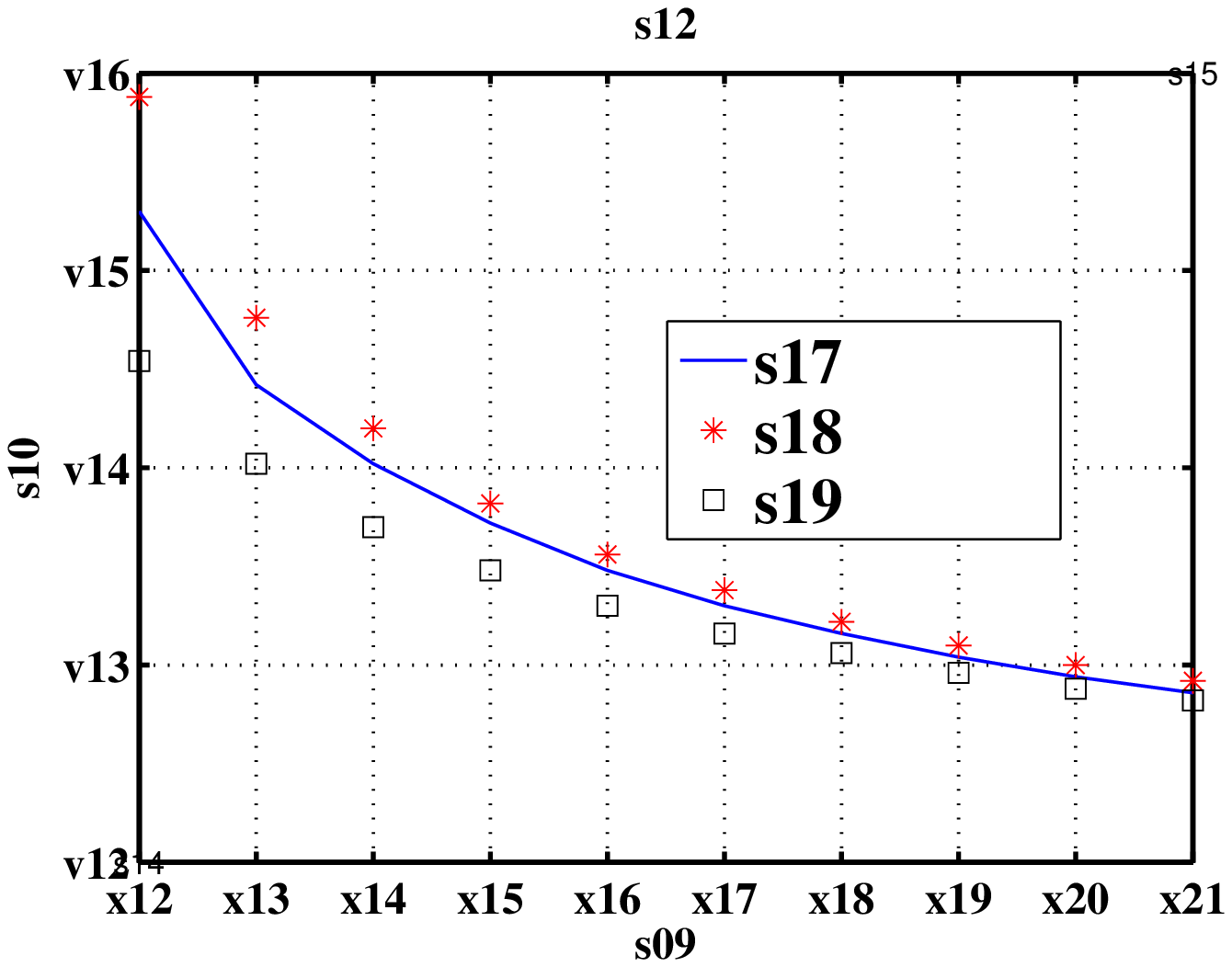}%
\end{psfrags}%
%

%
\hfill
%
%
%
\begin{psfrags}%
\psfragscanon%
%
\psfrag{s09}[t][t]{\color[rgb]{0,0,0}\setlength{\tabcolsep}{0pt}\scriptsize \begin{tabular}{c} Transmitter Intensity ($\lambda$) \end{tabular}}%
\psfrag{s10}[b][b]{\color[rgb]{0,0,0}\setlength{\tabcolsep}{0pt}\scriptsize \begin{tabular}{c} Outage Capacity and Bounds \\ (Nats per Second per Hertz) \end{tabular}}%
\psfrag{s12}[b][b]{\color[rgb]{0,0,0}\setlength{\tabcolsep}{0pt}\scriptsize \begin{tabular}{c} $G_2(t) = \frac{1}{1+t^{4}}$ \end{tabular}}%
\psfrag{s14}[][]{\color[rgb]{0,0,0}\setlength{\tabcolsep}{0pt}\scriptsize \begin{tabular}{c} $ $ \end{tabular}}%
\psfrag{s15}[][]{\color[rgb]{0,0,0}\setlength{\tabcolsep}{0pt}\scriptsize \begin{tabular}{c} $ $ \end{tabular}}%
\psfrag{s16}[l][l]{\color[rgb]{0,0,0}\scriptsize Lower Bound}%
\psfrag{s17}[l][l]{\color[rgb]{0,0,0}\scriptsize Simulated Rate}%
\psfrag{s18}[l][l]{\color[rgb]{0,0,0}\scriptsize Upper Bound}%
\psfrag{s19}[l][l]{\color[rgb]{0,0,0}\scriptsize Lower Bound}%
%
\psfrag{x01}[t][t]{\scriptsize 0}%
\psfrag{x02}[t][t]{\scriptsize 0.1}%
\psfrag{x03}[t][t]{\scriptsize 0.2}%
\psfrag{x04}[t][t]{\scriptsize 0.3}%
\psfrag{x05}[t][t]{\scriptsize 0.4}%
\psfrag{x06}[t][t]{\scriptsize 0.5}%
\psfrag{x07}[t][t]{\scriptsize 0.6}%
\psfrag{x08}[t][t]{\scriptsize 0.7}%
\psfrag{x09}[t][t]{\scriptsize 0.8}%
\psfrag{x10}[t][t]{\scriptsize 0.9}%
\psfrag{x11}[t][t]{\scriptsize 1}%
\psfrag{x12}[t][t]{\scriptsize 100}%
\psfrag{x13}[t][t]{\scriptsize 20}%
\psfrag{x14}[t][t]{\scriptsize 30}%
\psfrag{x15}[t][t]{\scriptsize 40}%
\psfrag{x16}[t][t]{\scriptsize 50}%
\psfrag{x17}[t][t]{\scriptsize 60}%
\psfrag{x18}[t][t]{\scriptsize 70}%
\psfrag{x19}[t][t]{\scriptsize 80}%
\psfrag{x20}[t][t]{\scriptsize 90}%
%
\psfrag{v01}[r][r]{\scriptsize 0}%
\psfrag{v02}[r][r]{\scriptsize 0.1}%
\psfrag{v03}[r][r]{\scriptsize 0.2}%
\psfrag{v04}[r][r]{\scriptsize 0.3}%
\psfrag{v05}[r][r]{\scriptsize 0.4}%
\psfrag{v06}[r][r]{\scriptsize 0.5}%
\psfrag{v07}[r][r]{\scriptsize 0.6}%
\psfrag{v08}[r][r]{\scriptsize 0.7}%
\psfrag{v09}[r][r]{\scriptsize 0.8}%
\psfrag{v10}[r][r]{\scriptsize 0.9}%
\psfrag{v11}[r][r]{\scriptsize 1}%
\psfrag{v12}[r][r]{\scriptsize 0.05}%
\psfrag{v13}[r][r]{\scriptsize 0.1}%
\psfrag{v14}[r][r]{\scriptsize 0.15}%
\psfrag{v15}[r][r]{\scriptsize 0.2}%
\psfrag{v16}[r][r]{\scriptsize 0.25}%
\psfrag{v17}[r][r]{\scriptsize 0.3}%
%
\includegraphics[scale=0.55]{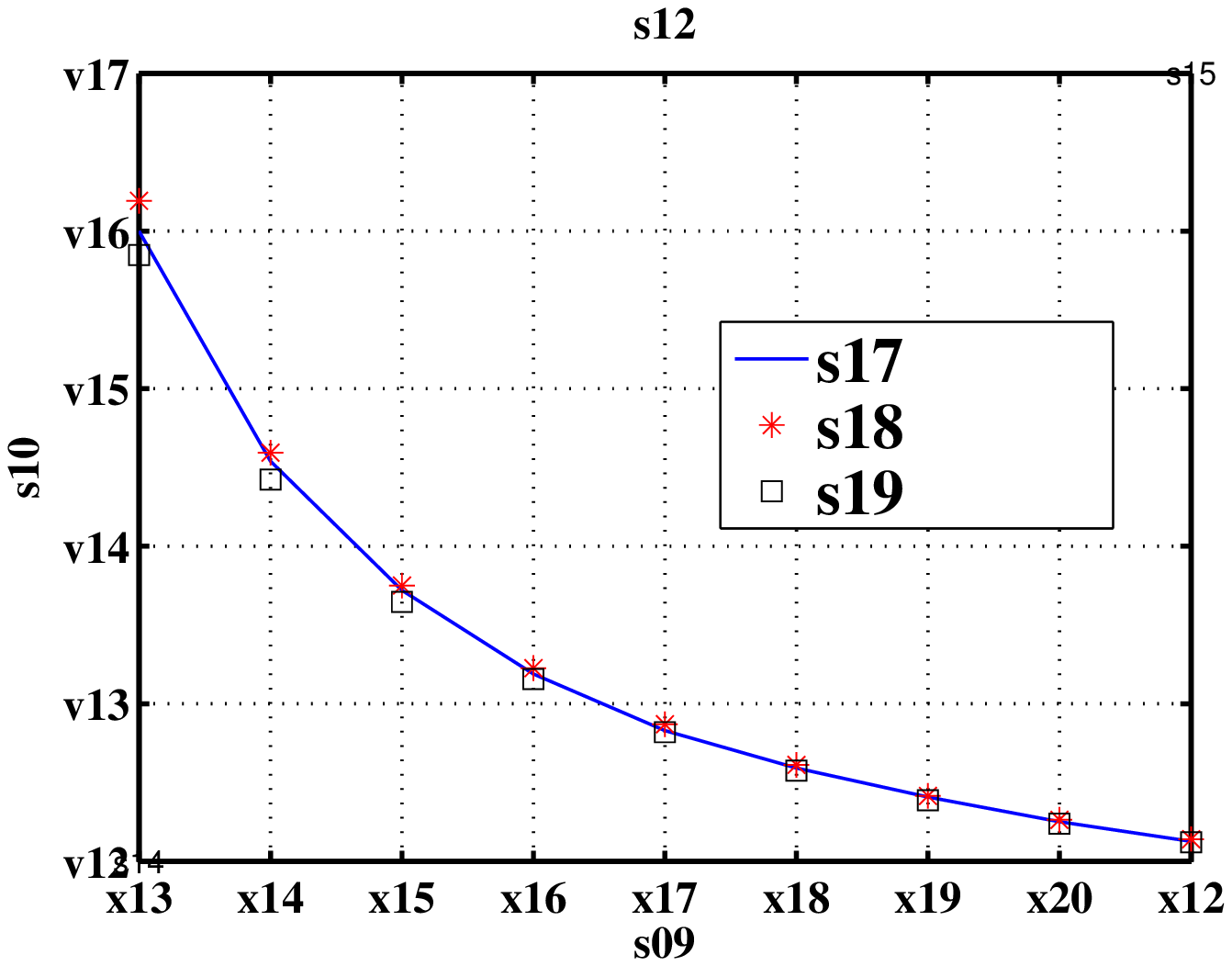}%
\end{psfrags}%
%

\end{center}
\end{minipage}
\caption{Upper and lower bounds on $C_{\lambda, {\rm outage}}\paren{\gamma}$ for the path-loss functions $G_1(t) = \frac{1}{\paren{1+t}^\alpha}$ (top figures) and $G_2(t) = \frac{1}{1+t^\alpha}$ (bottom figures).  Transmitters are distributed over $\R^2 - \mathcal{B}\paren{\vec{0}, \eta d}$ according to a stationary PPP of intensity $\lambda$.  All fading effects are modeled by using the Nakagami-$m$ fading model with $m$ parameter set to $5$.  For this model of spatial distribution of interfering transmitters, $c(x)$ is given as $c(x) = \frac{1}{\sqrt{2 \pi}}\frac{m_{H^3}}{\paren{m_{H^2}}^\frac32}\frac{\int_{\eta d}^\infty G^3(t) t \diff t}{\paren{\int_{\eta d}^\infty G^2(t) t \diff t}^\frac32} \min\paren{0.4785, \frac{31.935}{1+|x|^3}}$.  ($\snr = 20$ [dB], $\alpha = 4, d = 1, \pg = 100$ and $\eta = 0.5$.)} \label{Fig: Outage Capacity Bounds}
\end{figure*}

We plot $C_{\lambda, {\rm outage}}\paren{\gamma}$ and the corresponding bounds in Theorem \ref{Thm: Outage Capacity Bounds} as a function of $\lambda$ in Fig. \ref{Fig: Outage Capacity Bounds} for both path-loss functions.  For this simulation study, we assume that interfering transmitters are {\em uniformly} distributed over $\R^2 - \mathcal{B}\paren{\vec{0}, \eta d}$ according to a stationary PPP of intensity $\lambda$ [nodes per unit area], where $\eta \in [0, \infty)$.  Here, $\mathcal{B}\paren{\vec{0}, \eta d}$ can be interpreted as the {\em exclusion zone} around the RX in which no interfering transmitters are allowed, perhaps due to physical limitations ({\em i.e.,} small values of $\eta$), or the operation of the MAC layer ({\em i.e.,} CSMA-CA protocol), or the interference cancellation property \cite{WAYV07}.  $\eta$ is just a modeling parameter that allows us to control the radius of the exclusion zone, and the familiar stationary PPP model without any such holes can be recovered by setting it to zero.  

We set $\snr$ to $20$ [dB], and the path-loss exponent $\alpha$ to $4$.  All fading effects are modeled by means of the Nakagami-$m$ fading with unit mean power gain and $m=5$.  Similar observations continue to hold for other parameter selections.  We observe that our bounds closely approximate $C_{\lambda, {\rm outage}}\paren{\gamma}$ within {\em one} nats per second per hertz for moderate to high values of $\lambda$. Especially, for high values of $\lambda$, our bounds become very tight, and they {\em almost} coincide with $C_{\lambda, {\rm outage}}\paren{\gamma}$.  Considering the simulation and analytical results presented in Subsection \ref{Section: Approximation Bounds for Stationary PPPs}, this is an expected result since our Gaussian approximation for the interference power at the RX becomes more accurate in the dense network limit. 

\subsubsection{Scaling Behavior}  Our results could also provide a high level perspective about the detrimental effects of the network interference on the single link outage capacity.  With a slight abuse of notation, let $\gamma\paren{h, R} = 1 - F_\lambda\paren{\pg\paren{\frac{h G(d)}{\e{R}-1}-\snr^{-1}}}$.  $\gamma\paren{h, R}$ is the conditional probability of outage given $\widetilde{H} = h$ when the TX encodes data at rate $R$.  For simplicity, we will assume that the TX can track variations in $\widetilde{H}$, and adjust the rate of communication $R_\lambda\paren{h}$ as a function of observed values $h$ of $\widetilde{H}$.  For a given target outage probability $\gamma>0$, a reasonable rate selection policy is to choose $R_\lambda\paren{h}$ such that the outage probability is always equal to $\gamma$ for all $h$.\footnote{Optimum rate selection policy maximizing $\ES{R_\lambda\paren{\widetilde{H}}}$ subject to $\ES{\gamma\paren{\widetilde{H}, R_\lambda\paren{\widetilde{H}}}} \leq \gamma$ is an infinite dimensional, possibly non-convex depending on $F_\lambda$, functional optimization problem whose solution is out of the scope of this paper.}  Then, for any $\epsilon>0$, we can find $\bar{\lambda}$ such that 
\begin{eqnarray}
1-\gamma - \epsilon \leq \Psi\paren{\frac{\pg\paren{\frac{h G(d)}{\e{R_\lambda(h)} - 1} - \snr^{-1}} - \ES{I_\lambda(1)}}{\sqrt{\V{I_\lambda(1)}}}} \leq 1 - \gamma + \epsilon \nonumber 
\end{eqnarray}
for all $h>0$ and $\lambda \geq \bar{\lambda}$ since the interference distribution can be uniformly approximated by the normal distribution.  This implies 
$
\log\paren{1 + \frac{h G(d)}{\snr^{-1} + \frac{1}{\pg}\paren{\ES{I_\lambda(1)} + \sqrt{\V{I_\lambda(1)}} \Psi^{-1}\paren{1-\gamma+\epsilon}}}} \leq R_\lambda(h) \leq \log\paren{1 + \frac{h G(d)}{\snr^{-1} + \frac{1}{\pg}\paren{\ES{I_\lambda(1)} + \sqrt{\V{I_\lambda(1)}} \Psi^{-1}\paren{1-\gamma - \epsilon}}}}$, uniformly in $h$ for all $\lambda$ large enough.  Hence, we conclude that $\ES{R_\lambda\paren{\widetilde{H}}} = \TO{\frac{1}{\lambda}}$ as $\lambda \ra \infty$.  Similar but slightly more involved arguments, which are sketched in Appendix \ref{Appendix: Outage Probability Scaling}, also show that $C_{\lambda, {\rm outage}}\paren{\gamma} = \TO{\frac{1}{\lambda}}$ as $\lambda \ra \infty$.  Indeed, this is exactly the behavior observed in Fig. \ref{Fig: Outage Capacity Bounds}.  For example, a fivefold increase in $\lambda$ from $20$ to $100$ nodes per unit area results in a fivefold decrease in $C_{\lambda, {\rm outage}}\paren{\gamma}$ from $0.25$ to $0.05$ nats per second per hertz in the bottom righthand side performance figure in Fig. \ref{Fig: Outage Capacity Bounds}.

\subsection{Sum Capacity for Spatial Multiple Access Networks} \label{Subsection: Network Capacity}
We now illustrate another application of our results to bound the capacity of a {\em spatial} multiple access network in which transmitters are distributed according to a PPP with intensity parameter $\lambda$, and they all transmit to a common base station (BS) located at $\vec{0}$. The same assumptions in Section \ref{Section: Network Model} continue to hold for the spatial distribution of transmitters here.  This set-up requires us to interpret the test receiver node above as the common BS, and $I_\lambda(P)$ as the useful signal power for information flow from transmitters to the BS.  
Using the same notation above, the ergodic sum capacity of the network in nats per second per hertz is equal to
\begin{eqnarray}
C_{\lambda, {\rm sum}}\paren{\snr} &=& \ES{\log\paren{1+\snr \sum_{k \geq 1} H_k G\paren{\left\|\vec{X}_k\right\|_2}}} \nonumber \\ 
&=& \ES{\log\paren{1+I_\lambda\paren{\snr}}}, \label{Eqn: Capacity}
\end{eqnarray}
which is achievable by using complex Gaussian codebooks and successive interference cancellation receiver \cite{Tse05}.  Implicit in this formulation is that the communication delay requirement is much longer than the time scale of channel variations ({\em i.e.,} delay insensitive traffic) so that the BS can average over the fluctuations in the channel to achieve the communication rates in \eqref{Eqn: Capacity}.  
The next theorem provides the upper and lower bounds on $C_{\lambda, {\rm sum}}\paren{\snr}$. 
\begin{theorem} \label{Thm: Sum Capacity Bounds}
For the communication scenario above, $C_{\lambda, {\rm sum}}\paren{\snr}$ is upper and lower bounded as
\begin{eqnarray}
C_{\lambda, {\rm sum}}\paren{\snr} \geq \int_0^\infty 1 - \min\brparen{1, Q^+_\lambda\paren{\frac{\e{x}-1-\ES{I_\lambda\paren{\snr}}}{\sqrt{\V{I_\lambda\paren{\snr}}}}}} \diff x \nonumber
\end{eqnarray}
and
\begin{eqnarray}
C_{\lambda, {\rm sum}}\paren{\snr} \leq \int_0^\infty 1 - \max\brparen{0, Q^-_\lambda\paren{\frac{\e{x}-1-\ES{I_\lambda\paren{\snr}}}{\sqrt{\V{I_\lambda\paren{\snr}}}}}} \diff x, \nonumber 
\end{eqnarray}
where $Q_\lambda^\pm(x) = \Psi(x) \pm \frac{c(x)}{\sqrt{\lambda}}$, and $\Psi(x)$ and $c(x)$ are as given in Theorem \ref{Thm: Rates of Convergence}.       
\end{theorem}
\begin{IEEEproof}
Since $\log\paren{1+I_\lambda\paren{\snr}}$ is a positive random variable,  $C_{\lambda, {\rm sum}}\paren{\snr}$ is equal to
\begin{eqnarray}
C_{\lambda, {\rm sum}}\paren{\snr} &=& \int_0^\infty \PR{\log\paren{1+I_\lambda\paren{\snr}} > x} \diff x \nonumber \\
&=& \PR{I_\lambda\paren{\snr} > \e{x}-1} \diff x \nonumber \\
&=& \int_0^\infty 1 - \PR{\frac{I_\lambda\paren{\snr} - \ES{I_\lambda\paren{\snr}}}{\sqrt{\V{I_\lambda\paren{\snr}}}} \leq \frac{\e{x} - 1 - \ES{I_\lambda\paren{\snr}}}{\sqrt{\V{I_\lambda\paren{\snr}}}}} \diff x. \hspace{0.4cm} \label{Eqn: Sum Capacity Expression}
\end{eqnarray}
We complete the proof by substituting the bounds in Theorem \ref{Thm: Rates of Convergence} in \eqref{Eqn: Sum Capacity Expression} and using the natural bounds $0$ and $1$ on the probability.     
\end{IEEEproof}

\begin{figure*}[!t]
\begin{minipage}[t]{\textwidth}
\begin{center}
%
%
\begin{psfrags}%
\psfragscanon%
%
\psfrag{s09}[t][t]{\color[rgb]{0,0,0}\setlength{\tabcolsep}{0pt}\scriptsize \begin{tabular}{c}Transmitter Intensity ($\lambda$)\end{tabular}}%
\psfrag{s10}[b][b]{\color[rgb]{0,0,0}\setlength{\tabcolsep}{0pt}\scriptsize \begin{tabular}{c}Sum Capacity and Bounds \\ (Nats per Second per Hertz)\end{tabular}}%
\psfrag{s12}[b][b]{\color[rgb]{0,0,0}\setlength{\tabcolsep}{0pt}\scriptsize \begin{tabular}{c}$G_1(t) = \frac{1}{(1+t)^{4}}$\end{tabular}}%
\psfrag{s14}[][]{\color[rgb]{0,0,0}\setlength{\tabcolsep}{0pt}\scriptsize \begin{tabular}{c} $ $ \end{tabular}}%
\psfrag{s15}[][]{\color[rgb]{0,0,0}\setlength{\tabcolsep}{0pt}\scriptsize \begin{tabular}{c} $ $ \end{tabular}}%
\psfrag{s16}[l][l]{\color[rgb]{0,0,0}\scriptsize Lower Bound}%
\psfrag{s17}[l][l]{\color[rgb]{0,0,0}\scriptsize Simulated Rate}%
\psfrag{s18}[l][l]{\color[rgb]{0,0,0}\scriptsize Upper Bound}%
\psfrag{s19}[l][l]{\color[rgb]{0,0,0}\scriptsize Lower Bound}%
%
\psfrag{x01}[t][t]{\scriptsize 0}%
\psfrag{x02}[t][t]{ \scriptsize 0.1}%
\psfrag{x03}[t][t]{\scriptsize 0.2}%
\psfrag{x04}[t][t]{\scriptsize 0.3}%
\psfrag{x05}[t][t]{\scriptsize 0.4}%
\psfrag{x06}[t][t]{\scriptsize 0.5}%
\psfrag{x07}[t][t]{\scriptsize 0.6}%
\psfrag{x08}[t][t]{\scriptsize 0.7}%
\psfrag{x09}[t][t]{\scriptsize 0.8}%
\psfrag{x10}[t][t]{\scriptsize 0.9}%
\psfrag{x11}[t][t]{\scriptsize 1}%
\psfrag{x12}[t][t]{\scriptsize 1}%
\psfrag{x13}[t][t]{\scriptsize 2}%
\psfrag{x14}[t][t]{\scriptsize 3}%
\psfrag{x15}[t][t]{\scriptsize 4}%
\psfrag{x16}[t][t]{\scriptsize 5}%
\psfrag{x17}[t][t]{\scriptsize 6}%
\psfrag{x18}[t][t]{\scriptsize 7}%
\psfrag{x19}[t][t]{\scriptsize 8}%
\psfrag{x20}[t][t]{\scriptsize 9}%
\psfrag{x21}[t][t]{\scriptsize 10}%
%
\psfrag{v01}[r][r]{\scriptsize 0}%
\psfrag{v02}[r][r]{\scriptsize 0.1}%
\psfrag{v03}[r][r]{\scriptsize 0.2}%
\psfrag{v04}[r][r]{\scriptsize 0.3}%
\psfrag{v05}[r][r]{\scriptsize 0.4}%
\psfrag{v06}[r][r]{\scriptsize 0.5}%
\psfrag{v07}[r][r]{\scriptsize 0.6}%
\psfrag{v08}[r][r]{\scriptsize 0.7}%
\psfrag{v09}[r][r]{\scriptsize 0.8}%
\psfrag{v10}[r][r]{\scriptsize 0.9}%
\psfrag{v11}[r][r]{\scriptsize 1}%
\psfrag{v12}[r][r]{\scriptsize 0}%
\psfrag{v13}[r][r]{\scriptsize 0.5}%
\psfrag{v14}[r][r]{\scriptsize 1}%
\psfrag{v15}[r][r]{\scriptsize 1.5}%
\psfrag{v16}[r][r]{\scriptsize 2}%
\psfrag{v17}[r][r]{\scriptsize 2.5}%
\psfrag{v18}[r][r]{\scriptsize 3}%
%
\includegraphics[scale=0.55]{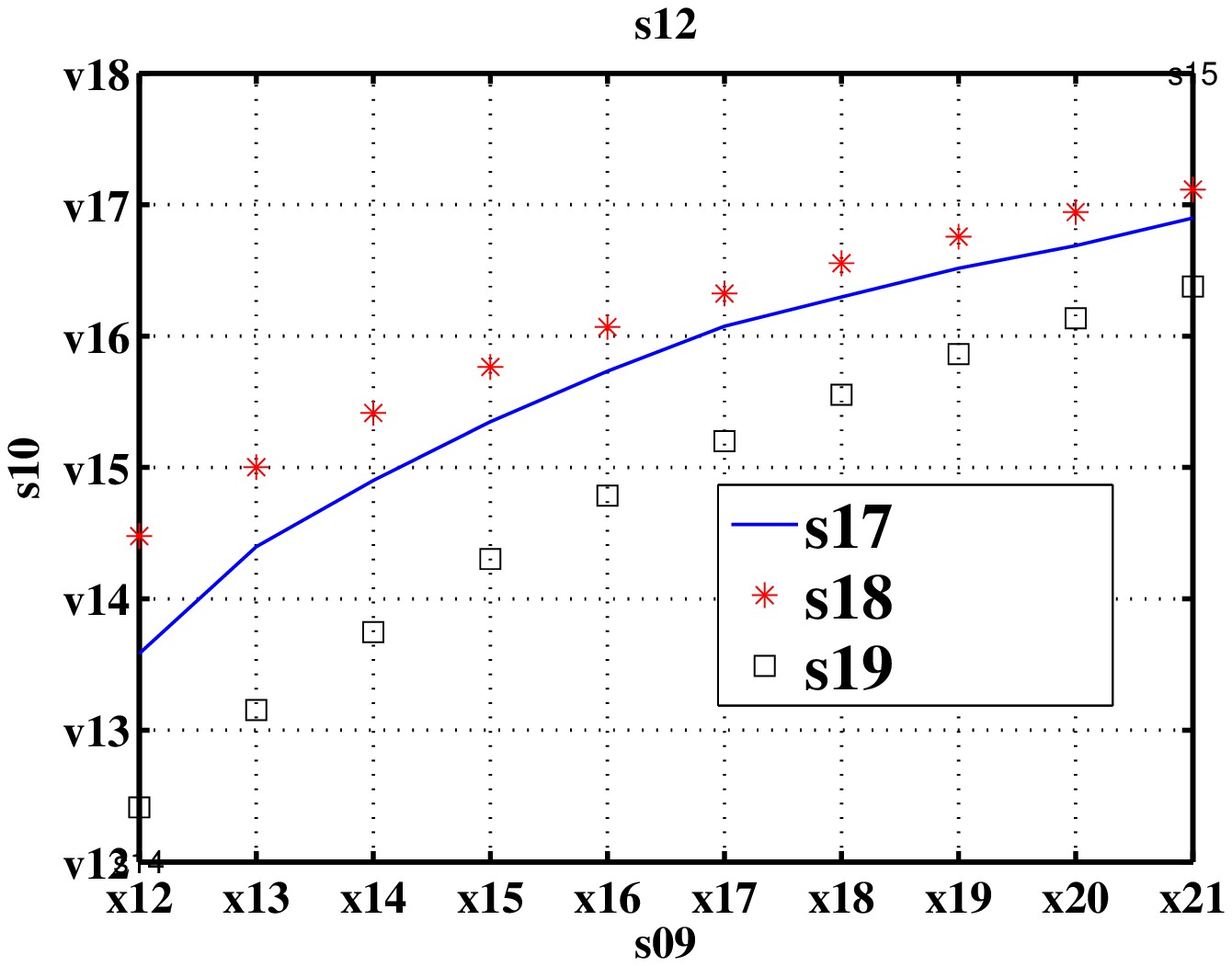}%
\end{psfrags}%
%

%
\hspace{\fill}
%
%
%
\begin{psfrags}%
\psfragscanon%
%
\psfrag{s09}[t][t]{\color[rgb]{0,0,0}\setlength{\tabcolsep}{0pt}\scriptsize \begin{tabular}{c}Transmitter Intensity ($\lambda$)\end{tabular}}%
\psfrag{s10}[b][b]{\color[rgb]{0,0,0}\setlength{\tabcolsep}{0pt}\scriptsize \begin{tabular}{c}Sum Capacity and Bounds \\ (Nats per Second per Hertz)\end{tabular}}%
\psfrag{s12}[b][b]{\color[rgb]{0,0,0}\setlength{\tabcolsep}{0pt}\scriptsize \begin{tabular}{c}$G_1(t) = \frac{1}{(1+t)^{4}}$\end{tabular}}%
\psfrag{s14}[][]{\color[rgb]{0,0,0}\setlength{\tabcolsep}{0pt}\scriptsize \begin{tabular}{c} $ $ \end{tabular}}%
\psfrag{s15}[][]{\color[rgb]{0,0,0}\setlength{\tabcolsep}{0pt}\scriptsize \begin{tabular}{c} $ $ \end{tabular}}%
\psfrag{s16}[l][l]{\color[rgb]{0,0,0}\scriptsize Lower Bound}%
\psfrag{s17}[l][l]{\color[rgb]{0,0,0}\scriptsize Simulated Rate}%
\psfrag{s18}[l][l]{\color[rgb]{0,0,0}\scriptsize Upper Bound}%
\psfrag{s19}[l][l]{\color[rgb]{0,0,0}\scriptsize Lower Bound}%
%
\psfrag{x01}[t][t]{\scriptsize 0}%
\psfrag{x02}[t][t]{\scriptsize 0.1}%
\psfrag{x03}[t][t]{\scriptsize 0.2}%
\psfrag{x04}[t][t]{\scriptsize 0.3}%
\psfrag{x05}[t][t]{\scriptsize 0.4}%
\psfrag{x06}[t][t]{\scriptsize 0.5}%
\psfrag{x07}[t][t]{\scriptsize 0.6}%
\psfrag{x08}[t][t]{\scriptsize 0.7}%
\psfrag{x09}[t][t]{\scriptsize 0.8}%
\psfrag{x10}[t][t]{\scriptsize 0.9}%
\psfrag{x11}[t][t]{\scriptsize 1}%
\psfrag{x12}[t][t]{\scriptsize 100}%
\psfrag{x13}[t][t]{\scriptsize 20}%
\psfrag{x14}[t][t]{\scriptsize 30}%
\psfrag{x15}[t][t]{\scriptsize 40}%
\psfrag{x16}[t][t]{\scriptsize 50}%
\psfrag{x17}[t][t]{\scriptsize 60}%
\psfrag{x18}[t][t]{\scriptsize 70}%
\psfrag{x19}[t][t]{\scriptsize 80}%
\psfrag{x20}[t][t]{\scriptsize 90}%
\psfrag{x21}[t][t]{\scriptsize 100}%
%
\psfrag{v01}[r][r]{\scriptsize 0}%
\psfrag{v02}[r][r]{\scriptsize 0.1}%
\psfrag{v03}[r][r]{\scriptsize 0.2}%
\psfrag{v04}[r][r]{\scriptsize 0.3}%
\psfrag{v05}[r][r]{\scriptsize 0.4}%
\psfrag{v06}[r][r]{\scriptsize 0.5}%
\psfrag{v07}[r][r]{\scriptsize 0.6}%
\psfrag{v08}[r][r]{\scriptsize 0.7}%
\psfrag{v09}[r][r]{\scriptsize 0.8}%
\psfrag{v10}[r][r]{\scriptsize 0.9}%
\psfrag{v11}[r][r]{\scriptsize 1}%
\psfrag{v12}[r][r]{\scriptsize 2}%
\psfrag{v13}[r][r]{\scriptsize 2.5}%
\psfrag{v14}[r][r]{\scriptsize 3}%
\psfrag{v15}[r][r]{\scriptsize 3.5}%
\psfrag{v16}[r][r]{\scriptsize 4}%
\psfrag{v17}[r][r]{\scriptsize 4.5}%
\psfrag{v18}[r][r]{\scriptsize 5}%
%
\includegraphics[scale=0.55]{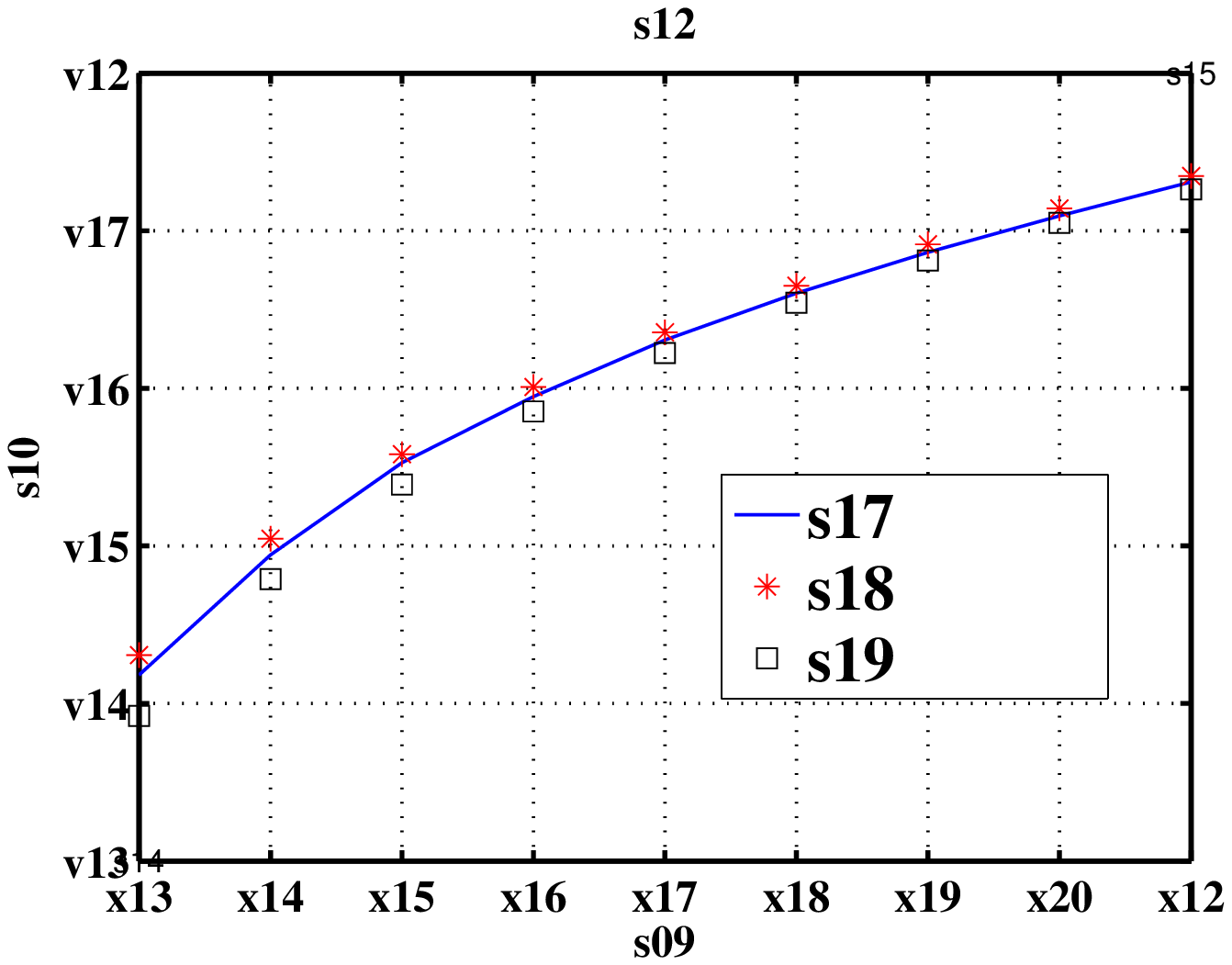}%
\end{psfrags}%
%

\end{center}
\end{minipage} 
\\ 
\begin{minipage}[t]{\textwidth}
\vspace{-2mm} 
\begin{center}
%
%
\begin{psfrags}%
\psfragscanon%
%
\psfrag{s09}[t][t]{\color[rgb]{0,0,0}\setlength{\tabcolsep}{0pt}\scriptsize \begin{tabular}{c}Transmitter Intensity ($\lambda$)\end{tabular}}%
\psfrag{s10}[b][b]{\color[rgb]{0,0,0}\setlength{\tabcolsep}{0pt}\scriptsize \begin{tabular}{c}Sum Capacity and Bounds \\ (Nats per Second per Hertz)\end{tabular}}%
\psfrag{s12}[b][b]{\color[rgb]{0,0,0}\setlength{\tabcolsep}{0pt}\scriptsize \begin{tabular}{c}$G_2(t) = \frac{1}{1+t^{4}}$\end{tabular}}%
\psfrag{s14}[][]{\color[rgb]{0,0,0}\setlength{\tabcolsep}{0pt}\scriptsize \begin{tabular}{c} $ $ \end{tabular}}%
\psfrag{s15}[][]{\color[rgb]{0,0,0}\setlength{\tabcolsep}{0pt}\scriptsize \begin{tabular}{c} $ $ \end{tabular}}%
\psfrag{s16}[l][l]{\color[rgb]{0,0,0}\scriptsize Lower Bound}%
\psfrag{s17}[l][l]{\color[rgb]{0,0,0}\scriptsize Simulated Rate}%
\psfrag{s18}[l][l]{\color[rgb]{0,0,0}\scriptsize Upper Bound}%
\psfrag{s19}[l][l]{\color[rgb]{0,0,0}\scriptsize Lower Bound}%
%
\psfrag{x01}[t][t]{\scriptsize 0}%
\psfrag{x02}[t][t]{\scriptsize 0.1}%
\psfrag{x03}[t][t]{\scriptsize 0.2}%
\psfrag{x04}[t][t]{\scriptsize 0.3}%
\psfrag{x05}[t][t]{\scriptsize 0.4}%
\psfrag{x06}[t][t]{\scriptsize 0.5}%
\psfrag{x07}[t][t]{\scriptsize 0.6}%
\psfrag{x08}[t][t]{\scriptsize 0.7}%
\psfrag{x09}[t][t]{\scriptsize 0.8}%
\psfrag{x10}[t][t]{\scriptsize 0.9}%
\psfrag{x11}[t][t]{\scriptsize 1}%
\psfrag{x12}[t][t]{\scriptsize 1}%
\psfrag{x13}[t][t]{\scriptsize 2}%
\psfrag{x14}[t][t]{\scriptsize 3}%
\psfrag{x15}[t][t]{\scriptsize 4}%
\psfrag{x16}[t][t]{\scriptsize 5}%
\psfrag{x17}[t][t]{\scriptsize 6}%
\psfrag{x18}[t][t]{\scriptsize 7}%
\psfrag{x19}[t][t]{\scriptsize 8}%
\psfrag{x20}[t][t]{\scriptsize 9}%
\psfrag{x21}[t][t]{\scriptsize 10}%
%
\psfrag{v01}[r][r]{\scriptsize 0}%
\psfrag{v02}[r][r]{\scriptsize 0.1}%
\psfrag{v03}[r][r]{\scriptsize 0.2}%
\psfrag{v04}[r][r]{\scriptsize 0.3}%
\psfrag{v05}[r][r]{\scriptsize 0.4}%
\psfrag{v06}[r][r]{\scriptsize 0.5}%
\psfrag{v07}[r][r]{\scriptsize 0.6}%
\psfrag{v08}[r][r]{\scriptsize 0.7}%
\psfrag{v09}[r][r]{\scriptsize 0.8}%
\psfrag{v10}[r][r]{\scriptsize 0.9}%
\psfrag{v11}[r][r]{\scriptsize 1}%
\psfrag{v12}[r][r]{\scriptsize 1}%
\psfrag{v13}[r][r]{\scriptsize 1.5}%
\psfrag{v14}[r][r]{\scriptsize 2}%
\psfrag{v15}[r][r]{\scriptsize 2.5}%
\psfrag{v16}[r][r]{\scriptsize 3}%
\psfrag{v17}[r][r]{\scriptsize 3.5}%
\psfrag{v18}[r][r]{\scriptsize 4}%
%
\includegraphics[scale=0.55]{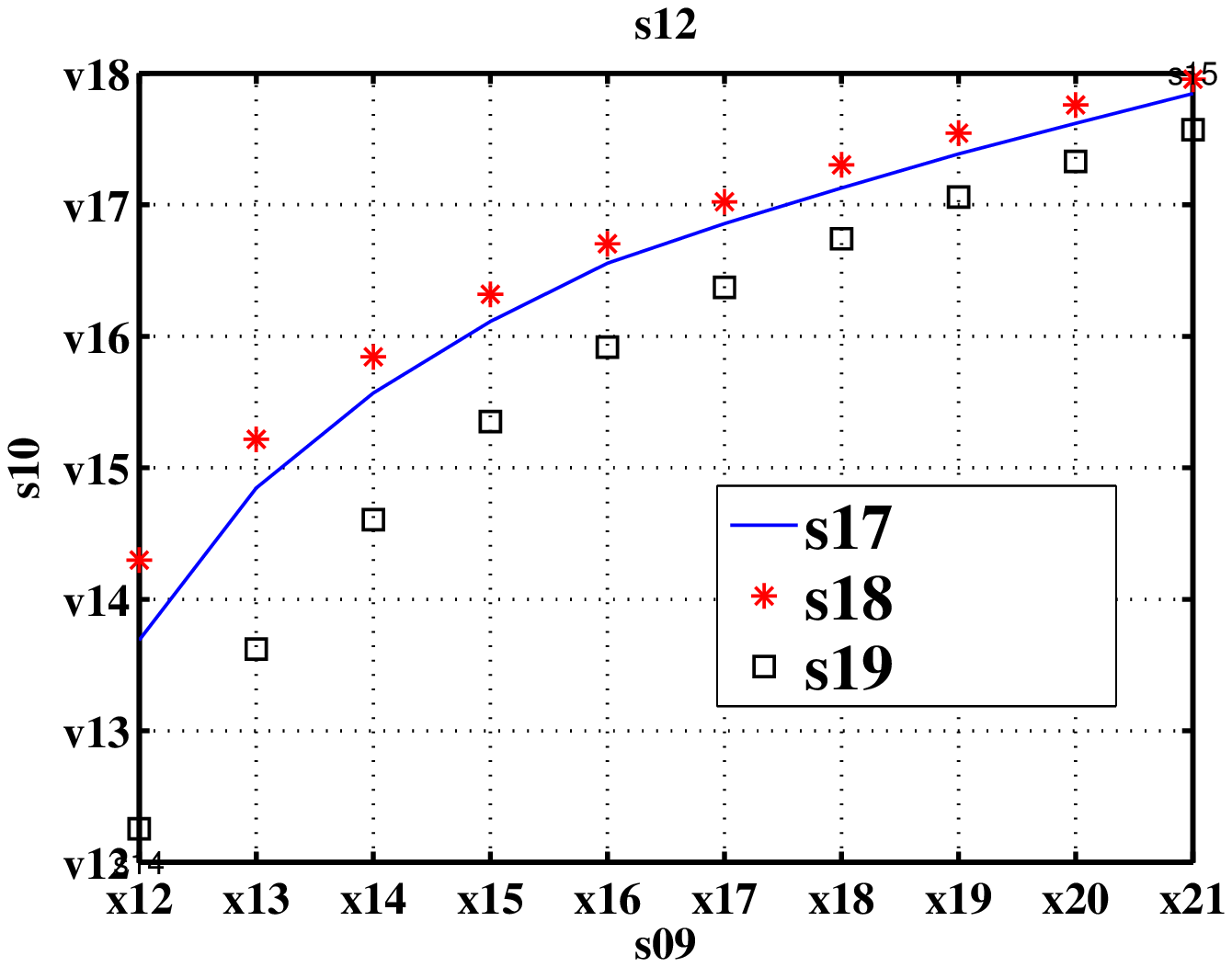}%
\end{psfrags}%
%

%
\hfill
%
%
%
\begin{psfrags}%
\psfragscanon%
%
\psfrag{s09}[t][t]{\color[rgb]{0,0,0}\setlength{\tabcolsep}{0pt}\scriptsize \begin{tabular}{c}Transmitter Intensity ($\lambda$)\end{tabular}}%
\psfrag{s10}[b][b]{\color[rgb]{0,0,0}\setlength{\tabcolsep}{0pt}\scriptsize \begin{tabular}{c}Sum Capacity and Bounds \\ (Nats per Second per Hertz)\end{tabular}}%
\psfrag{s12}[b][b]{\color[rgb]{0,0,0}\setlength{\tabcolsep}{0pt}\scriptsize \begin{tabular}{c}$G_2(t) = \frac{1}{1+t^{4}}$\end{tabular}}%
\psfrag{s14}[][]{\color[rgb]{0,0,0}\setlength{\tabcolsep}{0pt}\scriptsize \begin{tabular}{c} $ $ \end{tabular}}%
\psfrag{s15}[][]{\color[rgb]{0,0,0}\setlength{\tabcolsep}{0pt}\scriptsize \begin{tabular}{c} $ $ \end{tabular}}%
\psfrag{s16}[l][l]{\color[rgb]{0,0,0}\scriptsize Lower Bound}%
\psfrag{s17}[l][l]{\color[rgb]{0,0,0}\scriptsize Simulated Rate}%
\psfrag{s18}[l][l]{\color[rgb]{0,0,0}\scriptsize Upper Bound}%
\psfrag{s19}[l][l]{\color[rgb]{0,0,0}\scriptsize Lower Bound}%
%
\psfrag{x01}[t][t]{\scriptsize 0}%
\psfrag{x02}[t][t]{\scriptsize 0.1}%
\psfrag{x03}[t][t]{\scriptsize 0.2}%
\psfrag{x04}[t][t]{\scriptsize 0.3}%
\psfrag{x05}[t][t]{\scriptsize 0.4}%
\psfrag{x06}[t][t]{\scriptsize 0.5}%
\psfrag{x07}[t][t]{\scriptsize 0.6}%
\psfrag{x08}[t][t]{\scriptsize 0.7}%
\psfrag{x09}[t][t]{\scriptsize 0.8}%
\psfrag{x10}[t][t]{\scriptsize 0.9}%
\psfrag{x11}[t][t]{\scriptsize 1}%
\psfrag{x12}[t][t]{\scriptsize 100}%
\psfrag{x13}[t][t]{\scriptsize 20}%
\psfrag{x14}[t][t]{\scriptsize 30}%
\psfrag{x15}[t][t]{\scriptsize 40}%
\psfrag{x16}[t][t]{\scriptsize 50}%
\psfrag{x17}[t][t]{\scriptsize 60}%
\psfrag{x18}[t][t]{\scriptsize 70}%
\psfrag{x19}[t][t]{\scriptsize 80}%
\psfrag{x20}[t][t]{\scriptsize 90}%
\psfrag{x21}[t][t]{\scriptsize 100}%
%
\psfrag{v01}[r][r]{\scriptsize 0}%
\psfrag{v02}[r][r]{\scriptsize 0.1}%
\psfrag{v03}[r][r]{\scriptsize 0.2}%
\psfrag{v04}[r][r]{\scriptsize 0.3}%
\psfrag{v05}[r][r]{\scriptsize 0.4}%
\psfrag{v06}[r][r]{\scriptsize 0.5}%
\psfrag{v07}[r][r]{\scriptsize 0.6}%
\psfrag{v08}[r][r]{\scriptsize 0.7}%
\psfrag{v09}[r][r]{\scriptsize 0.8}%
\psfrag{v10}[r][r]{\scriptsize 0.9}%
\psfrag{v11}[r][r]{\scriptsize 1}%
\psfrag{v12}[r][r]{\scriptsize 3.5}%
\psfrag{v13}[r][r]{\scriptsize 4}%
\psfrag{v14}[r][r]{\scriptsize 4.5}%
\psfrag{v15}[r][r]{\scriptsize 5}%
\psfrag{v16}[r][r]{\scriptsize 5.5}%
\psfrag{v17}[r][r]{\scriptsize 6}%
\psfrag{v18}[r][r]{\scriptsize 6.5}%
%
\includegraphics[scale=0.55]{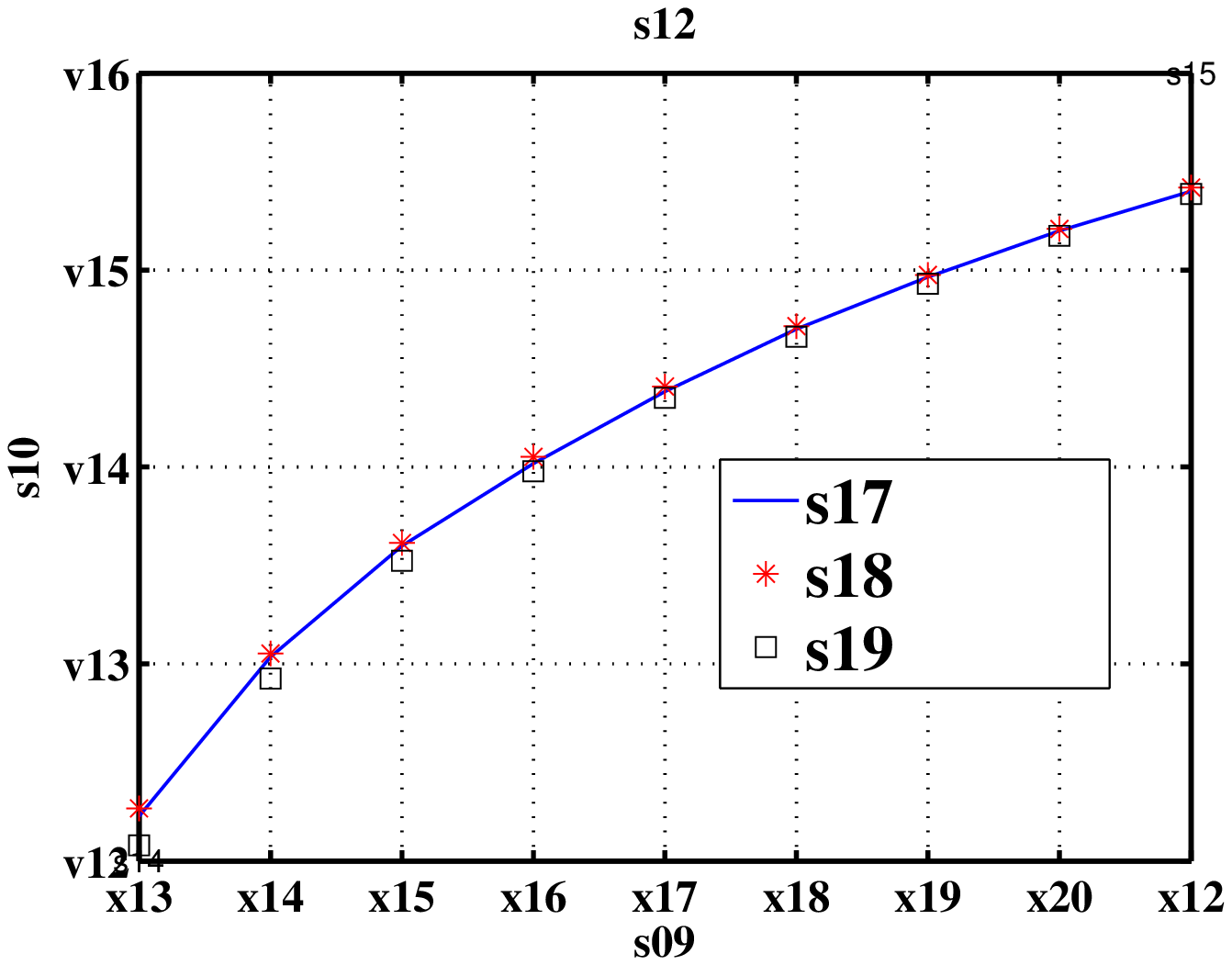}%
\end{psfrags}%
%

\end{center}
\end{minipage}
\caption{Upper and lower bounds on $C_{\lambda, {\rm sum}}\paren{\snr}$ for the path-loss functions $G_1(t) = \frac{1}{\paren{1+t}^\alpha}$ (top figures) and $G_2(t) = \frac{1}{1+t^\alpha}$ (bottom figures).  Transmitters are distributed over $\R^2$ according to a stationary PPP of intensity $\lambda$.  Fading effects are modeled by using the Nakagami-$m$ fading model with $m$ parameter set to $5$. For this model of spatial distribution of interfering transmitters, $c(x)$ is given as $c(x) = \frac{1}{\sqrt{2 \pi}}\frac{m_{H^3}}{\paren{m_{H^2}}^\frac32}\frac{\int_{0}^\infty G^3(t) t \diff t}{\paren{\int_{0}^\infty G^2(t) t \diff t}^\frac32} \min\paren{0.4785, \frac{31.935}{1+|x|^3}}$. ($\snr = 0$ [dB] and $\alpha = 4$)} \label{Fig: Sum Capacity Bounds}
\end{figure*}

We plot $C_{\lambda, {\rm sum}}\paren{\snr}$ and the corresponding bounds in Theorem \ref{Thm: Sum Capacity Bounds} as a function of $\lambda$ in Fig. \ref{Fig: Sum Capacity Bounds} for both path-loss functions.  For this simulation study, we assume that transmitters are {\em uniformly} distributed over $\R^2$ according to a stationary PPP of intensity $\lambda$ [nodes per unit area].  We set $\snr$ to $0$ [dB], and the path-loss exponent $\alpha$ to $4$.  Fading effects are modeled by means of the Nakagami-$m$ fading with unit mean power gain and $m=5$.  Similar observations continue to hold for other parameter selections.  Again, we observe that our bounds closely approximate $C_{\lambda, {\rm sum}}\paren{\snr}$ within {\em one} nats per second per hertz for moderate to high values of $\lambda$. Especially, for high values of $\lambda$, our bounds become very tight, and they {\em almost} coincide with $C_{\lambda, {\rm sum}}\paren{\snr}$.  
With increasing values of $\lambda$, it is observed that $C_{\lambda, {\rm sum}}\paren{\snr}$ grows logarithmically\footnote{Expected total received power at the BS grows linearly with $\lambda$, {\em i.e.,} see \eqref{Eqn: Interference Mean} in Appendix \ref{Appendix: Mean Variance Convergence}.}, which is also a property demonstrated by the derived bounds.                    

\section{Conclusions} \label{Section: Conclusions}
In this paper, we have focused on the statistical characterization of wireless multi-access interference distributions when transmitters are spatially distributed over the plane according to a Poisson point process, which is not necessarily stationary.  The signal propagation at the physical layer is modeled by means of a general bounded and power-law decaying path-loss function.  Other random wireless channel dynamics such as fading and shadowing are also incorporated in the employed signal propagation model. 

We have shown that the wireless multi-access interference distribution converges to the Gaussian distribution with the same mean and variance at a rate $\frac{c(x)}{\sqrt{\lambda}}$, where $\lambda$ is a modeling parameter controlling the ``intensity" of the planar Poisson point process generating transmitter locations, and $c(x)>0$ is a function which depends on the shape of the path-loss function and on the point $x \in \Re$ at which we want to estimate the interference distribution.  $c(x)$ decays to zero at a rate $|x|^{-3}$ as the absolute value of $x$ increases.  This behavior makes the derived bounds tight for any finite value of $\lambda$.  An explicit expression for $c(x)$ appearing in our approximation results has also been provided, {\em i.e.,} see Theorem \ref{Thm: Rates of Convergence}.  

We have performed an extensive numerical and simulation study to verify the derived theoretical bounds.  A very good statistical match between the simulated (centered and normalized) wireless multi-access interference distributions and the normal distribution with zero mean and variance one has been observed even for moderately small values of $\lambda$.  Since there are no closed form expressions available for the wireless multi-access interference distributions under general bounded path-loss models at the present, these results are expected to help researchers in the field by simplifying the derivation of closed form expressions for various performance bounds and metrics for wireless networks.  In particular, applications of our Gaussian approximation results have been illustrated to derive tight upper and lower bounds on the outage capacity for a given victim link in a spatial wireless network and those on the sum capacity for spatial wireless multiple-access networks.  It has been observed that the derived performance bounds can approximate these capacity metrics up to one nats per second per hertz for moderate to high values of $\lambda$.  With increasing values of $\lambda$, the approximation precision improves significantly, and the derived performance bounds almost overlap with the simulated outage capacity and sum capacity curves.   

\vspace{-0.1cm}
\section{Acknowledgements}
The authors thanks Stephen V. Hanly for his valuable comments on the earlier versions of this work. 

\vspace{-0.3cm}
\appendices
\section{Proof of Lemma \ref{Lemma: Non-degenerate Distribution}} \label{Appendix: Non-degenerate Distribution} 
We can find constant $B_1 > 0$ and $K>0$ such that $G(t) \leq K t^{-\alpha}$ for all $t \geq B_1$.  Then, 
\begin{eqnarray}
\lefteqn{\int_0^\infty \int_0^\infty \paren{1-\e{-s P h G(t)}} q(h) p(t) \diff t \diff h} \hspace{16cm} \nonumber \\ 
\lefteqn{\leq \int_0^\infty \int_0^{B_1} \paren{1-\e{-s P h G(t)}} q(h) p(t) \diff t \diff h + \int_0^\infty \int_{B_1}^\infty \paren{1-\e{-s P K h t^{-\alpha}}} q(h) p(t) \diff t \diff h} \hspace{14.5cm} \nonumber \\
\lefteqn{\leq \int_0^{B_1} p(t) \diff t + \int_0^\infty \int_{B_1}^\infty \paren{1-\e{-s P K h t^{-\alpha}}} q(h) p(t) \diff t \diff h.} \hspace{14.5cm} \label{Eqn: Finiteness of Integrals 1} 
\end{eqnarray}   

The first integral in \eqref{Eqn: Finiteness of Integrals 1} is finite since $\Lambda$ is locally finite.  To show the finiteness of the second integral, we divide it into two parts as follows.   
\begin{eqnarray}
\lefteqn{\int_0^\infty \int_{B_1}^\infty \paren{1-\e{-s P K h t^{-\alpha}}} q(h) p(t) \diff t \diff h} \hspace{16cm} \nonumber \\ 
\lefteqn{= \int_0^\frac{1}{sPK} \int_{B_1}^\infty \paren{1-\e{-s P K h t^{-\alpha}}} q(h) p(t) \diff t \diff h +  \int_\frac{1}{sPK}^\infty \int_{B_1}^\infty \paren{1-\e{-s P K h t^{-\alpha}}} q(h) p(t) \diff t \diff h.} \hspace{15.5cm} \label{Eqn: Finiteness of Integrals 2}
\end{eqnarray}

The first integral in \eqref{Eqn: Finiteness of Integrals 2} can be bounded as
\begin{eqnarray}
 \int_0^\frac{1}{sPK} \int_{B_1}^\infty \paren{1-\e{-s P K h t^{-\alpha}}} q(h) p(t) \diff t \diff h \leq \int_{B_1}^\infty \paren{1 - \e{-t^{-\alpha}}} p(t) \diff t, \nonumber
\end{eqnarray}
which is finite since $1 - \e{-t^{-\alpha}} = \BO{t^{-\alpha}}$ and $p(t) = \BO{t^{\alpha - 1 -\epsilon}}$ as $t \ra \infty$.  Hence, proving the finiteness of $\int_\frac{1}{sPK}^\infty \int_{B_1}^\infty \paren{1-\e{-s P K h t^{-\alpha}}} q(h) p(t) \diff t \diff h$ will complete the proof.  To this end, we need the following lemma.
\begin{lemma} \label{Lemma: An Exp Bound}
$1-\e{-a t^{-\alpha}} \leq 2 a \paren{1 - \e{-a}} t^{-\alpha}$ for all $a \geq 1$ and $t$ large enough. 
\end{lemma}
\begin{IEEEproof}
We let $f_t(a) = 1-\e{-a t^{-\alpha}}$ and $g_t(a) = 2 a \paren{1 - \e{-a}} t^{-\alpha}$.  For $a = 1$, we have $\lim_{t \ra \infty} \frac{f_t(1)}{t^{-\alpha}} = 1$ and $\lim_{t \ra \infty} \frac{g_t(1)}{t^{-\alpha}} = 2 \paren{1 - \e{-1}} > 1$.  Hence, there exists a constant $B_2>0$ such that $g_t(1) \geq f_t(1)$ for all $t \geq B_2$.  Fix an arbitrary $t$ greater than $B_2$.  Then, 
\begin{eqnarray}
\frac{\diff f_t(a)}{\diff a} = t^{-\alpha} \e{-a t^{-\alpha}} \mbox{ and } \frac{\diff g_t(a)}{\diff a} = 2t^{-\alpha}\paren{1+a\e{-a} - \e{-a}}. \nonumber
\end{eqnarray}
Thus, $g_t(a)$ grows faster than $f_t(a)$, implying that $g_t(a) \geq f_t(a)$ for all $a \geq 1$ and $t \geq B_2$.    
\end{IEEEproof}

By using Lemma \ref{Lemma: An Exp Bound}, we can upper bound the second integral in \eqref{Eqn: Finiteness of Integrals 2} as
\begin{eqnarray}
\lefteqn{\int_\frac{1}{sPK}^\infty \int_{B_1}^\infty \paren{1-\e{-s P K h t^{-\alpha}}} q(h) p(t) \diff t \diff h} \hspace{16cm} \nonumber \\ 
\lefteqn{\leq \int_{B_1}^{B_3} p(t) \diff t + \int_{B_3}^\infty  \int_{\frac{1}{sPK}}^\infty 2 s P K h \paren{1 - \e{-sPKh}} q(h) t^{-\alpha} p(t) \diff h \diff t} \hspace{14cm} \label{Eqn: Finiteness of Integrals 3}
\end{eqnarray}
for some positive constant $B_3$ large enough.  The first integral in \eqref{Eqn: Finiteness of Integrals 3} is finite due to local finiteness of $\Lambda$.  The second integral in \eqref{Eqn: Finiteness of Integrals 3} can be upper bounded by $2 s P K m_H \int_{B_3}^\infty t^{-\alpha} p(t) \diff t$, which is finite since $m_H < \infty$ and $p(t) = \BO{t^{\alpha - 1 - \epsilon}}$ as $t \ra \infty$. 

\section{Proof of Lemma \ref{Lemma: Distribution Convergence}} \label{Appendix: Distribution Convergence}
 We will show that $\Lap{I_n}(s)$ converges to $\Lap{I_\lambda}(s)$ pointwise as $n$ tends to infinity.  Since $U_{1, n}, \ldots, U_{\Ceil{\Lambda_n}, n}$ are independent, we have
\begin{eqnarray}
\Lap{I_n}(s) = \paren{1 - \frac{\lambda}{\Lambda_n} \int_0^\infty \int_0^n \paren{1 - \e{-s P h G(t)}} p(t) q(h) \diff t \diff h}^{\Ceil{\Lambda_n}}. \nonumber
\end{eqnarray}
We have $\int_0^\infty \int_0^n \paren{1 - \e{-s P h G(t)}} p(t) q(h) \diff t \diff h$ converging to $\int_0^\infty \int_0^\infty \paren{1 - \e{-s P h G(t)}} p(t) q(h) \diff t \diff h$ as $n$ tends to infinity, and $\int_0^\infty \int_0^\infty \paren{1 - \e{-s P h G(t)}} p(t) q(h) \diff t \diff h < \infty$ by Lemma \ref{Lemma: Non-degenerate Distribution}.  Hence, by observing that $\Lambda\paren{\R^2} = \infty$, we have    
\begin{eqnarray}
\lim_{n \ra \infty} \Lap{I_n}(s) &=& \exp\paren{-\lambda \int_0^\infty \int_0^\infty \paren{1 - \e{-s P h G(t)}} q(h) p(t) \diff t \diff h} \nonumber \\
&=& \Lap{I_\lambda}(s). \nonumber
\end{eqnarray}  

\section{Proof of Lemma \ref{Lemma: Mean Variance Convergence}} \label{Appendix: Mean Variance Convergence}
By using Campbell's Theorem \cite{Kingman93}, we have 
\begin{eqnarray}
\ES{I_\lambda} = \lambda P m_H \int_0^\infty G(t)p(t) \diff t \label{Eqn: Interference Mean}
\end{eqnarray} 
and 
\begin{eqnarray}
\V{I_\lambda} = \lambda P^2 m_{H^2} \int_0^\infty G^2(t)p(t) \diff t. \label{Eqn: Interference Variance}
\end{eqnarray}
Note that our assumptions on $G(t)$ and $p(t)$ ensure that $\ES{I_\lambda}$ and $\V{I_\lambda}$ are both finite.  Let $U_{k, n}$ be defined as in the proof of Lemma \ref{Lemma: Distribution Convergence}.  Let also $m_{k, n} = \ES{P H_k G\paren{U_{k, n}}}$ and $\sigma^2_{k, n} = \V{P H_k G\paren{U_{k, n}}}$.   Then, $\ES{I_n} = \Ceil{\Lambda_n} m_{1, n}$ and $\V{I_n} = \Ceil{\Lambda_n}\sigma^2_{1, n}$.  We can explicitly write $m_{1, n}$ as 
\begin{eqnarray}
m_{1, n} = \frac{\lambda P m_H}{\Lambda_n} \int_0^n G(t)p(t)\diff t, \nonumber
\end{eqnarray}
which implies that $\lim_{n \ra \infty} \ES{I_n} = \ES{I_\lambda}$ since $\Lambda\paren{\R^2} = \infty$.  Similarly, we have 
\begin{eqnarray}
\sigma^2_{1, n} = \frac{\lambda P^2 m_{H^2}}{\Lambda_n}\int_0^n G^2(t)p(t)\diff t - \frac{\lambda^2P^2m_H^2}{\Lambda_n^2}\paren{\int_0^nG(t)p(t)\diff t}^2, \nonumber
\end{eqnarray}
which implies that $\lim_{n \ra \infty}\V{I_n} = \V{I_\lambda}$.        

\section{Normal Approximation Results for a Non-stationary PPP} \label{Appendix: Normal Approximation Non-stationary PPP}
Our normal approximation bounds given in Theorem \ref{Thm: Rates of Convergence} are valid for both  stationary and non-stationary PPPs.  In this appendix, we will illustrate the validity and utility of these bounds for a non-stationary PPP.  To this end, we place the test receiver node to the origin, {\em i.e.,} $\vec{X}_o = \vec{0}$, and consider a non-stationary PPP, still denoted by $\Phi_\Lambda$, with mean measure {\em density} $f$ given as 
\begin{eqnarray}
f\paren{\vec{x}} = \left\{ \begin{array}{c c} \frac{\lambda}{\|\vec{x}\|_2^2}& \mbox{ if } \|\vec{x}\|_2 \geq r \\ 0& \mbox{ if } \|\vec{x}\|_2 < r \end{array} \right., \label{Eqn: Density 1 - Appendix}
\end{eqnarray}
where $r>0$ is a given positive constant.  As in Section \ref{Section: Applications}, $r$ can be interpreted as the radius of an {\em exclusion zone} in which no interfering transmitters are allowed, perhaps due to physical limitations, or the operation of the MAC layer ({\em i.e.,} CSMA-CA protocol), or the interference cancellation property \cite{WAYV07}.  Since the test receiver node is located at the origin, we have $T\paren{\vec{x}} = \|\vec{x}\|_2$, and 
\begin{eqnarray}
\Lambda \circ T^{-1}\paren{[0, t]} &=& \Lambda\brparen{\vec{x} \in \R^2: \|\vec{x}\|_2 \leq t} \nonumber \\
&=& \int_{\R^2} \frac{\lambda}{\|\vec{x}\|_2^2} \I{r \leq \|\vec{x}\|_2 \leq t} \diff \vec{x} \nonumber \\
&=& 2 \pi \lambda \paren{\log(t) - \log(r)} \nonumber
\end{eqnarray}
for $t \geq r$.  Hence, the transformed PPP $\sum_{k \geq 1} \delta_{T\paren{\vec{X}_k}}$ has the mean measure density $p_\lambda(t) = \lambda \frac{2 \pi}{t} \I{t \geq r}$.  Finally, we have the following theorem approximating the WMAI distribution in this case. 

\begin{figure*}[!t]
\begin{minipage}[t]{\textwidth}
\begin{center}
%
%
\begin{psfrags}%
\psfragscanon%
%
\psfrag{s09}[t][t]{\color[rgb]{0,0,0}\setlength{\tabcolsep}{0pt}\scriptsize \begin{tabular}{c}Centered and Normalized Interference Power\end{tabular}}%
\psfrag{s10}[b][b]{\color[rgb]{0,0,0}\setlength{\tabcolsep}{0pt}\scriptsize \begin{tabular}{c}Normal CDF and Bounds\end{tabular}}%
\psfrag{s12}[b][b]{\color[rgb]{0,0,0}\setlength{\tabcolsep}{0pt}\scriptsize \begin{tabular}{c}$G_1(t) = \frac{1}{(1+t)^{4}}$ (without fading) \end{tabular}}%
\psfrag{s14}[][]{\color[rgb]{0,0,0}\setlength{\tabcolsep}{0pt}\begin{tabular}{c} $ $ \end{tabular}}%
\psfrag{s15}[][]{\color[rgb]{0,0,0}\setlength{\tabcolsep}{0pt}\begin{tabular}{c} $ $\end{tabular}}%
\psfrag{s17}[l][l]{\color[rgb]{0,0,0}\scriptsize Normal}%
\psfrag{s18}[l][l]{\color[rgb]{0,0,0}\scriptsize $\lambda = 5$ (LB)}%
\psfrag{s19}[l][l]{\color[rgb]{0,0,0}\scriptsize $\lambda = 5$ (UB)}%
\psfrag{s20}[l][l]{\color[rgb]{0,0,0}\scriptsize $\lambda = 25$ (LB)}%
\psfrag{s21}[l][l]{\color[rgb]{0,0,0}\scriptsize $\lambda = 25$ (UB)}%
\psfrag{s22}[l][l]{\color[rgb]{0,0,0}\scriptsize $\lambda = 100$ (LB)}%
\psfrag{s23}[l][l]{\color[rgb]{0,0,0}\scriptsize $\lambda = 100$ (UB)}%
%
\psfrag{x01}[t][t]{\scriptsize 0}%
\psfrag{x02}[t][t]{\scriptsize 0.1}%
\psfrag{x03}[t][t]{\scriptsize 0.2}%
\psfrag{x04}[t][t]{\scriptsize 0.3}%
\psfrag{x05}[t][t]{\scriptsize 0.4}%
\psfrag{x06}[t][t]{\scriptsize 0.5}%
\psfrag{x07}[t][t]{\scriptsize 0.6}%
\psfrag{x08}[t][t]{\scriptsize 0.7}%
\psfrag{x09}[t][t]{\scriptsize 0.8}%
\psfrag{x10}[t][t]{\scriptsize 0.9}%
\psfrag{x11}[t][t]{\scriptsize 1}%
\psfrag{x12}[t][t]{\scriptsize -6}%
\psfrag{x13}[t][t]{\scriptsize -5}%
\psfrag{x14}[t][t]{\scriptsize -4}%
\psfrag{x15}[t][t]{\scriptsize -3}%
\psfrag{x16}[t][t]{\scriptsize -2}%
\psfrag{x17}[t][t]{\scriptsize -1}%
\psfrag{x18}[t][t]{\scriptsize 0}%
\psfrag{x19}[t][t]{\scriptsize 1}%
\psfrag{x20}[t][t]{\scriptsize 2}%
\psfrag{x21}[t][t]{\scriptsize 3}%
\psfrag{x22}[t][t]{\scriptsize 4}%
\psfrag{x23}[t][t]{\scriptsize 5}%
\psfrag{x24}[t][t]{\scriptsize 6}%
%
\psfrag{v01}[r][r]{\scriptsize 0}%
\psfrag{v02}[r][r]{\scriptsize 0.1}%
\psfrag{v03}[r][r]{\scriptsize 0.2}%
\psfrag{v04}[r][r]{\scriptsize 0.3}%
\psfrag{v05}[r][r]{\scriptsize 0.4}%
\psfrag{v06}[r][r]{\scriptsize 0.5}%
\psfrag{v07}[r][r]{\scriptsize 0.6}%
\psfrag{v08}[r][r]{\scriptsize 0.7}%
\psfrag{v09}[r][r]{\scriptsize 0.8}%
\psfrag{v10}[r][r]{\scriptsize 0.9}%
\psfrag{v11}[r][r]{\scriptsize 1}%
\psfrag{v12}[r][r]{\scriptsize 0}%
\psfrag{v13}[r][r]{\scriptsize 0.2}%
\psfrag{v14}[r][r]{\scriptsize 0.4}%
\psfrag{v15}[r][r]{\scriptsize 0.6}%
\psfrag{v16}[r][r]{\scriptsize 0.8}%
\psfrag{v17}[r][r]{\scriptsize 1}%
%
\includegraphics[scale=0.55]{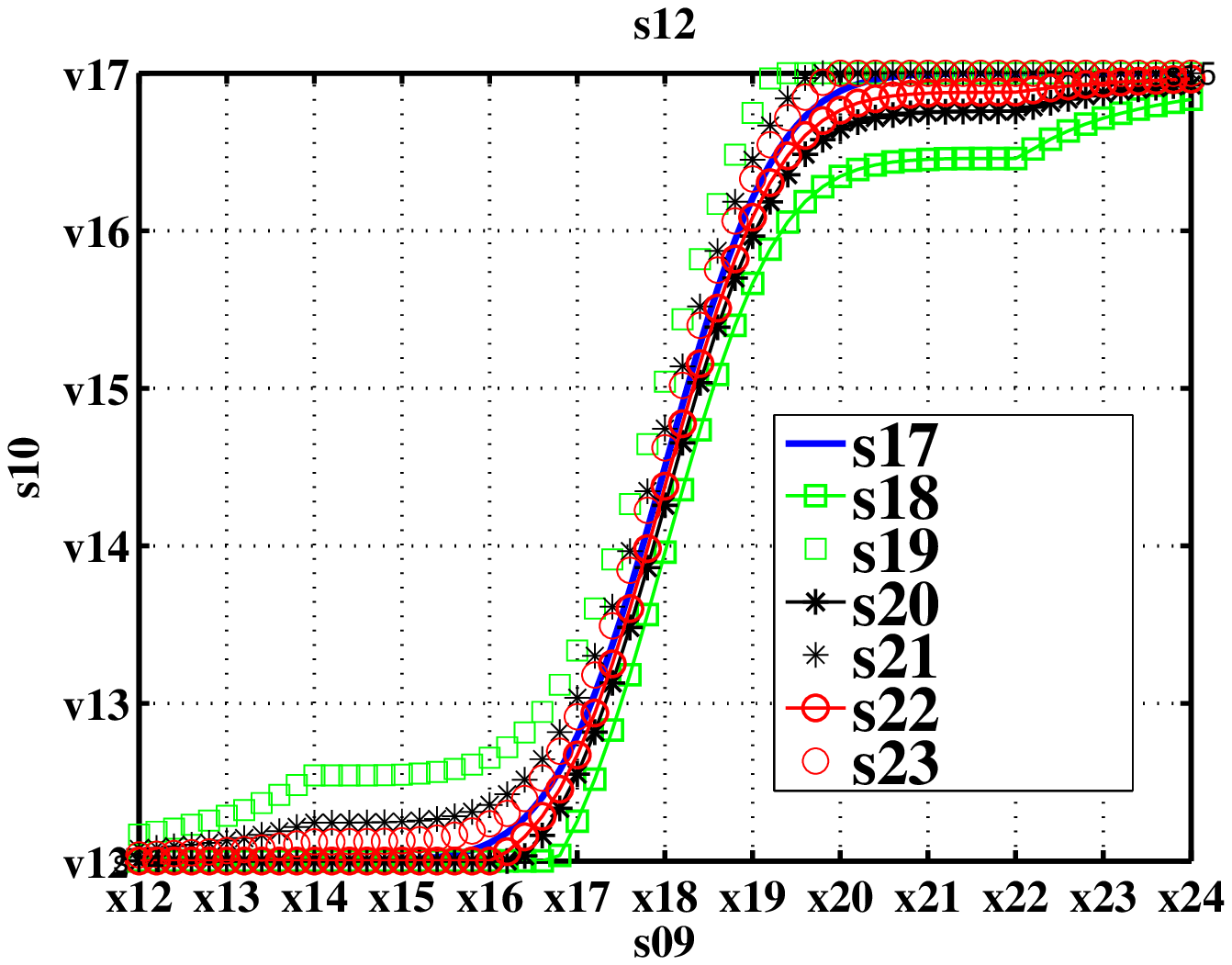}%
\end{psfrags}%
%

%
\hspace{\fill}
%
%
%
\begin{psfrags}%
\psfragscanon%
%
\psfrag{s09}[t][t]{\color[rgb]{0,0,0}\setlength{\tabcolsep}{0pt}\scriptsize \begin{tabular}{c}Centered and Normalized Interference Power\end{tabular}}%
\psfrag{s10}[b][b]{\color[rgb]{0,0,0}\setlength{\tabcolsep}{0pt}\scriptsize \begin{tabular}{c}Normal CDF and Bounds\end{tabular}}%
\psfrag{s12}[b][b]{\color[rgb]{0,0,0}\setlength{\tabcolsep}{0pt}\scriptsize \begin{tabular}{c}$G_2(t) = \frac{1}{1+t^{4}}$ (without fading) \end{tabular}}%
\psfrag{s14}[][]{\color[rgb]{0,0,0}\setlength{\tabcolsep}{0pt}\begin{tabular}{c} $ $ \end{tabular}}%
\psfrag{s15}[][]{\color[rgb]{0,0,0}\setlength{\tabcolsep}{0pt}\begin{tabular}{c} $ $ \end{tabular}}%
\psfrag{s17}[l][l]{\color[rgb]{0,0,0}\scriptsize Normal}%
\psfrag{s18}[l][l]{\color[rgb]{0,0,0}\scriptsize $\lambda = 5$ (LB)}%
\psfrag{s19}[l][l]{\color[rgb]{0,0,0}\scriptsize $\lambda = 5$ (UB)}%
\psfrag{s20}[l][l]{\color[rgb]{0,0,0}\scriptsize $\lambda = 25$ (LB)}%
\psfrag{s21}[l][l]{\color[rgb]{0,0,0}\scriptsize $\lambda = 25$ (UB)}%
\psfrag{s22}[l][l]{\color[rgb]{0,0,0}\scriptsize $\lambda = 100$ (LB)}%
\psfrag{s23}[l][l]{\color[rgb]{0,0,0}\scriptsize $\lambda = 100$ (UB)}%
%
\psfrag{x01}[t][t]{\scriptsize 0}%
\psfrag{x02}[t][t]{\scriptsize 0.1}%
\psfrag{x03}[t][t]{\scriptsize 0.2}%
\psfrag{x04}[t][t]{\scriptsize 0.3}%
\psfrag{x05}[t][t]{\scriptsize 0.4}%
\psfrag{x06}[t][t]{\scriptsize 0.5}%
\psfrag{x07}[t][t]{\scriptsize 0.6}%
\psfrag{x08}[t][t]{\scriptsize 0.7}%
\psfrag{x09}[t][t]{\scriptsize 0.8}%
\psfrag{x10}[t][t]{\scriptsize 0.9}%
\psfrag{x11}[t][t]{\scriptsize 1}%
\psfrag{x12}[t][t]{\scriptsize -6}%
\psfrag{x13}[t][t]{\scriptsize -5}%
\psfrag{x14}[t][t]{\scriptsize -4}%
\psfrag{x15}[t][t]{\scriptsize -3}%
\psfrag{x16}[t][t]{\scriptsize -2}%
\psfrag{x17}[t][t]{\scriptsize -1}%
\psfrag{x18}[t][t]{\scriptsize 0}%
\psfrag{x19}[t][t]{\scriptsize 1}%
\psfrag{x20}[t][t]{\scriptsize 2}%
\psfrag{x21}[t][t]{\scriptsize 3}%
\psfrag{x22}[t][t]{\scriptsize 4}%
\psfrag{x23}[t][t]{\scriptsize 5}%
\psfrag{x24}[t][t]{\scriptsize 6}%
%
\psfrag{v01}[r][r]{\scriptsize 0}%
\psfrag{v02}[r][r]{\scriptsize 0.1}%
\psfrag{v03}[r][r]{\scriptsize 0.2}%
\psfrag{v04}[r][r]{\scriptsize 0.3}%
\psfrag{v05}[r][r]{\scriptsize 0.4}%
\psfrag{v06}[r][r]{\scriptsize 0.5}%
\psfrag{v07}[r][r]{\scriptsize 0.6}%
\psfrag{v08}[r][r]{\scriptsize 0.7}%
\psfrag{v09}[r][r]{\scriptsize 0.8}%
\psfrag{v10}[r][r]{\scriptsize 0.9}%
\psfrag{v11}[r][r]{\scriptsize 1}%
\psfrag{v12}[r][r]{\scriptsize 0}%
\psfrag{v13}[r][r]{\scriptsize 0.2}%
\psfrag{v14}[r][r]{\scriptsize 0.4}%
\psfrag{v15}[r][r]{\scriptsize 0.6}%
\psfrag{v16}[r][r]{\scriptsize 0.8}%
\psfrag{v17}[r][r]{\scriptsize 1}%
%
\includegraphics[scale=0.55]{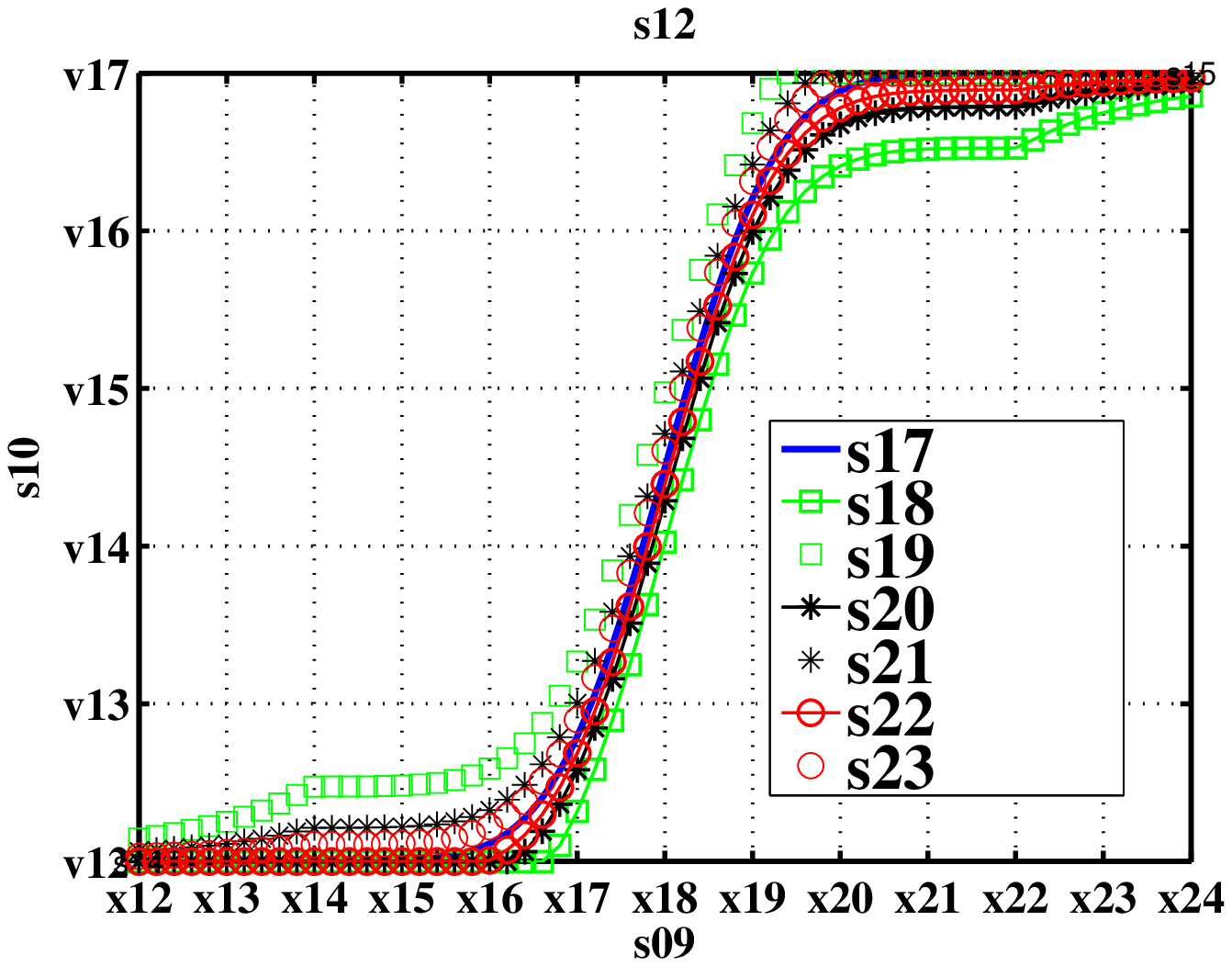}%
\end{psfrags}%
%

\end{center}
\end{minipage} 
\\ 
\begin{minipage}[t]{\textwidth}
\vspace{-2mm} 
\begin{center}
%
%
\begin{psfrags}%
\psfragscanon%
%
\psfrag{s09}[t][t]{\color[rgb]{0,0,0}\setlength{\tabcolsep}{0pt}\scriptsize \begin{tabular}{c}Centered and Normalized Interference Power\end{tabular}}%
\psfrag{s10}[b][b]{\color[rgb]{0,0,0}\setlength{\tabcolsep}{0pt}\scriptsize \begin{tabular}{c}Normal CDF and Bounds\end{tabular}}%
\psfrag{s12}[b][b]{\color[rgb]{0,0,0}\setlength{\tabcolsep}{0pt}\scriptsize \begin{tabular}{c}$G_1(t) = \frac{1}{(1+t)^4}$ (with fading)\end{tabular}}%
\psfrag{s14}[][]{\color[rgb]{0,0,0}\setlength{\tabcolsep}{0pt}\begin{tabular}{c} $ $ \end{tabular}}%
\psfrag{s15}[][]{\color[rgb]{0,0,0}\setlength{\tabcolsep}{0pt}\begin{tabular}{c} $ $ \end{tabular}}%
\psfrag{s17}[l][l]{\color[rgb]{0,0,0}\scriptsize Normal}%
\psfrag{s18}[l][l]{\color[rgb]{0,0,0}\scriptsize $\lambda = 5$ (LB)}%
\psfrag{s19}[l][l]{\color[rgb]{0,0,0}\scriptsize $\lambda = 5$ (UB)}%
\psfrag{s20}[l][l]{\color[rgb]{0,0,0}\scriptsize $\lambda = 25$ (LB)}%
\psfrag{s21}[l][l]{\color[rgb]{0,0,0}\scriptsize $\lambda = 25$ (UB)}%
\psfrag{s22}[l][l]{\color[rgb]{0,0,0}\scriptsize $\lambda = 100$ (LB)}%
\psfrag{s23}[l][l]{\color[rgb]{0,0,0}\scriptsize $\lambda = 100$ (UB)}%
%
\psfrag{x01}[t][t]{\scriptsize 0}%
\psfrag{x02}[t][t]{\scriptsize 0.1}%
\psfrag{x03}[t][t]{\scriptsize 0.2}%
\psfrag{x04}[t][t]{\scriptsize 0.3}%
\psfrag{x05}[t][t]{\scriptsize 0.4}%
\psfrag{x06}[t][t]{\scriptsize 0.5}%
\psfrag{x07}[t][t]{\scriptsize 0.6}%
\psfrag{x08}[t][t]{\scriptsize 0.7}%
\psfrag{x09}[t][t]{\scriptsize 0.8}%
\psfrag{x10}[t][t]{\scriptsize 0.9}%
\psfrag{x11}[t][t]{\scriptsize 1}%
\psfrag{x12}[t][t]{\scriptsize -6}%
\psfrag{x13}[t][t]{\scriptsize -5}%
\psfrag{x14}[t][t]{\scriptsize -4}%
\psfrag{x15}[t][t]{\scriptsize -3}%
\psfrag{x16}[t][t]{\scriptsize -2}%
\psfrag{x17}[t][t]{\scriptsize -1}%
\psfrag{x18}[t][t]{\scriptsize 0}%
\psfrag{x19}[t][t]{\scriptsize 1}%
\psfrag{x20}[t][t]{\scriptsize 2}%
\psfrag{x21}[t][t]{\scriptsize 3}%
\psfrag{x22}[t][t]{\scriptsize 4}%
\psfrag{x23}[t][t]{\scriptsize 5}%
\psfrag{x24}[t][t]{\scriptsize 6}%
%
\psfrag{v01}[r][r]{\scriptsize 0}%
\psfrag{v02}[r][r]{\scriptsize 0.1}%
\psfrag{v03}[r][r]{\scriptsize 0.2}%
\psfrag{v04}[r][r]{\scriptsize 0.3}%
\psfrag{v05}[r][r]{\scriptsize 0.4}%
\psfrag{v06}[r][r]{\scriptsize 0.5}%
\psfrag{v07}[r][r]{\scriptsize 0.6}%
\psfrag{v08}[r][r]{\scriptsize 0.7}%
\psfrag{v09}[r][r]{\scriptsize 0.8}%
\psfrag{v10}[r][r]{\scriptsize 0.9}%
\psfrag{v11}[r][r]{\scriptsize 1}%
\psfrag{v12}[r][r]{\scriptsize 0}%
\psfrag{v13}[r][r]{\scriptsize 0.2}%
\psfrag{v14}[r][r]{\scriptsize 0.4}%
\psfrag{v15}[r][r]{\scriptsize 0.6}%
\psfrag{v16}[r][r]{\scriptsize 0.8}%
\psfrag{v17}[r][r]{\scriptsize 1}%
%
\includegraphics[scale=0.55]{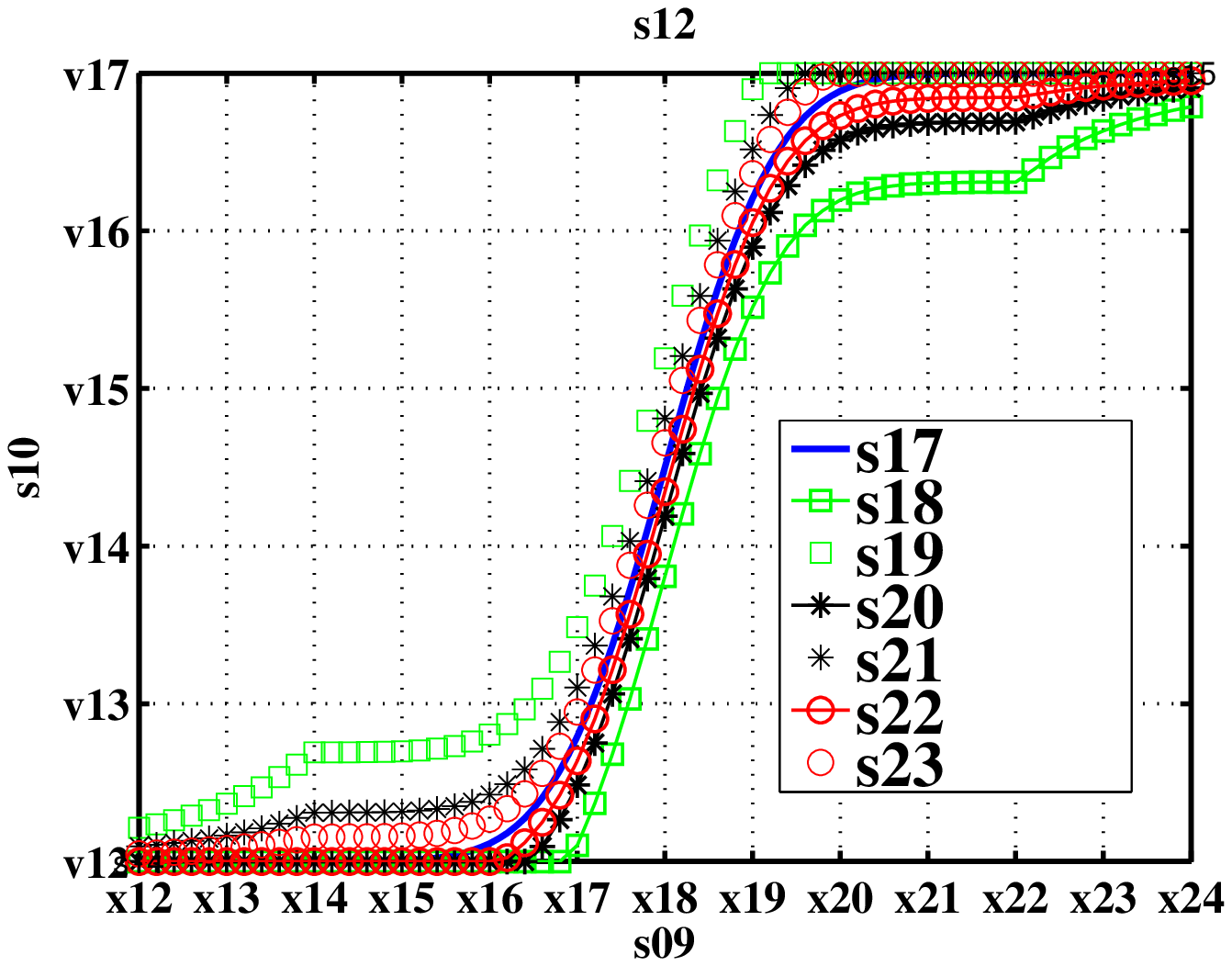}%
\end{psfrags}%
%

%
\hfill
%
%
%
\begin{psfrags}%
\psfragscanon%
%
\psfrag{s09}[t][t]{\color[rgb]{0,0,0}\setlength{\tabcolsep}{0pt}\scriptsize \begin{tabular}{c}Centered and Normalized Interference Power\end{tabular}}%
\psfrag{s10}[b][b]{\color[rgb]{0,0,0}\setlength{\tabcolsep}{0pt}\scriptsize \begin{tabular}{c}Normal CDF and Bounds\end{tabular}}%
\psfrag{s12}[b][b]{\color[rgb]{0,0,0}\setlength{\tabcolsep}{0pt}\scriptsize \begin{tabular}{c}$G_2(t) = \frac{1}{1+t^{4}}$ (with fading)\end{tabular}}%
\psfrag{s14}[][]{\color[rgb]{0,0,0}\setlength{\tabcolsep}{0pt}\begin{tabular}{c} $ $ \end{tabular}}%
\psfrag{s15}[][]{\color[rgb]{0,0,0}\setlength{\tabcolsep}{0pt}\begin{tabular}{c} $ $ \end{tabular}}%
\psfrag{s17}[l][l]{\color[rgb]{0,0,0}\scriptsize Normal}%
\psfrag{s18}[l][l]{\color[rgb]{0,0,0}\scriptsize $\lambda = 5$ (LB)}%
\psfrag{s19}[l][l]{\color[rgb]{0,0,0}\scriptsize $\lambda = 5$ (UB)}%
\psfrag{s20}[l][l]{\color[rgb]{0,0,0}\scriptsize $\lambda = 25$ (LB)}%
\psfrag{s21}[l][l]{\color[rgb]{0,0,0}\scriptsize $\lambda = 25$ (UB)}%
\psfrag{s22}[l][l]{\color[rgb]{0,0,0}\scriptsize $\lambda = 100$ (LB)}%
\psfrag{s23}[l][l]{\color[rgb]{0,0,0}\scriptsize $\lambda = 100$ (UB)}%
%
\psfrag{x01}[t][t]{\scriptsize 0}%
\psfrag{x02}[t][t]{\scriptsize 0.1}%
\psfrag{x03}[t][t]{\scriptsize 0.2}%
\psfrag{x04}[t][t]{\scriptsize 0.3}%
\psfrag{x05}[t][t]{\scriptsize 0.4}%
\psfrag{x06}[t][t]{\scriptsize 0.5}%
\psfrag{x07}[t][t]{\scriptsize 0.6}%
\psfrag{x08}[t][t]{\scriptsize 0.7}%
\psfrag{x09}[t][t]{\scriptsize 0.8}%
\psfrag{x10}[t][t]{\scriptsize 0.9}%
\psfrag{x11}[t][t]{\scriptsize 1}%
\psfrag{x12}[t][t]{\scriptsize -6}%
\psfrag{x13}[t][t]{\scriptsize -5}%
\psfrag{x14}[t][t]{\scriptsize -4}%
\psfrag{x15}[t][t]{\scriptsize -3}%
\psfrag{x16}[t][t]{\scriptsize -2}%
\psfrag{x17}[t][t]{\scriptsize -1}%
\psfrag{x18}[t][t]{\scriptsize 0}%
\psfrag{x19}[t][t]{\scriptsize 1}%
\psfrag{x20}[t][t]{\scriptsize 2}%
\psfrag{x21}[t][t]{\scriptsize 3}%
\psfrag{x22}[t][t]{\scriptsize 4}%
\psfrag{x23}[t][t]{\scriptsize 5}%
\psfrag{x24}[t][t]{6}%
%
\psfrag{v01}[r][r]{\scriptsize 0}%
\psfrag{v02}[r][r]{\scriptsize 0.1}%
\psfrag{v03}[r][r]{\scriptsize 0.2}%
\psfrag{v04}[r][r]{\scriptsize 0.3}%
\psfrag{v05}[r][r]{\scriptsize 0.4}%
\psfrag{v06}[r][r]{\scriptsize 0.5}%
\psfrag{v07}[r][r]{\scriptsize 0.6}%
\psfrag{v08}[r][r]{\scriptsize 0.7}%
\psfrag{v09}[r][r]{\scriptsize 0.8}%
\psfrag{v10}[r][r]{\scriptsize 0.9}%
\psfrag{v11}[r][r]{\scriptsize 1}%
\psfrag{v12}[r][r]{\scriptsize 0}%
\psfrag{v13}[r][r]{\scriptsize 0.2}%
\psfrag{v14}[r][r]{\scriptsize 0.4}%
\psfrag{v15}[r][r]{\scriptsize 0.6}%
\psfrag{v16}[r][r]{\scriptsize 0.8}%
\psfrag{v17}[r][r]{\scriptsize 1}%
%
\includegraphics[scale=0.55]{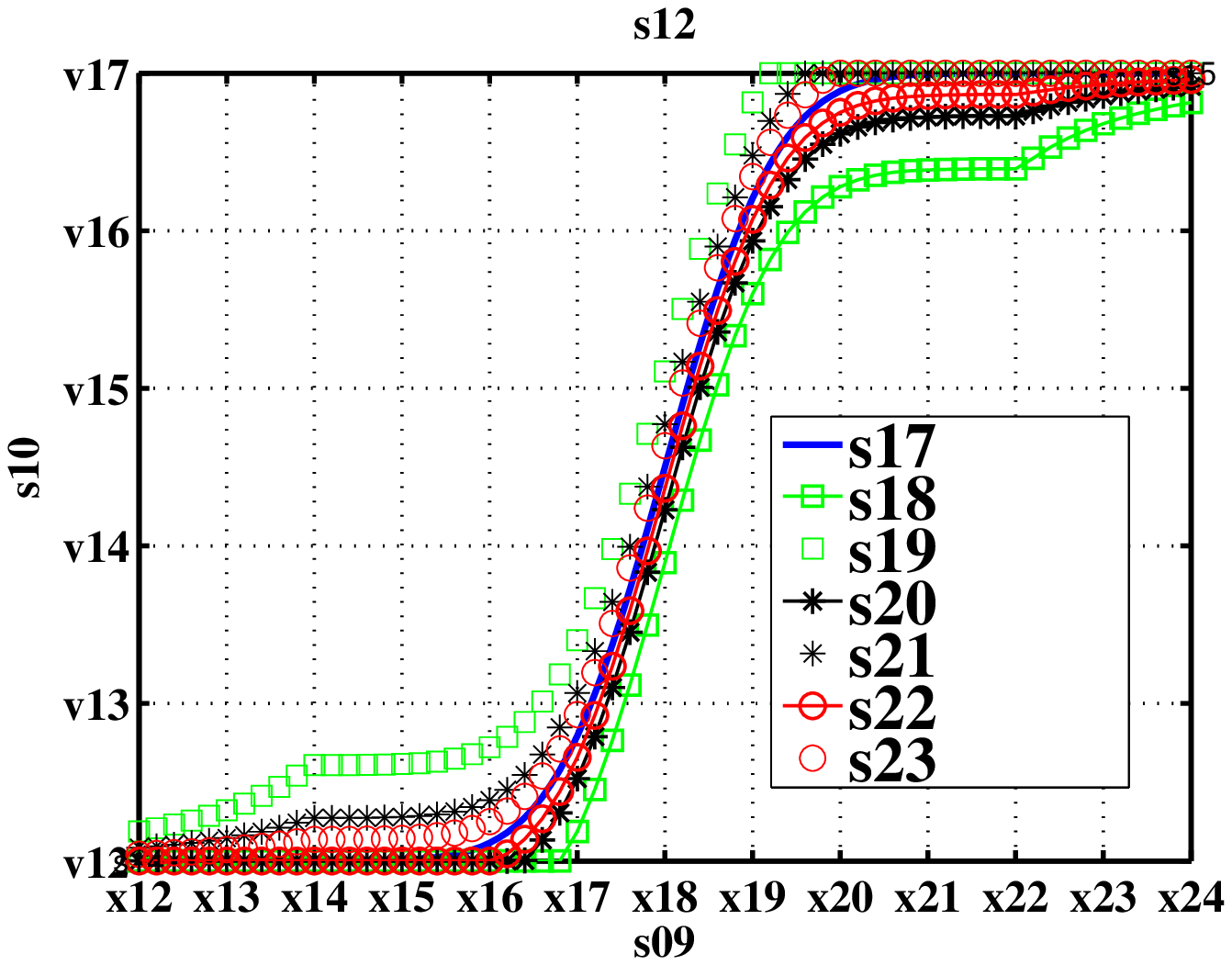}%
\end{psfrags}%
%

\end{center}
\end{minipage}
\caption{Upper and lower bounds on the centered and normalized WMAI CDFs for the path-loss functions $G_1(t) = \frac{1}{\paren{1+t}^\alpha}$ (lefthand side figures) and $G_2(t) = \frac{1}{1+t^\alpha}$ (righthand side figures).  Interferers are distributed over $\R^2$ according to a non-stationary PPP with mean measure density given as in \eqref{Eqn: Density 1 - Appendix}.  The effect of fading is also illustrated in the bottom figures by assuming Nakagami-$m$ fading with $m$ parameter set to $5$. ($\alpha = 4$ and $r=0.5$)} \label{Fig: CDF Bounds - Nonstationary}
\end{figure*}

\begin{theorem} \label{Thm: Rates of Convergence - Non-stationary}
Assume $\Phi_\Lambda$ is a PPP with mean measure density $f$ given as in \eqref{Eqn: Density 1 - Appendix}.  Then, for all $x \in \Re$,
\begin{eqnarray}
\left| \PR{\frac{I_\lambda - \ES{I_\lambda}}{\sqrt{\V{I_\lambda}}} \leq x} - \Psi(x)\right| \leq \frac{c(x)}{\sqrt{\lambda}}, \label{Eqn: Uniform Rates Non-Stationary PPP}
\end{eqnarray}
where $\Psi(x) = \frac{1}{\sqrt{2 \pi}}\int_{-\infty}^x \e{-\frac{t^2}{2}} dt$ and $c(x) = \frac{1}{\sqrt{2 \pi}}\frac{m_{H^3}}{\paren{m_{H^2}}^\frac32}\frac{\int_r^\infty G^3(t) \frac{1}{t} dt}{\paren{\int_r^\infty G^2(t)\frac{1}{t} dt}^\frac32} \min\paren{0.4785, \frac{31.935}{1+|x|^3}}$. 
\end{theorem}
\begin{IEEEproof}
Directly follows from Theorem \ref{Thm: Rates of Convergence} after substituting $\frac{2\pi}{t}\I{t\geq r}$ for $p(t)$.  
\end{IEEEproof} 

\begin{figure*}[!t]
\begin{minipage}[t]{\textwidth}
\begin{center}
%
%
\begin{psfrags}%
\psfragscanon%
%
\psfrag{s13}[t][t]{\color[rgb]{0,0,0}\setlength{\tabcolsep}{0pt}\scriptsize \begin{tabular}{c} $ $ \end{tabular}}%
\psfrag{s14}[b][b]{\color[rgb]{0,0,0}\setlength{\tabcolsep}{0pt}\scriptsize \begin{tabular}{c} \hspace{-25mm} Cumulative Distribution Function \end{tabular}}%
\psfrag{s16}[b][b]{\color[rgb]{0,0,0}\setlength{\tabcolsep}{0pt}\scriptsize \begin{tabular}{c} $G_1(t) = \frac{1}{(1+t)^\alpha}$ (without fading) \end{tabular}}%
\psfrag{s18}[][]{\color[rgb]{0,0,0}\setlength{\tabcolsep}{0pt}\scriptsize \begin{tabular}{c} $ $ \end{tabular}}%
\psfrag{s19}[][]{\color[rgb]{0,0,0}\setlength{\tabcolsep}{0pt}\scriptsize \begin{tabular}{c} $ $ \end{tabular}}%
\psfrag{s20}[l][l]{\color[rgb]{0,0,0}\scriptsize $\lambda = 10$}%
\psfrag{s21}[l][l]{\color[rgb]{0,0,0}\scriptsize Normal}%
\psfrag{s22}[l][l]{\color[rgb]{0,0,0}\scriptsize $\lambda = 0.1$}%
\psfrag{s23}[l][l]{\color[rgb]{0,0,0}\scriptsize $\lambda = 1$}%
\psfrag{s24}[l][l]{\color[rgb]{0,0,0}\scriptsize $\lambda = 10$}%
\psfrag{s25}[t][t]{\color[rgb]{0,0,0}\scriptsize \setlength{\tabcolsep}{0pt}\begin{tabular}{c} Centered and Normalized Interference Power \end{tabular}}%
\psfrag{s26}[b][b]{\color[rgb]{0,0,0}\setlength{\tabcolsep}{0pt}\scriptsize \begin{tabular}{c} $ $ \end{tabular}}%
\psfrag{s30}[][]{\color[rgb]{0,0,0}\setlength{\tabcolsep}{0pt}\scriptsize \begin{tabular}{c} $ $ \end{tabular}}%
\psfrag{s31}[][]{\color[rgb]{0,0,0}\setlength{\tabcolsep}{0pt}\scriptsize \begin{tabular}{c} $ $ \end{tabular}}%
\psfrag{s32}[l][l]{\color[rgb]{0,0,0}\scriptsize $\lambda = 10$}%
\psfrag{s33}[l][l]{\color[rgb]{0,0,0}\scriptsize Normal}%
\psfrag{s34}[l][l]{\color[rgb]{0,0,0}\scriptsize $\lambda = 0.1$}%
\psfrag{s35}[l][l]{\color[rgb]{0,0,0}\scriptsize $\lambda = 1$}%
\psfrag{s36}[l][l]{\color[rgb]{0,0,0}\scriptsize $\lambda = 10$}%
%
\psfrag{x01}[t][t]{\scriptsize 0}%
\psfrag{x02}[t][t]{\scriptsize 0.1}%
\psfrag{x03}[t][t]{\scriptsize 0.2}%
\psfrag{x04}[t][t]{\scriptsize 0.3}%
\psfrag{x05}[t][t]{\scriptsize 0.4}%
\psfrag{x06}[t][t]{\scriptsize 0.5}%
\psfrag{x07}[t][t]{\scriptsize 0.6}%
\psfrag{x08}[t][t]{\scriptsize 0.7}%
\psfrag{x09}[t][t]{\scriptsize 0.8}%
\psfrag{x10}[t][t]{\scriptsize 0.9}%
\psfrag{x11}[t][t]{\scriptsize 1}%
\psfrag{x12}[t][t]{\scriptsize -3}%
\psfrag{x13}[t][t]{\scriptsize -2}%
\psfrag{x14}[t][t]{\scriptsize -1}%
\psfrag{x15}[t][t]{\scriptsize 0}%
\psfrag{x16}[t][t]{\scriptsize 1}%
\psfrag{x17}[t][t]{\scriptsize 2}%
\psfrag{x18}[t][t]{\scriptsize 3}%
\psfrag{x19}[t][t]{\scriptsize -3}%
\psfrag{x20}[t][t]{\scriptsize -2}%
\psfrag{x21}[t][t]{\scriptsize -1}%
\psfrag{x22}[t][t]{\scriptsize 0}%
\psfrag{x23}[t][t]{\scriptsize 1}%
\psfrag{x24}[t][t]{\scriptsize 2}%
\psfrag{x25}[t][t]{\scriptsize 3}%
%
\psfrag{v01}[r][r]{\scriptsize 0}%
\psfrag{v02}[r][r]{\scriptsize 0.1}%
\psfrag{v03}[r][r]{\scriptsize 0.2}%
\psfrag{v04}[r][r]{\scriptsize 0.3}%
\psfrag{v05}[r][r]{\scriptsize 0.4}%
\psfrag{v06}[r][r]{\scriptsize 0.5}%
\psfrag{v07}[r][r]{\scriptsize 0.6}%
\psfrag{v08}[r][r]{\scriptsize 0.7}%
\psfrag{v09}[r][r]{\scriptsize 0.8}%
\psfrag{v10}[r][r]{\scriptsize 0.9}%
\psfrag{v11}[r][r]{\scriptsize 1}%
\psfrag{v12}[r][r]{\scriptsize 0}%
\psfrag{v13}[r][r]{\scriptsize 0.2}%
\psfrag{v14}[r][r]{\scriptsize 0.4}%
\psfrag{v15}[r][r]{\scriptsize 0.6}%
\psfrag{v16}[r][r]{\scriptsize 0.8}%
\psfrag{v17}[r][r]{\scriptsize 1}%
\psfrag{v18}[r][r]{\scriptsize 0}%
\psfrag{v19}[r][r]{\scriptsize 0.2}%
\psfrag{v20}[r][r]{\scriptsize 0.4}%
\psfrag{v21}[r][r]{\scriptsize 0.6}%
\psfrag{v22}[r][r]{\scriptsize 0.8}%
\psfrag{v23}[r][r]{\scriptsize 1}%
%
\setlength{\unitlength}{1cm}
\begin{picture}(8, 6)(0, 0)
\put(0, 0){\includegraphics[scale=0.55]{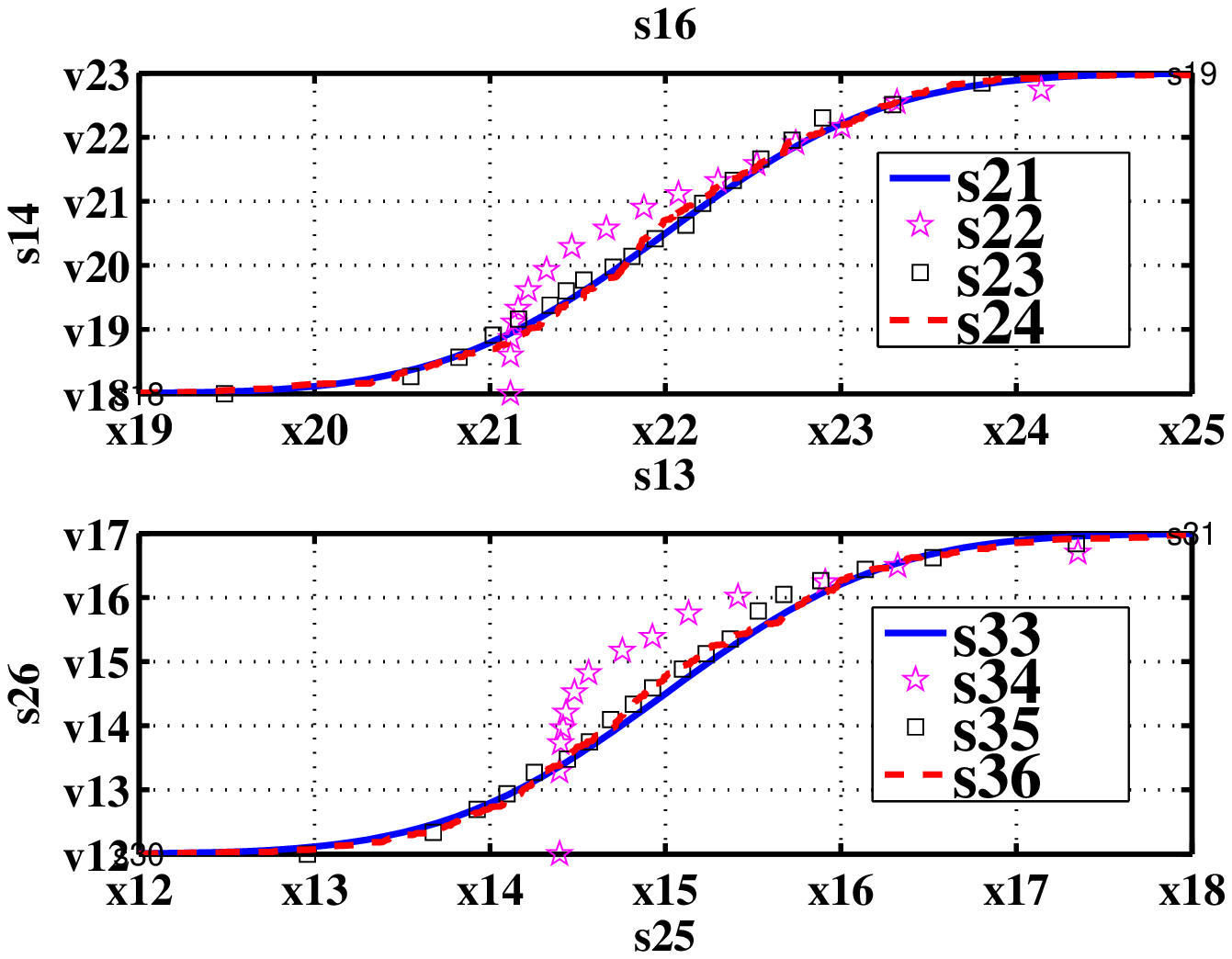}}%
\put(1.92, 4.35){\color[rgb]{0, 0, 0} \scriptsize \colorbox{white}{\fbox{$\alpha = 3$}}}
\put(1.92, 1.54){\color[rgb]{0, 0, 0} \scriptsize \colorbox{white}{\fbox{$\alpha = 5$}}}
\end{picture}
\end{psfrags}%
%

%
\hspace{\fill}
%
%
%
\begin{psfrags}%
\psfragscanon%
%
\psfrag{s13}[t][t]{\color[rgb]{0,0,0}\setlength{\tabcolsep}{0pt}\scriptsize \begin{tabular}{c} $ $ \end{tabular}}%
\psfrag{s14}[b][b]{\color[rgb]{0,0,0}\setlength{\tabcolsep}{0pt}\scriptsize \begin{tabular}{c}\hspace{-25mm} Cumulative Distribution Function \end{tabular}}%
\psfrag{s16}[b][b]{\color[rgb]{0,0,0}\setlength{\tabcolsep}{0pt}\scriptsize \begin{tabular}{c} $G_2(t) = \frac{1}{1+t^\alpha}$ (without fading) \end{tabular}}%
\psfrag{s18}[][]{\color[rgb]{0,0,0}\setlength{\tabcolsep}{0pt}\scriptsize \begin{tabular}{c} $ $ \end{tabular}}%
\psfrag{s19}[][]{\color[rgb]{0,0,0}\setlength{\tabcolsep}{0pt}\scriptsize \begin{tabular}{c} $ $ \end{tabular}}%
\psfrag{s20}[l][l]{\color[rgb]{0,0,0}\scriptsize $\lambda = 10$}%
\psfrag{s21}[l][l]{\color[rgb]{0,0,0}\scriptsize Normal}%
\psfrag{s22}[l][l]{\color[rgb]{0,0,0}\scriptsize $\lambda = 0.1$}%
\psfrag{s23}[l][l]{\color[rgb]{0,0,0}\scriptsize $\lambda = 1$}%
\psfrag{s24}[l][l]{\color[rgb]{0,0,0}\scriptsize $\lambda = 10$}%
\psfrag{s25}[t][t]{\color[rgb]{0,0,0}\setlength{\tabcolsep}{0pt}\scriptsize \begin{tabular}{c} Centered and Normalized Interference Power \end{tabular}}%
\psfrag{s26}[b][b]{\color[rgb]{0,0,0}\setlength{\tabcolsep}{0pt}\scriptsize \begin{tabular}{c} $ $ \end{tabular}}%
\psfrag{s30}[][]{\color[rgb]{0,0,0}\setlength{\tabcolsep}{0pt}\scriptsize \begin{tabular}{c} $ $ \end{tabular}}%
\psfrag{s31}[][]{\color[rgb]{0,0,0}\setlength{\tabcolsep}{0pt}\scriptsize \begin{tabular}{c} $ $ \end{tabular}}%
\psfrag{s32}[l][l]{\color[rgb]{0,0,0}\scriptsize $\lambda = 10$}%
\psfrag{s33}[l][l]{\color[rgb]{0,0,0}\scriptsize Normal}%
\psfrag{s34}[l][l]{\color[rgb]{0,0,0}\scriptsize $\lambda = 0.1$}%
\psfrag{s35}[l][l]{\color[rgb]{0,0,0}\scriptsize $\lambda = 1$}%
\psfrag{s36}[l][l]{\color[rgb]{0,0,0}\scriptsize $\lambda = 10$}%
%
\psfrag{x01}[t][t]{\scriptsize 0}%
\psfrag{x02}[t][t]{\scriptsize 0.1}%
\psfrag{x03}[t][t]{\scriptsize 0.2}%
\psfrag{x04}[t][t]{\scriptsize 0.3}%
\psfrag{x05}[t][t]{\scriptsize 0.4}%
\psfrag{x06}[t][t]{\scriptsize 0.5}%
\psfrag{x07}[t][t]{\scriptsize 0.6}%
\psfrag{x08}[t][t]{\scriptsize 0.7}%
\psfrag{x09}[t][t]{\scriptsize 0.8}%
\psfrag{x10}[t][t]{\scriptsize 0.9}%
\psfrag{x11}[t][t]{\scriptsize 1}%
\psfrag{x12}[t][t]{\scriptsize -3}%
\psfrag{x13}[t][t]{\scriptsize -2}%
\psfrag{x14}[t][t]{\scriptsize -1}%
\psfrag{x15}[t][t]{\scriptsize 0}%
\psfrag{x16}[t][t]{\scriptsize 1}%
\psfrag{x17}[t][t]{\scriptsize 2}%
\psfrag{x18}[t][t]{\scriptsize 3}%
\psfrag{x19}[t][t]{\scriptsize -3}%
\psfrag{x20}[t][t]{\scriptsize -2}%
\psfrag{x21}[t][t]{\scriptsize -1}%
\psfrag{x22}[t][t]{\scriptsize 0}%
\psfrag{x23}[t][t]{\scriptsize 1}%
\psfrag{x24}[t][t]{\scriptsize 2}%
\psfrag{x25}[t][t]{\scriptsize 3}%
%
\psfrag{v01}[r][r]{\scriptsize 0}%
\psfrag{v02}[r][r]{\scriptsize 0.1}%
\psfrag{v03}[r][r]{\scriptsize 0.2}%
\psfrag{v04}[r][r]{\scriptsize 0.3}%
\psfrag{v05}[r][r]{\scriptsize 0.4}%
\psfrag{v06}[r][r]{\scriptsize 0.5}%
\psfrag{v07}[r][r]{\scriptsize 0.6}%
\psfrag{v08}[r][r]{\scriptsize 0.7}%
\psfrag{v09}[r][r]{\scriptsize 0.8}%
\psfrag{v10}[r][r]{\scriptsize 0.9}%
\psfrag{v11}[r][r]{\scriptsize 1}%
\psfrag{v12}[r][r]{\scriptsize 0}%
\psfrag{v13}[r][r]{\scriptsize 0.2}%
\psfrag{v14}[r][r]{\scriptsize 0.4}%
\psfrag{v15}[r][r]{\scriptsize 0.6}%
\psfrag{v16}[r][r]{\scriptsize 0.8}%
\psfrag{v17}[r][r]{\scriptsize 1}%
\psfrag{v18}[r][r]{\scriptsize 0}%
\psfrag{v19}[r][r]{\scriptsize 0.2}%
\psfrag{v20}[r][r]{\scriptsize 0.4}%
\psfrag{v21}[r][r]{\scriptsize 0.6}%
\psfrag{v22}[r][r]{\scriptsize 0.8}%
\psfrag{v23}[r][r]{\scriptsize 1}%
%
\setlength{\unitlength}{1cm}
\begin{picture}(8, 6)(0, 0)
\includegraphics[scale=0.55]{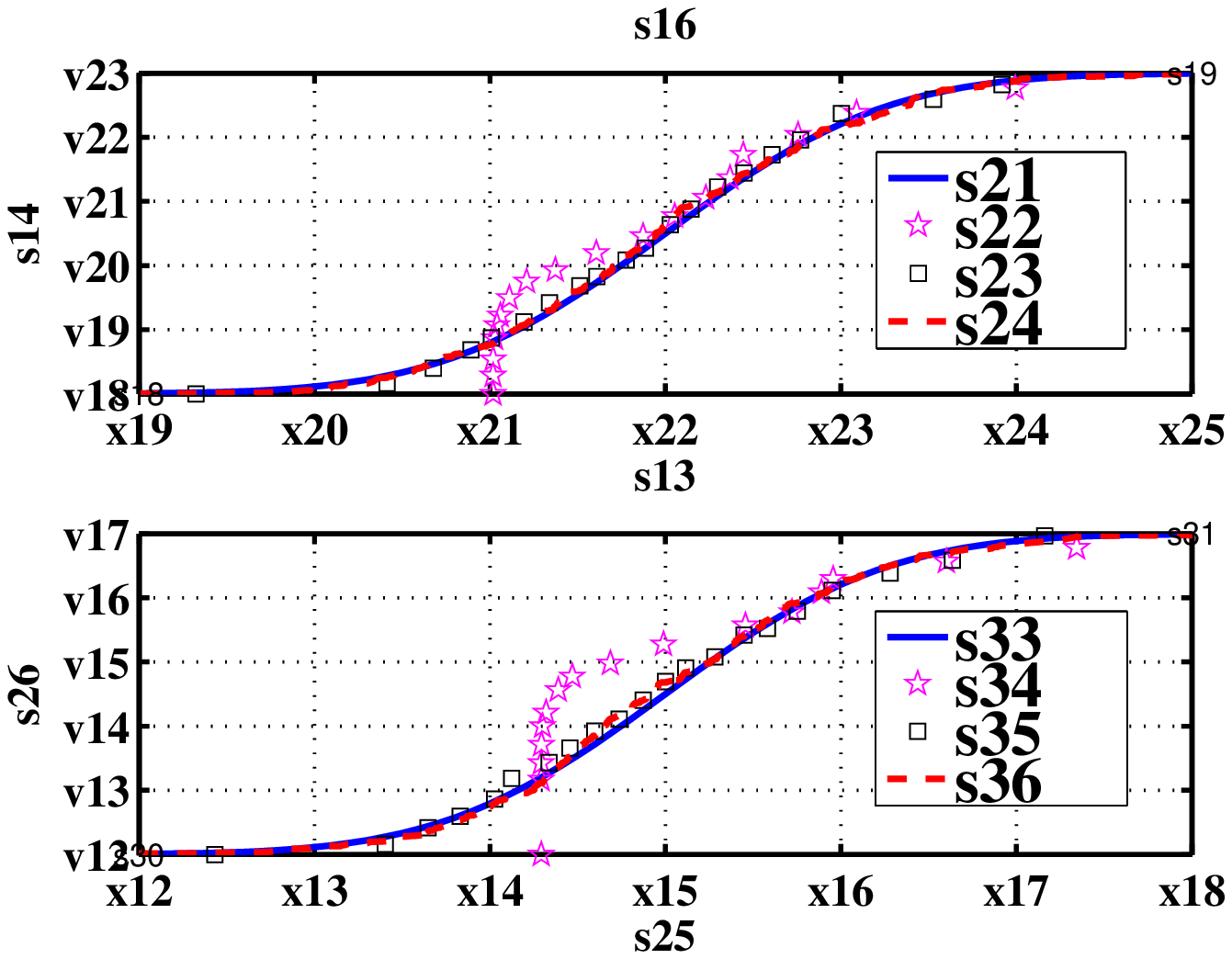}%
\put(-5.67, 4.35){\color[rgb]{0, 0, 0} \scriptsize \colorbox{white}{\fbox{$\alpha = 3$}}}
\put(-5.67, 1.54){\color[rgb]{0, 0, 0} \scriptsize \colorbox{white}{\fbox{$\alpha = 5$}}}
\end{picture}
\end{psfrags}%
%

\end{center}
\end{minipage} 
\\ 
\begin{minipage}[t]{\textwidth}
\vspace{-2mm} 
\begin{center}
%
%
\begin{psfrags}%
\psfragscanon%
%
\psfrag{s13}[t][t]{\color[rgb]{0,0,0}\setlength{\tabcolsep}{0pt}\scriptsize \begin{tabular}{c} $ $ \end{tabular}}%
\psfrag{s14}[b][b]{\color[rgb]{0,0,0}\setlength{\tabcolsep}{0pt}\scriptsize \begin{tabular}{c} \hspace{-25mm} Cumulative Distribution Function \end{tabular}}%
\psfrag{s16}[b][b]{\color[rgb]{0,0,0}\setlength{\tabcolsep}{0pt}\scriptsize \begin{tabular}{c} $G_1(t) = \frac{1}{(1+t)^\alpha}$ (with fading) \end{tabular}}%
\psfrag{s18}[][]{\color[rgb]{0,0,0}\setlength{\tabcolsep}{0pt}\scriptsize \begin{tabular}{c} $ $ \end{tabular}}%
\psfrag{s19}[][]{\color[rgb]{0,0,0}\setlength{\tabcolsep}{0pt}\scriptsize \begin{tabular}{c} $ $ \end{tabular}}%
\psfrag{s20}[l][l]{\color[rgb]{0,0,0}\scriptsize $\lambda = 10$}%
\psfrag{s21}[l][l]{\color[rgb]{0,0,0}\scriptsize Normal}%
\psfrag{s22}[l][l]{\color[rgb]{0,0,0}\scriptsize $\lambda = 0.1$}%
\psfrag{s23}[l][l]{\color[rgb]{0,0,0}\scriptsize $\lambda = 1$}%
\psfrag{s24}[l][l]{\color[rgb]{0,0,0}\scriptsize $\lambda = 10$}%
\psfrag{s25}[t][t]{\color[rgb]{0,0,0}\scriptsize \setlength{\tabcolsep}{0pt}\begin{tabular}{c} Centered and Normalized Interference Power \end{tabular}}%
\psfrag{s26}[b][b]{\color[rgb]{0,0,0}\setlength{\tabcolsep}{0pt}\scriptsize \begin{tabular}{c} $ $ \end{tabular}}%
\psfrag{s30}[][]{\color[rgb]{0,0,0}\setlength{\tabcolsep}{0pt}\scriptsize \begin{tabular}{c} $ $ \end{tabular}}%
\psfrag{s31}[][]{\color[rgb]{0,0,0}\setlength{\tabcolsep}{0pt}\scriptsize \begin{tabular}{c} $ $ \end{tabular}}%
\psfrag{s32}[l][l]{\color[rgb]{0,0,0}\scriptsize $\lambda = 10$}%
\psfrag{s33}[l][l]{\color[rgb]{0,0,0}\scriptsize Normal}%
\psfrag{s34}[l][l]{\color[rgb]{0,0,0}\scriptsize $\lambda = 0.1$}%
\psfrag{s35}[l][l]{\color[rgb]{0,0,0}\scriptsize $\lambda = 1$}%
\psfrag{s36}[l][l]{\color[rgb]{0,0,0}\scriptsize $\lambda = 10$}%
%
\psfrag{x01}[t][t]{\scriptsize 0}%
\psfrag{x02}[t][t]{\scriptsize 0.1}%
\psfrag{x03}[t][t]{\scriptsize 0.2}%
\psfrag{x04}[t][t]{\scriptsize 0.3}%
\psfrag{x05}[t][t]{\scriptsize 0.4}%
\psfrag{x06}[t][t]{\scriptsize 0.5}%
\psfrag{x07}[t][t]{\scriptsize 0.6}%
\psfrag{x08}[t][t]{\scriptsize 0.7}%
\psfrag{x09}[t][t]{\scriptsize 0.8}%
\psfrag{x10}[t][t]{\scriptsize 0.9}%
\psfrag{x11}[t][t]{\scriptsize 1}%
\psfrag{x12}[t][t]{\scriptsize -3}%
\psfrag{x13}[t][t]{\scriptsize -2}%
\psfrag{x14}[t][t]{\scriptsize -1}%
\psfrag{x15}[t][t]{\scriptsize 0}%
\psfrag{x16}[t][t]{\scriptsize 1}%
\psfrag{x17}[t][t]{\scriptsize 2}%
\psfrag{x18}[t][t]{\scriptsize 3}%
\psfrag{x19}[t][t]{\scriptsize -3}%
\psfrag{x20}[t][t]{\scriptsize -2}%
\psfrag{x21}[t][t]{\scriptsize -1}%
\psfrag{x22}[t][t]{\scriptsize 0}%
\psfrag{x23}[t][t]{\scriptsize 1}%
\psfrag{x24}[t][t]{\scriptsize 2}%
\psfrag{x25}[t][t]{\scriptsize 3}%
%
\psfrag{v01}[r][r]{\scriptsize 0}%
\psfrag{v02}[r][r]{\scriptsize 0.1}%
\psfrag{v03}[r][r]{\scriptsize 0.2}%
\psfrag{v04}[r][r]{\scriptsize 0.3}%
\psfrag{v05}[r][r]{\scriptsize 0.4}%
\psfrag{v06}[r][r]{\scriptsize 0.5}%
\psfrag{v07}[r][r]{\scriptsize 0.6}%
\psfrag{v08}[r][r]{\scriptsize 0.7}%
\psfrag{v09}[r][r]{\scriptsize 0.8}%
\psfrag{v10}[r][r]{\scriptsize 0.9}%
\psfrag{v11}[r][r]{\scriptsize 1}%
\psfrag{v12}[r][r]{\scriptsize 0}%
\psfrag{v13}[r][r]{\scriptsize 0.2}%
\psfrag{v14}[r][r]{\scriptsize 0.4}%
\psfrag{v15}[r][r]{\scriptsize 0.6}%
\psfrag{v16}[r][r]{\scriptsize 0.8}%
\psfrag{v17}[r][r]{\scriptsize 1}%
\psfrag{v18}[r][r]{\scriptsize 0}%
\psfrag{v19}[r][r]{\scriptsize 0.2}%
\psfrag{v20}[r][r]{\scriptsize 0.4}%
\psfrag{v21}[r][r]{\scriptsize 0.6}%
\psfrag{v22}[r][r]{\scriptsize 0.8}%
\psfrag{v23}[r][r]{\scriptsize 1}%
%
\setlength{\unitlength}{1cm}
\begin{picture}(8, 6)(0, 0)
\put(0, 0){\includegraphics[scale=0.55]{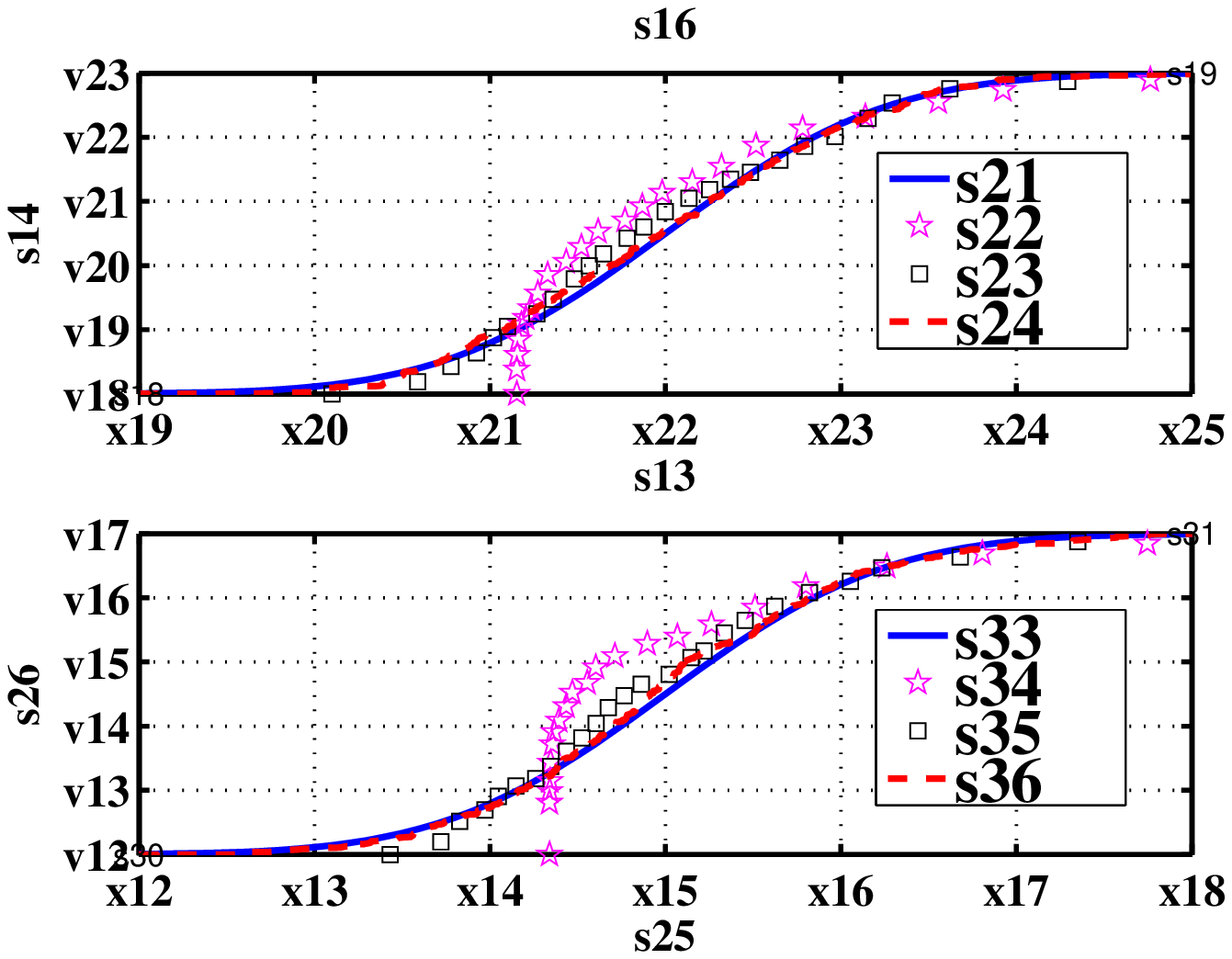}}%
\put(1.92, 4.35){\color[rgb]{0, 0, 0} \scriptsize \colorbox{white}{\fbox{$\alpha = 3$}}}
\put(1.92, 1.54){\color[rgb]{0, 0, 0} \scriptsize \colorbox{white}{\fbox{$\alpha = 5$}}}
\end{picture}
\end{psfrags}%
%

%
\hfill
%
%
%
\begin{psfrags}%
\psfragscanon%
%
\psfrag{s13}[t][t]{\color[rgb]{0,0,0}\setlength{\tabcolsep}{0pt}\scriptsize \begin{tabular}{c} $ $ \end{tabular}}%
\psfrag{s14}[b][b]{\color[rgb]{0,0,0}\setlength{\tabcolsep}{0pt}\scriptsize \begin{tabular}{c}\hspace{-25mm} Cumulative Distribution Function \end{tabular}}%
\psfrag{s16}[b][b]{\color[rgb]{0,0,0}\setlength{\tabcolsep}{0pt}\scriptsize \begin{tabular}{c} $G_2(t) = \frac{1}{1+t^\alpha}$ (with fading) \end{tabular}}%
\psfrag{s18}[][]{\color[rgb]{0,0,0}\setlength{\tabcolsep}{0pt}\scriptsize \begin{tabular}{c} $ $ \end{tabular}}%
\psfrag{s19}[][]{\color[rgb]{0,0,0}\setlength{\tabcolsep}{0pt}\scriptsize \begin{tabular}{c} $ $ \end{tabular}}%
\psfrag{s20}[l][l]{\color[rgb]{0,0,0}\scriptsize $\lambda = 10$}%
\psfrag{s21}[l][l]{\color[rgb]{0,0,0}\scriptsize Normal}%
\psfrag{s22}[l][l]{\color[rgb]{0,0,0}\scriptsize $\lambda = 0.1$}%
\psfrag{s23}[l][l]{\color[rgb]{0,0,0}\scriptsize $\lambda = 1$}%
\psfrag{s24}[l][l]{\color[rgb]{0,0,0}\scriptsize $\lambda = 10$}%
\psfrag{s25}[t][t]{\color[rgb]{0,0,0}\setlength{\tabcolsep}{0pt}\scriptsize \begin{tabular}{c} Centered and Normalized Interference Power \end{tabular}}%
\psfrag{s26}[b][b]{\color[rgb]{0,0,0}\setlength{\tabcolsep}{0pt}\scriptsize \begin{tabular}{c} $ $ \end{tabular}}%
\psfrag{s30}[][]{\color[rgb]{0,0,0}\setlength{\tabcolsep}{0pt}\scriptsize \begin{tabular}{c} $ $ \end{tabular}}%
\psfrag{s31}[][]{\color[rgb]{0,0,0}\setlength{\tabcolsep}{0pt}\scriptsize \begin{tabular}{c} $ $ \end{tabular}}%
\psfrag{s32}[l][l]{\color[rgb]{0,0,0}\scriptsize $\lambda = 10$}%
\psfrag{s33}[l][l]{\color[rgb]{0,0,0}\scriptsize Normal}%
\psfrag{s34}[l][l]{\color[rgb]{0,0,0}\scriptsize $\lambda = 0.1$}%
\psfrag{s35}[l][l]{\color[rgb]{0,0,0}\scriptsize $\lambda = 1$}%
\psfrag{s36}[l][l]{\color[rgb]{0,0,0}\scriptsize $\lambda = 10$}%
%
\psfrag{x01}[t][t]{\scriptsize 0}%
\psfrag{x02}[t][t]{\scriptsize 0.1}%
\psfrag{x03}[t][t]{\scriptsize 0.2}%
\psfrag{x04}[t][t]{\scriptsize 0.3}%
\psfrag{x05}[t][t]{\scriptsize 0.4}%
\psfrag{x06}[t][t]{\scriptsize 0.5}%
\psfrag{x07}[t][t]{\scriptsize 0.6}%
\psfrag{x08}[t][t]{\scriptsize 0.7}%
\psfrag{x09}[t][t]{\scriptsize 0.8}%
\psfrag{x10}[t][t]{\scriptsize 0.9}%
\psfrag{x11}[t][t]{\scriptsize 1}%
\psfrag{x12}[t][t]{\scriptsize -3}%
\psfrag{x13}[t][t]{\scriptsize -2}%
\psfrag{x14}[t][t]{\scriptsize -1}%
\psfrag{x15}[t][t]{\scriptsize 0}%
\psfrag{x16}[t][t]{\scriptsize 1}%
\psfrag{x17}[t][t]{\scriptsize 2}%
\psfrag{x18}[t][t]{\scriptsize 3}%
\psfrag{x19}[t][t]{\scriptsize -3}%
\psfrag{x20}[t][t]{\scriptsize -2}%
\psfrag{x21}[t][t]{\scriptsize -1}%
\psfrag{x22}[t][t]{\scriptsize 0}%
\psfrag{x23}[t][t]{\scriptsize 1}%
\psfrag{x24}[t][t]{\scriptsize 2}%
\psfrag{x25}[t][t]{\scriptsize 3}%
%
\psfrag{v01}[r][r]{\scriptsize 0}%
\psfrag{v02}[r][r]{\scriptsize 0.1}%
\psfrag{v03}[r][r]{\scriptsize 0.2}%
\psfrag{v04}[r][r]{\scriptsize 0.3}%
\psfrag{v05}[r][r]{\scriptsize 0.4}%
\psfrag{v06}[r][r]{\scriptsize 0.5}%
\psfrag{v07}[r][r]{\scriptsize 0.6}%
\psfrag{v08}[r][r]{\scriptsize 0.7}%
\psfrag{v09}[r][r]{\scriptsize 0.8}%
\psfrag{v10}[r][r]{\scriptsize 0.9}%
\psfrag{v11}[r][r]{\scriptsize 1}%
\psfrag{v12}[r][r]{\scriptsize 0}%
\psfrag{v13}[r][r]{\scriptsize 0.2}%
\psfrag{v14}[r][r]{\scriptsize 0.4}%
\psfrag{v15}[r][r]{\scriptsize 0.6}%
\psfrag{v16}[r][r]{\scriptsize 0.8}%
\psfrag{v17}[r][r]{\scriptsize 1}%
\psfrag{v18}[r][r]{\scriptsize 0}%
\psfrag{v19}[r][r]{\scriptsize 0.2}%
\psfrag{v20}[r][r]{\scriptsize 0.4}%
\psfrag{v21}[r][r]{\scriptsize 0.6}%
\psfrag{v22}[r][r]{\scriptsize 0.8}%
\psfrag{v23}[r][r]{\scriptsize 1}%
%
\setlength{\unitlength}{1cm}
\begin{picture}(8, 6)(0, 0)
\includegraphics[scale=0.55]{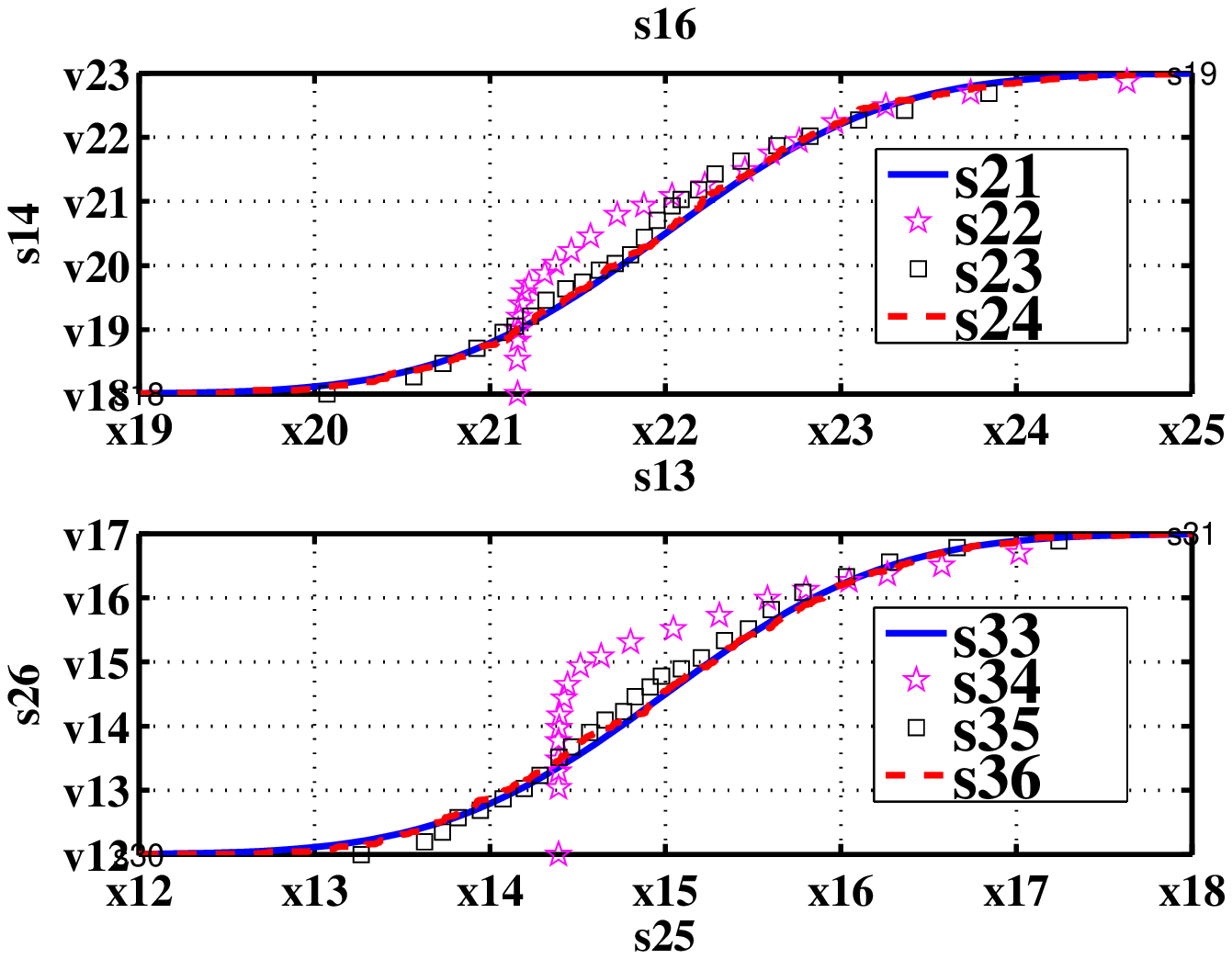}%
\put(-5.67, 4.35){\color[rgb]{0, 0, 0} \scriptsize \colorbox{white}{\fbox{$\alpha = 3$}}}
\put(-5.67, 1.54){\color[rgb]{0, 0, 0} \scriptsize \colorbox{white}{\fbox{$\alpha = 5$}}}
\end{picture}
\end{psfrags}%
%

\end{center}
\end{minipage}
\caption{Comparison of the simulated centered and normalized WMAI CDFs with the normal CDF for the path-loss functions $G_1(t) = \frac{1}{\paren{1+t}^\alpha}$ (lefthand side figures) and $G_2(t) = \frac{1}{1+t^\alpha}$ (righthand side figures).  Interferers are distributed over $\R^2$ according to a non-stationary PPP with mean measure density given as in \eqref{Eqn: Density 1 - Appendix}.  The effect of fading is also illustrated in the bottom figures by assuming Nakagami-$m$ fading with $m$ parameter set to $1$. ($r=0.5$)} \label{Fig: Simulated CDFs - Nonstationary}
\end{figure*}

In Figs. \ref{Fig: CDF Bounds - Nonstationary} and \ref{Fig: Simulated CDFs - Nonstationary}, we plot our numerically computed Gaussian approximation bounds and simulation results, respectively.  Since the key messages conveyed by these figures are similar to those explained in Section \ref{Section: Approximation Bounds for Stationary PPPs}, we do not repeat them here again.  However, several remarks are in order. Gaussian approximation bounds and simulation results given for $G_1(t)$ resemble to those given for $G_2(t)$ much more closely when compared to such bounds and simulation results given in Figs. \ref{Fig: CDF Bounds} and \ref{Fig: Simulated CDFs}.  This is because the path-loss dependent constants, {\em i.e.,} $\frac{\int_{r}^\infty G_i^3(t) \frac{1}{t} dt}{\paren{\int_{r}^\infty G_i^2(t) \frac{1}{t} dt}^\frac32}$, $i=1,2$, are more similar to each other with this particular choice of non-stationary spatial distribution of interfering transmitters over $\R^2$.  For example, this constant is equal to $1.27$ for $G_1(t)$ and equal to $1.11$ for $G_2(t)$ when $\alpha$ is set to $4$.  This also explains why we observe a significantly better statistical fit between the WMAI distribution and the normal distribution for $G_1(t)$ in this case, {\em i.e.,} see the corresponding constant in Table \ref{Table: 1}.  Finally, it is also noteworthy to mention that the WMAI distributions become closer to the Normal distribution for smaller values of $r$ as the path-loss dependent constants decrease when $r$ decreases.   

\section{Scaling Behavior of $C_{\lambda, {\rm outage}}(\gamma)$} \label{Appendix: Outage Probability Scaling}
In this appendix, we will briefly sketch the proof ideas leading to $C_{\lambda, {\rm outage}}\paren{\gamma} = \TO{\frac{1}{\lambda}}$ as $\lambda \ra \infty$.  Firstly, we observe that $\gamma\paren{h, R}$ is a non-increasing and continuous function of $h$ that is equal to one for $h \in \sqparen{0, \frac{\snr^{-1}\paren{\e{R}-1}}{G(d)}}$ and approaching to zero as $h$ grows to infinity for any fixed value of $R$ under our assumptions in Section \ref{Section: Network Model}.  Hence, we can find $h_\lambda^\star \in (0, \infty)$ such that  $\gamma\paren{h_\lambda^\star, C_{\lambda, {\rm outage}}\paren{\gamma}}$ is equal to $\gamma$. Secondly, we show that $0< \underline{h} = \liminf_{\lambda \ra \infty} h_\lambda^\star \leq \limsup_{\lambda \ra \infty} h_\lambda^\star = \overline{h} < \infty$.  The proof of this assertion is based on proof by contradiction.  For example, if $\limsup_{\lambda \ra \infty} h_\lambda^\star = \infty$, then it follows that $\limsup_{\lambda \ra \infty} \gamma\paren{h, C_{\lambda, {\rm outage}}\paren{\gamma}} = \gamma$ for almost all $h$ with respect to the distribution of $\widetilde{H}$.  By using Theorem \ref{Thm: Rates of Convergence}, this result implies that for any given $\epsilon > 0$, there exist positive constants $B_1$ and $B_2$ such that
\begin{eqnarray}
h_\lambda^\star \leq B_1 \frac{\snr^{-1}+\frac{1}{\pg}\paren{\ES{I_\lambda(1)} + \sqrt{\V{I_\lambda(1)}} \Psi^{-1}\paren{1-\gamma+\epsilon}}}{\snr^{-1}+\frac{1}{\pg}\paren{\ES{I_\lambda(1)} + \sqrt{\V{I_\lambda(1)}} \Psi^{-1}\paren{1-\gamma-\epsilon}}}
\end{eqnarray}
for all $\lambda \geq B_2$. Hence, $\overline{h} \leq B_1 < \infty$, which is a contradiction.  A similar contradiction shows $\underline{h} > 0$.  Finally, it follows that $\log\paren{1+\frac{G(d) \underline{h}}{\snr^{-1} + \frac{1}{\pg}\paren{\ES{I_\lambda(1)} + \sqrt{\V{I_\lambda(1)}} \Psi^{-1}\paren{1-\gamma+\epsilon}}}} \leq C_{\lambda, {\rm outage}}\paren{\gamma} \leq \log\paren{1+\frac{G(d) \overline{h}}{\snr^{-1} + \frac{1}{\pg}\paren{\ES{I_\lambda(1)} + \sqrt{\V{I_\lambda(1)}} \Psi^{-1}\paren{1-\gamma - \epsilon}}}}$ for all $\lambda$ large enough.

{}
     

\begin{thebibliography}{}

\bibitem{ITU08}
International Telecommunication Union, ``Requirements related to technical performance for IMT-Advanced radio interface(s)," {\em ITU-R M.2134 Technical Report}, Available Online: http://www.itu.int/publ/R-REP-M.2134-2008/en, Nov. 2008.

\bibitem{JSAC09-Tutorial} 
M. Haenggi, J. G. Andrews, F. Baccelli, O. Dousse and M. Franceschetti, ``Stochastic geometry and random graphs for the analysis and design of wireless networks," {\em IEEE J. Sel. Areas Commun.}, vol. 27, no. 7, pp. 1029-1046, Sept 2009.

\bibitem{HG09} 
M. Haenggi and R. K. Ganti, ``Interference in large wireless networks," {\em Foundations and Trends in Networking}, vol. 3, no. 2, pp. 127-248, 2009. 

\bibitem{WPS09}
M. Z. Win, P. C. Pinto and L. A. Shepp, ``A mathematical theory of network interference and its applications," {\em Proc. IEEE}, vol. 97, no. 2, pp. 205-230, Feb 2009. 

\bibitem{WA06}
S. P. Weber and J. G. Andrews, ``Bounds on the SIR distribution for a class of channel models in ad hoc networks," in {\em Proc. $49$th IEEE Global Telecommunications Conference}, San Francisco, CA, Dec 2006.

\bibitem{InaltekinHanly10} 
H. Inaltekin and S. V. Hanly, ``On the rates of convergence of the wireless multi-access interference distribution to the normal distribution," in {\em Proc. $6$th Workshop on Spatial Stochastic Models for Wireless Networks}, Avignon, France, June 2010.

\bibitem{AY10}
M. Aljuaid and H. Yanikomeroglu, ``Investigating the Gaussian convergence of the distribution of the aggregate interference power in large wireless networks," {\em IEEE Trans. Veh. Technol.}, vol. 59, no. 9, pp. 4418-4424,  Nov 2010. 

\bibitem{Musa78}
S. Musa and W. Wasylkiwskyj, ``Co-channel interference of spread spectrum systems in a multiple user environment," {\em IEEE Trans. Commun.}, vol. 26, no. 10, pp. 1405-1413, Oct 1978. 

\bibitem{Sousa90}
E. S. Sousa and J. A. Silvester, ``Optimum transmission ranges in a direct sequence spread spectrum multihop packet radio network," {\it IEEE J. Sel. Areas Commun.}, vol. 8, no. 5, pp. 762-771, June 1990.

\bibitem{Sousa92}
E. S. Sousa, ``Performance of a spread spectrum packet radio network in a Poisson field of interferers," \textit{IEEE Trans. Inf. Theory}, vol. 38, no. 6, pp. 1743-1754, Nov 1992. 

\bibitem{IH98}
J. Ilow and D. Hatzinakos, ``Analytic alpha-stable noise modeling in a Poisson field of interferers or scatterers," {\it IEEE Trans. Signal Process.}, vol. 46, no. 6, pp. 1601Ð1611, June 1998.

\bibitem{CH01} 
C. C. Chan and S. V. Hanly, ``Calculating the outage probability in a CDMA network with spatial Poisson traffic," \textit{IEEE Trans. Veh. Technol.}, vol. 50, no. 1, pp. 183-204, Jan 2001.


\bibitem{Gubner96}
J. A. Gubner, ``Computation of shot-noise probability distributions and densities," {\it SIAM Journal of Scientific Computing}, vol. 17, no. 3, pp. 750Ð761, May 1996.

\bibitem{Inaltekin09}
H. Inaltekin, M. Chiang, H. V. Poor and S. B. Wicker, ``On unbounded path-loss models: effects of singularity on wireless network performance," {\em IEEE J. Sel. Areas Commun.}, vol. 27, no. 7, pp. 1078-1092, Sept 2009.

\bibitem{GH09}
R. K. Ganti and M. Haenggi, ``Interference and outage in clustered wireless ad hoc networks," {\it IEEE Trans. Inf. Theory}, vol. 55, pp. 4067-4086, Sept 2009.

\bibitem{Lowen90}
S. B. Lowen. and M. C. Teich, ``Power-law shot noise," \textit{IEEE Trans. on Info. Theory}, vol. IT-36, no. 6, pp. 1302-1318, Nov. 1990.

\bibitem{HS85}
L. Heinrich and V. Schmidt, ``Normal convergence of multidimensional shot noise and rates of this convergence," {\it Advances in Applied Probability}, vol. 17, no. 4, pp. 709-730, Dec. 1985.

\bibitem{Kingman93} 
J. F. C. Kingman, \textit{Poisson Processes}, Clarendon Press, Oxford, 1993. 

\bibitem{Tyurin10}
I. S. Tyurin, ``Refinement of the upper bounds of the constants in Lyapunov's theorem," {\em Communications of the Moscow Mathematical Society}, vol. 65, no. 3, pp. 586-588, 2010.    


\bibitem{Paditz89}
L. Paditz, ``On the analytical structure of the constant in the nonuniform version of the Esseen inequality," {\em Statistics}, vol. 20, no. 3, pp. 453-464, 1989. 



\bibitem{Stuber96}
G. St\"uber, {\it Principles of Mobile Communication}, Kluwer Academic Publishers, Boston, 1996.


\bibitem{BBM06} F. Baccelli, B. Blaszczyszyn and P. Muhlethaler, ``An Aloha protocol for multihop mobile wireless networks," {\em IEEE Trans. Inf. Theory}, vol. 52, no.2, pp. 421-436, Feb 2006.  

\bibitem{WAYV07} 
S. P. Weber, J. G. Andrews, X. Yang and G. Veciana, ``Transmission capacity of wireless ad hoc networks with successive interference cancellation," {\it IEEE Trans. Inf. Theory}, vol. 53, no. 8, pp. 2799-2814, Aug 2007.

\bibitem{Tse05}
D. Tse and P. Viswanath., \textit{Fundamentals of Wireless Communication}, Cambridge University Press, New York, NY, 2005.





















\end{thebibliography}
\end{document}